\documentclass[11pt,a4paper,draft]{article} 
\usepackage{a4wide}
\usepackage{amsmath,amssymb,amsfonts,amsthm,amsxtra,latexsym}

\newtheorem{theorem}{Theorem}

\newtheorem{proposition}{Proposition}
\newtheorem{lemma}{Lemma} 
\newtheorem{claim}{Claim}
\newtheorem{remark}{Remark}
\numberwithin{theorem}{section}
\numberwithin{claim}{section}
\numberwithin{equation}{section}
\numberwithin{lemma}{section}
\numberwithin{proposition}{section}
\DeclareMathOperator{\sech}{sech}
\newcommand{\R}{\mathbb{R}}

\newcommand{\N}{\mathbb{N}}

\newcommand{\pd}{\partial}
\newcommand{\wqs}{\widetilde{Q}_\sigma}
\newcommand{\ys}{y_\sigma}
\begin{document}
\title{Description of the inelastic collision of two solitary waves \\ for the BBM equation
\footnote{
This research was supported in part by the Agence Nationale de la Recherche
(ANR ONDENONLIN).}}
\author{Yvan Martel$^{(1)}$, Frank Merle$^{(2)}$ and
Tetsu Mizumachi $^{(3)}$}
\date{(1) Universit\'e de Versailles Saint-Quentin-en-Yvelines,
 Math\'ematiques, \\
 45, av. des Etats-Unis,
 78035 Versailles cedex, France\\  
 martel@math.uvsq.fr\\ \quad \\ 
(2)
Universit\'e de Cergy-Pontoise, IHES and CNRS,
Math\'ematiques   \\
2, av. Adolphe Chauvin,
95302 Cergy-Pontoise cedex, France \\
 Frank.Merle@math.u-cergy.fr\\ \quad \\
(3)
Kyushu University, Faculty of Mathematics,
\\ Hakozaki 6-10-1, Fukuoka 812-8581, Japan
\\ mizumati@math.kyushu-u.ac.jp}
\maketitle
\begin{abstract}
        We prove that the collision of two solitary waves of the BBM equation is inelastic
        but almost elastic in the case where one solitary wave  is small in the energy space.
        We show precise estimates of the nonzero residue due to the collision.
        Moreover, we give a precise description of the collision phenomenon
        (change of size  of the solitary waves and shifts in their trajectories).
        
        To prove these results, we extend the
        method introduced in \cite{MMcol1} and
        \cite{MMcol2} for the generalized KdV equation, in particular in the quartic case.
        The main argument is the construction of an explicit approximate solution (in a certain sense) in       the collision region. 
\end{abstract}
\section{Introduction}
This paper concerns the Benjamin--Bona--Mahoney equation (BBM), also called
Regularized Long Range equation (RLW)
\begin{equation}\label{eq:BBM}
(1-\pd_x^2)\pd_t u + \pd_x (u + u^2)=0,
\quad (t,x)\in\R\times \R.
\end{equation}
The BBM equation was introduced by Peregrine \cite{Pe}
and Benjamin et al. \cite{BBM} as an alternative to the KdV equation
\begin{equation}\label{kdv}
\pd_t u + \pd_x(\pd^2_{x} u + u^2)=0,
\end{equation}
for the description of  the undirectional propagation of long waves of small amplitude in water.
We recall some aspects of  comparison between these two models (see also \cite{BBM} and \cite{BPS}).

First, recall that the local Cauchy problem in the energy space is much easier
for the BBM equation than for the KdV equation. 
This property was a main motivation to the introduction of the BBM model, in physical situations where it is as relevant as the KdV equation.
Indeed, 
it is easily established that the Cauchy problem for the BBM equation is globally well-posed in $H^1$, see \cite{BBM}. For the KdV equation, the same result is true but it relies on very delicate estimates on the Airy group, see Kenig, Ponce and Vega \cite{KPV}.

Second, whereas it is well-known that the KdV equation is completely integrable (see Section \ref{se:13} for more details), 
the BBM equation is not an integrable model. This lack of the algebraic structure prohibits to extend to the BBM equation the Inverse Scattering Theory developed for the KdV equation and its remarkable applications
(existence of  explicit multi-solitons, decomposition of smooth and decaying solutions, etc.). 
However, when studying the qualitative behavior of solutions  without using the integrability (for example, to study dynamic properties in the energy space), the difficulty of the KdV equation and of BBM equation seem comparable.
For example, the strategy developed in \cite{MM1} for studying the asymptotic stability of solitons for the KdV equation (and its subcritical generalizations, see Section \ref{se:13}), could be adapted (in a nontrivial way) to the BBM equation by Mizumachi \cite{Mi} and El Dika \cite{Di}. See also \cite{DiMa} which extended the results of \cite{MMT}, \cite{Ma2}.

The present paper is another illustration of the similar structures of the generalized KdV equations (gKdV) and the BBM equation, extending the methods of \cite{MMcol1} and \cite{MMcol2} concerning the collision of two solitons of the gKdV equations to the BBM equation. In particular, the results are analogous with the quartic gKdV equation
which is not integrable.

 We point out that although the method of proof is similar to the one of \cite{MMcol1}, the computations are different. The computations are in some sense easier because of the lower order of the nonlinearity, but one has to deal with two parameters, since the equation is not scaling invariant. When the two soliton are small, the equation is close to the KdV model, but since we deal with any possible size for the main soliton, our results are not perturbation of results on the original KdV equation (see Remark \ref{rk:th1}).

\subsection{Solitary waves of the BBM equation}
First, recall that  $H^1$ solutions $u(t)$ of \eqref{eq:BBM} satisfy the following conservation
laws:
\begin{equation}
  \label{eq:cv}
E(u(t))=\int\left(\frac12u^2+\frac13u^3\right)dx=E(u_0),
\qquad N(u(t))= \frac 12 \int (u_x^2+u^2)dx = N(u_0).  
\end{equation}

Next, recall
that equation \eqref{eq:BBM} has a two-parameter family of solitary wave solutions
$ \{\varphi_{c}(x-ct-x_0)\,|\, c>1,x_0\in\R\},$
where $\varphi_{c}$ is satisfies
\begin{equation}
  \label{eq:1}
c\varphi_{c}''-(c-1)\varphi_{c}+\varphi_{c}^2=0\quad\text{on $\R$,}  
\end{equation}
or equivalently
\begin{equation}\label{eq:varphic}
\varphi_{c}(x) = (c-1) Q\left(\sqrt{\frac{c-1}{c}}x\right)
\quad \text{where} \quad 
Q(x)= \frac32\sech^2\left(x\right) \quad \text{solves} \quad Q''+Q^2=Q.
\end{equation}
Such a solution $\varphi_c$ to \eqref{eq:1} is a critical point of 
$S_c(u)=cN(u)-E(u)$. 
This variational characterization allows one to apply the Lyapunov
stability theory (see Grillakis et al. \cite{GSS}, Weinstein \cite{We3} and references therein).
Since $dN(\varphi_c)/dc>0$ for
any $c>1$, all the solitary wave $\varphi_c$  are orbitally
stable in the following sense.

\medskip

\noindent\textbf{Stability of solitary waves (\cite{We3})}
\textit{For any
$\varepsilon>0$, there exists  $\delta>0$ s.t., for all $c>1$, $u_0\in H^1$,}
\begin{equation}\label{eq:stab}
\inf_{y\in\R}\| u_0-\varphi_c(\cdot-y)\|_{H^1}<\delta
\quad \Rightarrow \quad \sup_{t\ge0}\inf_{y\in\R}
\|u(t,\cdot)-\varphi_c(\cdot-y)\|_{H^1}<\varepsilon.
\end{equation}

Next, we recall  that the family of
solitary waves is asymptotically stable in $H^1$.

\medskip

\noindent\textbf{Asymptotic stability (\cite{Mi2, Di, Ma})}
\textit{There exists  $\delta>0$ s.t., for all $c>1$, $u_0\in H^1$,}
\begin{equation}\label{eq:asym}
\inf_{y\in\R}\| u_0-\varphi_c(\cdot-y)\|_{H^1}<\delta
\quad \Rightarrow \quad 
\lim_{t\to +\infty}\|u(t,\cdot)-\varphi_{c^+}(\cdot-\rho(t))\|_{H^1(x>\frac 12 (1+c) t)}= 0,
\end{equation}
\emph{for some $c^+$ close to $c$ and some function $\rho(t)$ such that $\lim_{t\to +\infty}\rho'(t)= c^+$.}

\medskip

We refer to the main results of \cite{Mi2}, \cite{Di} and to Theorem 2 in \cite{Ma}.
See Miller and Weinstein \cite{MW} for previous results concerning the asymptotic stability of solitary waves for BBM in weighted spaces, and under some restriction on $c>1$.

Moreover, recall that  similar results of stability and asymptotic stability in $H^1$ hold for solutions close to the sum of several decoupled solitary waves  (see \cite{DiMa}, \cite{Ma} and Section \ref{sec:3} of the present paper).
Finally, in \cite{DiMa}, the following existence and uniqueness result is proved.

\medskip

\noindent\textbf{Asymptotic multi-solitary waves (\cite{DiMa})} \emph{Let $N\geq 1$,
$1<c_N<\ldots<c_1$ and $x_1,\ldots,x_N\in \R$. There exists a unique $H^1$ solution
$u(t)$ of \eqref{eq:BBM} such that}
\begin{equation}\label{eq:pure}
\lim_{t\to -\infty} \bigg\| u(t,\cdot)-\sum_{j=1}^N \varphi_{c_j} (\cdot-x_j-c_j t)\bigg\|_{H^1} =0.
\end{equation}

\medskip

Such solutions behave asymptotically as $t\to +\infty$ as the sum of $N$ solitons with different speeds. By the symmetry $x\to -x$, $t\to -t$ of the BBM equation, there exist solutions with similar behavior as $t\to +\infty$. But the global behavior of such solutions is unknown, i.e. during and after the collision of the various solitons.

\subsection{Statement of the problem and numerical predictions}
In this paper, we consider the problem of the collision of two solitary waves for 
\eqref{eq:BBM}.
We focus on the following more precise questions: let $1<c_2<c_1$ and 
let $u(t)$ be the $H^1$ solution of \eqref{eq:BBM} such that 
\begin{equation}\label{eq:2pure}
\lim_{t\to -\infty} \bigg\| u(t,\cdot)-\sum_{j=1,2} \varphi_{c_j} (\cdot-c_j t)\bigg\|_{H^1} =0.
\end{equation}
What is the behavior of $u(t)$ during and after the collision? In particular,
do we recover two solitary waves after the collision? In this case, how the two solitary waves are changed by the collision? Is the collision elastic (zero residue) or inelastic (small but nonzero residue)? 

There has been a lot of numerical work on the  BBM equation, which has become a kind of test  
problem for numerical schemes for nonlinear wave equations.
Numerical simulations generally predict that the two solitons are preserved by the collision, but a small residue appears (inelastic collision). 
See Abdulloev et al. \cite{ABM}, Eilbeck and McGuire \cite{E}, Bona et al. \cite{BPS}, Kalisch and Bona \cite{KB}. This is in contrast with the integrable case, for which the collision is completely elastic.

From the experiment point of view, or from the numerical point of view for more elaborate systems, it also seems that the inelastic but close to elastic collision is the most established conjecture, see Craig et al. \cite{CGHHS} and  Hammack \cite{HHGY}.

\subsection{Collision problem for the generalized KdV equations}\label{se:13}

In the integrable case, i.e. for the Korteweg de Vries equation (KdV)
(which corresponds in some sense to the limit $c\to 1$ in \eqref{eq:varphic}, see Remark \ref{re:np}),
the collision of several solitons of different velocities is described  by the  explicit pure multi-soliton solutions
(see \cite{HIROTA}, \cite{WT} as well as \cite{Miura} and references therein). It is well-known that
any two solitons with different velocities collide elastically, with no size change, and only suffer an explicit shift on their trajectory due to the collision.

Recall that 
Fermi, Pasta and Ulam \cite{FPU} and Zabusky and Kruskal  \cite{KZ} presented the first numerical experiments related to soliton collision. Then, Lax \cite{LAX1} introduced  a mathematical framework to study these problems,
known now as complete integrability. Other developments
appeared, such as the Inverse Scattering Transform (see for example the review paper by  Miura \cite{Miura}).
This nonlinear transformation led to one of the most striking property of the KdV and mKdV equations  which is the existence of explicit
$N$-soliton solutions (Hirota \cite{HIROTA}, Wadati and Toda, \cite{WT}). Moreover, one consequence of the inverse scattering transform is the so-called decomposition result (Kruskal \cite{KRUSKAL}, Eckhaus and Schuur \cite{EcSc}, \cite{Schuur},
Cohen \cite{Cohen}): any smooth and decaying solution of \eqref{kdv} eventually decomposes as $t\to +\infty$ as the sum of a finite number of solitons.

\medskip

Consider now the generalized KdV equations
\begin{equation}\label{gkdv}
\partial_t u + \partial_x (\partial_x^2 u + f(u))=0.
\end{equation}
The integrable cases correspond to $f(u)=u^2$ and $f(u)=u^3$.
For other nonlinearities, the Inverse Scattering Transform is not applicable and no explicit multi-soliton are known.

The question of collision of two solitary waves for \eqref{gkdv} was
addressed recently from another point of view by Martel and Merle (\cite{MMcol1,MMcol2}).
Recall  that solitary waves of \eqref{gkdv}, also called solitons, are solutions of the form $u(t,x)=Q_c(x-ct)$ where $c>0$, and $Q_c$ solves
$Q_c''+f(Q_c)=cQ_c$.

Under general conditions of $f$, there exists $c_*>0$ (note that $c_*=+\infty$ for
the subcritical power case, i.e. $f(u)=u^2$, $u^3$, $u^4$) such that for all $0<c<c^*$, the soliton $Q_c(x-ct)$ is stable and asymptotically stable. In \cite{MMcol1}, \cite{MMcol2}, the collision of two solitons $Q_{c_1}(x-c_1t)$, $Q_{c_2}(x-c_2t)$ such that $0<c_2<c_1<c^*$ is studied under the assumption 
$$
0<c_2\ll 1,
$$
which means that one solitary wave is small in the energy space.
Under this condition, for a general nonlinearity $f$, it was proved that the collision of two stable solitary waves
of \eqref{gkdv} is elastic up to a possible residual, small compared 
to the solitary waves (see Theorems 1, 2 and 3 in \cite{MMcol2}). 
Moreover,  monotonicity properties are obtained: the size of the large soliton does not decrease and the size of the small soliton does not increase through the collision.

In the quartic case, i.e. $f(u)=u^4$ in \eqref{gkdv}, the description of the collision was  much refined.
Indeed, in \cite{MMcol1}, it was proved in that case that the residue is nonzero, with a precise estimate on its size,
which implies that the collision is inelastic, but close to elastic.
As a consequence,  the  monotonicity properties are strict: the size of the large soliton increases and the size of the small soliton decreases through the collision, with explicit lower and upper bounds on the discrepancies.
Moreover,  the first orders of the shifts resulting from the collision
could be computed explicitly.

\subsection{Main results}

In the present paper, we extend the method introduced in \cite{MMcol1}, \cite{MMcol2} to the BBM equation. We consider two solitary waves (hereafter called solitons) $\varphi_{c_1}(x-c_1t)$, $\varphi_{c_2}(x-c_2t)$
in the case where $1<c_2<c_1$ and $c_2$ is close to $1$, so that by \eqref{eq:varphic}
the function $\varphi_{c_2}$ is small in $H^1$. In this context, we prove that the collision of
the two solitons is not elastic but almost elastic.

The main result of this paper is the following theorem.

\begin{theorem}\label{th:1}
Let $c_1>c_2>1$ and let $u(t)$ be the unique solution of \eqref{eq:BBM} such that
\begin{equation}\label{eq:th:1:0}
\lim_{t\to -\infty} \bigg\| u(t,\cdot)-\sum_{j=1,2} \varphi_{c_j} (\cdot-c_j t)\bigg\|_{H^1} =0.
\end{equation}
There exists $\epsilon_0=\epsilon_0(c_1)>0$ such that if $0<c_2-1<\epsilon_0$,
then there exist $c_1^+>c_2^+>1$, $\rho_1(t)$, $\rho_2(t)$ and    $T_0,K>0$ such that
$$
w^+(t,x)=u(t,x)-\sum_{j=1,2} \varphi_{c_j^+} (x-\rho_j(t))
$$
satisfies
\begin{equation}\label{eq:th:1:1}
\lim_{t\to +\infty} \|w^+(t)\|_{H^1(x> \frac 12 (1+c_2) t)} =0,
\end{equation}
and
\begin{equation}\label{eq:th:1:2}
\tfrac 1K (c_2-1)^{\frac {11}2} \leq  {c_1^+} - {c_1}  \leq K (c_2-1)^{\frac 92},\qquad
\tfrac 1K (c_2-1)^{5} \leq c_2- {c_2^+} \leq K (c_2-1) ^{4},
\end{equation}
\begin{equation}\label{eq:th:1:3}
\tfrac 1K (c_2-1)^{\frac {11}4} \leq \|\partial_x w^+(t)\|_{L^2}+\sqrt{c_2-1}\|w^+(t)\|_{L^2}\leq K (c_2-1)^{\frac  94}, \qquad \text{for $t\geq T_0$}.
\end{equation}
\end{theorem}

\begin{remark}\label{re:ze}
Theorem \ref{th:1} implies that  there exists no
pure $2$-soliton solution corresponding to the speeds $c_1$, $c_2$ in this regime. 
Indeed, by \eqref{eq:th:1:1} and the lower bound in \eqref{eq:th:1:3}, the perturbative term $w^+(t)$ does not go to zero in the region $x<\frac 12 (1+c_2) t$.

Thus, the conclusion of Theorem \ref{th:1} matches the numerical predictions mentioned in Section 1.2.

In a different spirit, let us mention that Bryan and Stuart \cite{BS} proved the nonexistence of a family  of multi-solitons of the BBM equation which would be an analytic continuation of the multi-solitons of the KdV equation.
\end{remark}

\begin{remark}
It is an open problem to understand the exact asymptotic behavior of $w^+(t)$ as $t\to +\infty$. Recall for the quartic gKdV equation (i.e. $f(u)=u^4$ in \eqref{gkdv}), Tao's paper \cite{Tao2} implies in a similar situation that $w^+(t)$ is purely dispersive.
For the BBM equation, it is probably not the case.
\end{remark}

\begin{remark}\label{rk:th1}
For $c_1>1$ small, the BBM equation is close to the KdV equation, see Remark~\ref{re:np}. But in Theorem \ref{th:1}, we allow any value of $c_1>1$, which means that the results are not perturbative of the KdV case.
\end{remark}

\begin{remark}\label{re:zu}
The size of the perturbation $w^+(t)$ is controlled by \eqref{eq:th:1:3}.
This is to compare with $\|\partial_x \varphi_{c_2}\|_{L^2}+\sqrt{c_2-1}\|\varphi_{c_2}\|_{L^2}
\sim K (c_2-1)^{\frac 54}.$

In some sense, the BBM problem is less degenerate than the quartic gKdV which is closer to the critical case (i.e. $f(u)=u^5$ in \eqref{gkdv}, for which all solitons have the same size in $L^2$).

Note that 
the estimate \eqref{eq:th:1:3} is not sharp, see Remark \ref{rk:no}. A similar gap in the estimates is observed in Theorem 1.1 of \cite{MMcol1}.

We point out that 
as a consequence of the our analysis of the collision for the BBM equation, one can obtain results similar to Theorems 1.2 and 1.3 in \cite{MMcol1}. In particular, thanks to Lemma \ref{lem:z}, on can construct a symmetric solution of \eqref{eq:BBM} (i.e. verifying $v(-t,-x)=v(t,x)$) with a sharp estimate of the perturbation.
Moreover, it is quite clear from the stability results of Section 3 that the behavior of the solution $u(t)$ of Theorem \ref{th:1} is stable with respect to $H^1$ perturbations (see Theorem 1.3 in \cite{MMcol1}).

See \eqref{eq:Dt} for the first order of the shifts  on the solitons after and before the collision region.
\end{remark}

The main argument of the proof of Theorem \ref{th:1} is the construction of an approximate solution of the BBM equation which describes the collision of two solitons of speeds $1<c_2<c_1$
in a large but finite time interval $[-T,T]$ containing the collision, similarly as in \cite{MMcol1}, Section~2.
See Section~2 of the present paper.

Then, large time stability arguments are used to compare this approximate solution to the solution $u(t)$ on $[-T,T]$, see Section 3.1.
Then, for $|t|>T$, the solitons are decoupled and 
large time asymptotic arguments are used to describe the asymptotic behavior of the solution
as $t\to +\infty$. These arguments are refinements of the ones of \cite{Mi2}, \cite{Di}, \cite{DiMa}, \cite{Ma}, see also \cite{MW} for a previous work on asymptotic stability. See Section 3.2.

The proof of Theorem \ref{th:1} is given in Section 4.

Appendices A, B, C and D are devoted to the proof of some technical results.

\medskip

\noindent\textbf{Acknowledgements.} We would like to thank Professor J.~C.\ Eilbeck
for indicating us several references.

\section{Construction of an approximate 2-soliton solution}\label{sec:2}

The objective of this section is to construct an approximate solution of
the BBM equation, which describes the collision of two solitons $\varphi_{c_1}$, $\varphi_{c_2}$
in the case where $1<c_2-1<\epsilon_0$ is small.
The main result of this section is Proposition \ref{cor:1} in Section \ref{se:26}.

The interest of changing the variable is to reduce ourselves to some simple algebra in the function $Q$ (see \eqref{eq:varphic}) similar to the one in \cite{MMcol1} for the KdV equation
($u^2$ nonlinearity). 

\subsection{Reduction of the problem}
Let 
\begin{equation}\label{eq:dl}
c_1>1 \quad \text{and}\quad \lambda=\frac{c_1-1}{c_1}\in(0,1).
\end{equation}
We introduce the following change of variable
\begin{equation}
  \label{eq:chv}
x'=\lambda^{1/2}\left(x-\frac{t}{1-\lambda}\right),
\quad t'=\frac{\lambda^{3/2}}{1-\lambda}t,\quad
 z(t',x')=\frac{1-\lambda}{\lambda}u(t,x).      
\end{equation}
If $u(t,x)$ is a solution to \eqref{eq:BBM} then $z(t',x')$ satisfies
\begin{equation}\label{eq:BBM0}
(1-\lambda \pd_{x'}^2) \pd_{t'} z + \pd_{x'} (\pd_{x'}^2 z - z +z^2)=0.
\end{equation}
\begin{remark}\label{re:np}
Observe that when $c_1\to 1$ so that  $\lambda\to 0$, the above equation 
converges to the KdV equation. For $0<\lambda<1$, the problem is not perturbative of the KdV case.
\end{remark}

\begin{claim}\label{cl:2.1}
  \begin{itemize}
  \item [\rm (i)]
Let $c>1$. By the change of variable \eqref{eq:chv},
a solitary wave solution $\varphi_c(x-ct)$ to \eqref{eq:BBM}
-- see \eqref{eq:varphic} -- is transformed into $\wqs(\ys),$ a solution of
\eqref{eq:BBM0} where
\begin{align*}
& 
\wqs(x):=\sigma \theta_\sigma Q(\sqrt{\sigma}x),\quad Q(x)=\frac{3}{2}\sech^2(x/2),\\
& \sigma=\frac{c-1}{c\lambda},\quad
\theta_\sigma=\frac{1-\lambda}{1-\lambda\sigma},\quad
\mu_\sigma=\frac{1-\sigma}{1-\lambda\sigma},\quad
\ys=x'+\mu_\sigma t'.
\end{align*}
Especially if $c=c_1$, then $\mu_\sigma=0$, $\ys=x'$ and $\wqs(\ys)=Q(x')$
and
\begin{equation}
\label{eq:Q}
 Q''-Q+Q^2=0,\quad
 (Q')^2 + \frac 23 Q^3 = Q^2 \quad\text{on $\R$}.
 \end{equation}
\item [\rm (ii)] Moreover, $\wqs$ satisfies the following
  \begin{equation}
    \label{eq:qsigma}
\wqs''=\sigma\wqs-\frac{1}{\theta_\sigma}\wqs^2,\quad
(\wqs')^2=\sigma\wqs^2-
\frac{2}{3\theta_\sigma}\wqs^3.
  \end{equation}
  \begin{equation}
  \label{eq:qsigmanorm}
\|\wqs\|_{L^\infty}\sim {(1-\lambda)}\sigma \|Q\|_{L^\infty},
\quad \|\wqs\|_{L^2}\sim  (1-\lambda)\sigma^{3/4} \|Q\|_{L^2},
\quad \text{for $\sigma>0$ small}.
\end{equation}
  \end{itemize}
\end{claim}

\begin{proof}
Note  that $\varphi_c(x)=(c-1)Q\left(\sqrt{\frac {c-1}c} x\right)$.
First, we have, for any $x'$, $t'$, $\delta$,
$$
\wqs(x'+\mu_\sigma t' + \delta)
=\sigma \theta_\sigma Q(\sqrt{\sigma \lambda} ( x- \tfrac {1}{1-\lambda} t + \tfrac \lambda {1-\lambda}  \mu_\sigma t + \tfrac \delta {\sqrt{\lambda}} )).
$$
But by direct computations, $1-\lambda \sigma = \frac 1c$ so that
$\sigma \theta_\sigma = \frac {1-\lambda}\lambda (c-1)$ and
$\frac 1{1-\lambda} - \mu_\sigma \frac \lambda {1-\lambda} = c$.
Thus,
\begin{equation}\label{eq:v1}
\frac \lambda{1-\lambda} \wqs(\ys+\delta)=\varphi_c(x-ct+\tfrac 1 {\sqrt{\lambda}} \delta).
\end{equation}
For $c=c_1$, we have $\sigma=1$, $\theta_\sigma=1$, $\mu_\sigma =0$.

The equation of $\wqs$ is quite clear from the equation of $Q$.
Estimates \eqref{eq:qsigmanorm} are also straightforward.

Finally, for future reference, we compute $(\wqs^2)'(x'+\mu_\sigma t' + \delta)$ in terms of 
$\varphi_c$.  As before, we have
\begin{align}
&(\wqs^2)'(\ys + \delta)
  =\sigma^2  \theta_\sigma^2 \sigma^{\frac 12} (Q^2)'(\sqrt{\sigma \lambda} ( x- c t  + \tfrac \delta {\sqrt{\lambda}} ))\label{eq:v2}\\
& = \left(\frac {1-\lambda}{\lambda}\right)^2 (c-1)^2 \sqrt{\frac {c-1}c} \frac 1 {\sqrt{\lambda}} (Q^2)'(\sqrt{\sigma \lambda} ( x- c t + \tfrac \delta {\sqrt{\lambda}} ))
=\frac {(1-\lambda)^2}{\lambda^{\frac 52}} (\varphi_c^2)' (x-ct + \tfrac \delta {\sqrt{\lambda}} ).\nonumber
\end{align}
\end{proof}

\subsection{Decomposition of the approximate solution}\label{sec:2-1}
We construct an approximate   solution
$z(t,x)$ of 
\begin{equation}\label{eq:BBM2}
(1-\lambda \pd_x^2) \pd_t z + \pd_x (\pd_x^2 z - z +z^2)=0,
\end{equation}
 which is a superposition of the function $Q$,
a small soliton $\widetilde{Q}_{\sigma}$ and an error term $w(x,t)$.
As in \cite{MMcol1}, we introduce the new coordinates and the approximate solution under the following form
\begin{align}
& \ys=x+\mu_\sigma t,\quad y=x-\alpha(\ys),\quad \alpha(s)=\int_0^s\beta(r)dr,\quad
\beta(s)=\sum_{(k,l)\in\Sigma_0}a_{k,l}\sigma^l\wqs^k(s),\label{eq:df1}\\
& z(t,x)=Q(y)+\wqs(\ys)+w(t,x), \label{eq:df2}
\end{align}
where
\begin{align}
& w(t,x)=\sum_{(k,l)\in\Sigma_0}\sigma^l\left\{
A_{k,l}(y)\wqs^k(\ys)+B_{k,l}(y)(\wqs^k)'(\ys)\right\},
\\& \Sigma_0:=\{(k,l)\in \N^*\times \N\,:\, k+l\le 3\},\label{eq:df3}
\end{align}
and $a_{k,l}$,
$A_{k,l}$, $b_{k,l}$ are to be determined. Let
\begin{equation}\label{eq:df4}
  S(z)=(1-\lambda\pd_x^2) \pd_t z+\pd_x(\pd_x^2z-z+z^2).
\end{equation}
Define the operator $L$ by  
$$
LA=-A''+A-2QA
$$
(see Lemma \ref{surL} for properties of $L$).
We claim
\begin{proposition}
\label{prop:decomp}
The following holds
  \begin{align*}
 S(z)   & =\sum_{(k,l)\in\Sigma_0}\sigma^l\wqs^k(\ys)\big(
a_{k,l}((\lambda-3)Q''-Q^2)' -(LA_{k,l})' +F_{k,l} \big)(y)
\\ &   +\sum_{(k,l)\in \Sigma_0}\sigma^l(\wqs^k)'(\ys) \big(
(3-\lambda)A_{k,l}''+2QA_{k,l}+a_{k,l}(2\lambda-3)Q''-(LB_{k,l})' +G_{k,l}\big)(y)
\\ &   +\mathcal{E}(t,x),
  \end{align*}
where 
\begin{align*}
F_{1,0}&=2Q',\quad G_{1,0}=2Q,\\
F_{1,1}&= (3-2\lambda)A_{1,0}'+(3-\lambda)B_{1,0}''+2QB_{1,0}
+\lambda(\lambda-1)a_{1,0}Q''',\\
G_{1,1}&=  \lambda(1-\lambda)A_{1,0}''+(3-2\lambda)B_{1,0}'+
2a_{1,0}\lambda(1-\lambda)Q'',\\
 F_{2,0}&=a_{1,0}\{(\lambda-3)A_{1,0}''-2QA_{1,0}-Q\}'
+(3-2\lambda)a_{1,0}^2Q'''\\ &
+(A_{1,0}^2)'-\frac{2}{1-\lambda}QB_{1,0}-\frac{1}{1-\lambda}A_{1,0}'
+\frac{\lambda-3}{1-\lambda}B_{1,0}'',\\
 G_{2,0}&= \frac{a_{1,0}}{2}\left\{(6\lambda-9)A_{1,0}'
+(\lambda-3)B_{1,0}''-2QB_{1,0}\right\}' \\ &
+A_{1,0}^2+(A_{1,0}B_{1,0})'+A_{1,0}+\frac{\lambda-2}{1-\lambda}B_{1,0}'
+\frac{3}{2}(1-\lambda)a_{1,0}^2Q'',
\end{align*}
where the following holds:

(i) For all $(k,l)\in \Sigma_0$ such that $k+l=3$,  $F_{k,l}$, $G_{k,l}$ depend on 
$A_{k',l'}$, $B_{k',l'}$ for $1\leq k'+l'\leq 2$.
Moreover, if $A_{k',l'}$ are even and $B_{k',l'}$ are odd then 
$F_{k,l}$ are odd and $G_{k,l}$ are even.

(ii) If the functions $A_{k,l}$, $B_{k,l}$ are bounded then the rest term 
 $\mathcal{E}(t,x)$ satisfies
\begin{equation}\label{eq:E}
|\mathcal{E}(t,x)|\leq K \sigma^3 \wqs(\ys).
\end{equation}
\end{proposition}
Proposition \ref{prop:decomp} is a consequence of
 Lemmas \ref{lem:SQ}--\ref{lem:Sint} in Appendix B.

\begin{remark}
Note that to prove the main results of this paper,  the explicit expressions of $F_{k,l}$ and $G_{k,l}$
for the cases $(k,l)=(1,1)$ and $k+l=3$ are not needed. Note also that for $\lambda =0$ we recover the KdV case,
see Proposition 2.1 with $p=2$ in \cite{MMcol1}.
\end{remark}

Proposition \ref{prop:decomp} means that if the system
\begin{equation}
  \label{eq:system1}
  (\Omega_{k,l})\quad 
  \left\{
    \begin{aligned}
& (LA_{k,l})'=a_{k,l}\{(\lambda-3)Q''-Q^2\}'+F_{k,l},\\
& (LB_{k,l})'=(3-\lambda)A_{k,l}''+2QA_{k,l}+a_{k,l}(2\lambda-3)Q''+G_{k,l}
    \end{aligned}
\right.
\end{equation}
is solved for every $(k,l)\in\Sigma_0$, 
then $S(z)=\mathcal{E}$ is small.

\subsection{Explicit resolution of the systems $(\Omega_{1,0})$ and $(\Omega_{2,0})$}
For $(\Omega_{1,0})$ and $(\Omega_{2,0})$, we look for explicit solutions such that
\begin{equation}\label{eq:st}
A_{k,l}=\widetilde A_{k,l} + \gamma_{k,l},\quad
B_{k,l}=\widetilde B_{k,l} + b_{k,l}\varphi,
\end{equation}
where $\widetilde A_{k,l}\in \mathcal{Y}$ is even and 
$\widetilde B_{k,l}\in \mathcal{Y}$ is odd
(see Section \ref{sec:32} for a justification of this choice).

First, we recall some preliminary notation and results from \cite{MMcol1}.
We denote by $\mathcal{Y}$ the set of $C^\infty$ functions $f$ such that
\begin{equation}\label{eq:Y}
\forall j\in \mathbb{N},~\exists K_j,~r_j>0,~Ê\forall x \in \mathbb{R},\quad
|f^{(j)}(x)|\leq K_j (1+|x|)^{r_j} e^{-|x|}.
\end{equation}

\begin{lemma}[Lemma 2.2 from \cite{MMcol1}]\label{surL}  The operator ${L}$ defined in $L^2(\mathbb{R})$ by
    $
        {L} f= -f''+f-2Q f
    $
    is self-adjoint and satisfies the following properties:
    \begin{itemize}
        \item[{\rm (i)}] First eigenfunction : ${L} Q^{\frac 3 2} = - \frac 54  Q^{\frac 3 2}$;
        \item[{\rm (ii)}] Second eigenfunction : ${L} Q'=0$; the kernel of ${L}$ is 
        $\{\lambda Q', \lambda \in \mathbb{R}\}$;
        \item[{\rm (iii)}] For any   function $h \in L^2(\mathbb{R})$ orthogonal to $Q'$ for the $L^2$ scalar product, 
        there exists a unique function $f \in H^2(\mathbb{R})$ orthogonal to $Q'$ such that ${L} f=h$; moreover,
        if $h$ is even (respectively, odd), then $f$ is even (respectively, odd).
        \item[{\rm (iv)}] Suppose that $f\in H^2(\mathbb{R})$ is such that ${L} f \in \mathcal{Y}$.
            Then, $f\in \mathcal{Y}$.
    \end{itemize}
\end{lemma}

\begin{claim}[Claim 2.1 in \cite{MMcol1}]
\label{cl:phi} Let $\varphi=-\frac{Q'}{Q}$. Then
  \begin{itemize}
  \item [\rm (a)] $\lim_{x\to\pm\infty}\varphi(x)=\pm 1$,
  \item [\rm (b)] $\varphi'\in \mathcal{Y}$,
   $1-\varphi^2 \in \mathcal{Y}$,
  \item [\rm (c)] $(L\varphi)'=2Q-\frac{5}{3}Q^2=\frac{1}{3}Q+\frac{5}{3}Q'',$
   $\varphi'=\frac 13 Q$,
  \item [\rm (d)] $\int \varphi Q'=-\frac13\int Q^2$.
  \end{itemize}
\end{claim}

We continue with the general computation of $a_{k,l}$, 
assuming that \eqref{eq:st} holds. Let
\begin{equation}\label{eq:defP}\begin{split}
& P_{\lambda}=2 Q + \left(\frac {3-\lambda} 2\right) xQ',\quad
 P=P_{1}=2 Q + xQ',\\
& \text{so that} \quad L P= - 2Q, \quad L P_{\lambda}= - ((3-\lambda) Q'' + 2Q^2)
\end{split}\end{equation}
\begin{claim}\label{lem:deta}
Assume \eqref{eq:st} and $(a_{k,l}, A_{k,l}, B_{k,l})$ solves
$(\Omega_{k,l})$. Then,
\begin{equation}\label{eq:akl}\begin{split}
a_{k,l}  = - \frac {20} {15+10\lambda -\lambda^2 }  \frac 1{\int Q^2}
\left\{ - \gamma_{k,l} \int P_\lambda +\int G_{k,l} Q
+ \int F_{k,l} \int_0^{x} P_\lambda\right\}. 
\end{split}\end{equation}
\end{claim}
\begin{proof}
First, we multiply the equation of $B_{k,l}$ by $Q$ and use $L Q'=0$. We get
\begin{align*}
a_{k,l} (2\lambda -3) \int (Q')^2 & =  
\int ((3-\lambda) Q'' + 2 Q^2) A_{k,l} + \int G_{k,l} Q\\
&= - \int (LA_{k,l}) P_\lambda + \int G_{k,l} Q.
\end{align*}
Second, we multiply the equation of $A_{k,l}$ by $\int_0^x P_\lambda(y) dy$.
We obtain
\begin{align*}
\int (L A_{k,l})' \int_0^x P_\lambda &= -\int (L A_{k,l}) P_\lambda 
+\gamma_{k,l} \int P_\lambda\\
& = - a_{k,l} \int ((\lambda -3) Q'' -Q^2) P_\lambda + \int F_{k,l} \int_0^x P_\lambda.
\end{align*}
Thus, combining the two identities, we get:
\begin{align*}
&a_{k,l} \left\{
(2\lambda - 3) \int(Q')^2 + \int ((\lambda-3) Q''- Q^2) P_\lambda\right\}\\
& = - \gamma_{k,l} \int P_\lambda + \int G_{k,l} Q + \int F_{k,l} \int_0^x P_\lambda. 
\end{align*}
Now, we need only observe that using Claim \ref{cl:qint}, we get
\begin{equation}\label{eq:on}
(2\lambda - 3) \int(Q')^2 + \int ((\lambda-3) Q''- Q^2) P_\lambda
=- \frac {15+ 10 \lambda - \lambda^2}{20} \int Q^2 
\end{equation}
\end{proof}

\begin{lemma}[Resolution of $(\Omega_{1,0})$]\label{lem:omega10}
The following is solution of $(\Omega_{1,0})$:
\begin{equation}\label{eq:solOmega10}\begin{split}
 &  a_{1,0}=\frac{10(1+\lambda)}{15+10\lambda-\lambda^2},\quad 
A_{1,0}=-(yQ'+2Q)-a_{1,0}\left(\frac{\lambda-3}{2}yQ' -Q\right),\\
& B_{1,0}=\frac{3-\lambda}{4}y^2Q'+yQ-a_{1,0}\left\{
\frac{(\lambda-3)^2}{8}y^2Q'+\frac{3-\lambda}{2}yQ\right\}+b_{1,0}\varphi+\kappa Q',
\\ &
b_{1,0}=\lim_{+\infty} B_{1,0}=3\left(\frac{\lambda+1}{2}a_{1,0}-1\right)
=\frac{-30+18\lambda^2}{15+10\lambda-\lambda^2}.
\end{split}\end{equation}
\end{lemma}
\begin{remark}
We obtain as required $A_{1,0}\in \mathcal{Y}$ (i.e. $\gamma_{1,0}=0$) 
and $B_{1,0}$ as in \eqref{eq:st}.
The constant $\kappa$ in the expression of $B_{1,0}$ above is a free parameter
that we will fix such that $\int B_{1,0} Q'=0$ for convenience in some future
computations. 
By \eqref{eq:solOmega10}, we have
\begin{align*}
\int Q'B_{1,0}&= \int\left(\frac{3-\lambda}4y^2(Q')^2-\frac12Q^2\right)
-a_{1,0} \left\{\frac{(\lambda-3)^2}8\int x^2(Q')^2-\frac{3-\lambda}{4}
\int Q^2\right\}
\\ & +b_{1,0}\int\varphi Q'+\kappa\int Q'^2.
\end{align*}
Thus by \eqref{eq:Q} and Claims \ref{cl:qint} and \ref{cl:phi},
\begin{align*}
\kappa=\frac{5}{3}b_{1,0}+5\left(\frac{\lambda-3}{4}a_{1,0}+\frac12\right)
+\frac{\pi^2}{24}\left\{(\lambda-3)^2a_{1,0}+2(\lambda-3)\right\}.  
\end{align*}
\end{remark}
\begin{proof}[Proof of Lemma \ref{lem:omega10}]
First, we determine $a_{1,0}$. We look for $A_{1,0}$ in $L^2$, i.e.
$\gamma_{1,0}=0$ in the assumption \eqref{eq:st}.
Recall that from Proposition  \ref{prop:decomp}, 
$F_{1,1}=2Q'$, $G_{1,1}=2Q$.
Thus, from Claim \ref{lem:deta}, we obtain:
\begin{equation}
  \label{eq:a10}
a_{1,0}= - \frac{20}{15+10\lambda-\lambda^2}{\frac 1{\int Q^2}}
\left( 2 \int Q^2 -2 \int Q P_\lambda \right)
=\frac{10(1+\lambda)}{15+10\lambda-\lambda^2}.
\end{equation}

We look for $A_{1,0}\in \mathcal{Y}$, even.
By integration, \eqref{eq:system1} with $(k,l)=(1,0)$ is translated into the system
  \begin{align}
\label{eq:A10}
& LA_{1,0}=2Q+a_{1,0}\{(\lambda-3)Q''-Q^2\}, \\
\label{eq:B10}
& (LB_{1,0})'=(3-\lambda)A_{1,0}''+2QA_{1,0}+a_{1,0}(2\lambda-3)Q''+2Q.
  \end{align}
Since $LQ=-Q^2$ and $L(yQ')=-2Q''$, 
\begin{equation}
  \label{eq:sola10}
A_{1,0}=-(yQ'+2Q)-a_{1,0}V_\lambda,\quad\text{where}\quad
V_\lambda=\frac{\lambda-3}{2}yQ' -Q,
\quad 
LV_\lambda = (3-\lambda)Q''+Q^2,
\end{equation}
is an even solution to \eqref{eq:A10} in $\mathcal{Y}$.

Next, we find $B_{1,0}$.
Assuming that  $B_{1,0}$ is smooth and bounded, it follows from \eqref{eq:B10} 
and $LQ'=0$ that $((LB_{1,0})',Q)_{L^2} =0$.
By \eqref{eq:B10}, \eqref{eq:sola10} and \eqref{eq:Q},
\begin{equation}
\label{eq:B10'}
\begin{split}
 (LB_{1,0})'
 &= (\lambda-3)(yQ'+2Q)''-2(yQQ'+2Q^2)+2Q
\\ & +a_{1,0}\{(\lambda-3)V_\lambda''-2QV_\lambda+(2\lambda-3)Q''\}
\\&= \{(\lambda-3)yQ''+(3\lambda-6)Q'-yQ^2\}'-Q
\\ &+ a_{1,0}\left[\left\{\frac{(\lambda-3)^2}{2}yQ''
+\frac{\lambda^2-5\lambda+8}{2}Q'
+\frac{3-\lambda}{2}yQ^2\right\}'+\frac{\lambda+1}{2}Q\right].
\end{split}
\end{equation}
Set
\begin{equation}
  \label{eq:b10}
b_{1,0}=3\left(\frac{\lambda+1}{2}a_{1,0}-1\right)
=\frac{-30+18\lambda^2}{15+10\lambda-\lambda^2},\quad
\widetilde B_{1,0}= B_{1,0} - b_{1,0} \varphi,
\quad \varphi=-\frac{Q'}{Q}.
\end{equation}
Then by Claim \ref{cl:phi} and \eqref{eq:B10'}, we have
\begin{align*}
(L\widetilde{B}_{1,0})'&= (LB_{1,0})'-b_{1,0}(L\varphi)'
\\ &= \{(\lambda-3)yQ''+(3\lambda-1)Q'-yQ^2\}'
\\ &+ a_{1,0}\left\{
\frac{(\lambda-3)^2}{2}yQ''+\frac{\lambda^2-10\lambda+3}{2}Q'
+\frac{3-\lambda}{2}yQ^2 \right\}'.
\end{align*}
In view of \eqref{eq:Q}, we have 
  \begin{equation}
\label{eq:yq}
yQ''=-\frac{1}{4}L(y^2Q')-\frac{1}{2}Q'\quad\text{and}\quad yQ^2=-L(yQ)-2Q'.
  \end{equation}  
Using \eqref{eq:a10} and \eqref{eq:yq}, we see that
\begin{align*}
\widetilde  B_{1,0} + b_{1,0} \varphi &= \frac{3-\lambda}{4}y^2Q'+yQ-a_{1,0}\left\{
\frac{(\lambda-3)^2}{8}y^2Q'+\frac{3-\lambda}{2}yQ\right\}+b_{1,0}\varphi
\end{align*}
is solution to \eqref{eq:B10} as well as 
$B_{1,0}=\widetilde B_{1,0} +b_{1,0}\varphi + \kappa Q'$,
for any constant $\kappa$.
\end{proof}

\begin{lemma}[Resolution of $(\Omega_{2,0})$]\label{lem:b20}
There exists a solution $a_{2,0}$, $A_{2,0}$, $B_{2,0}$ of 
$(\Omega_{2,0})$ such that $A_{2,0}$ is even, $B_{2,0}$ is odd and
$$ \lim_{+\infty} A_{2,0}= - \frac 1{2(1-\lambda)} b_{1,0}^2=\gamma_{2,0},\quad
A_{2,0} -\gamma_{2,0} \in \mathcal{Y},$$
$$ \lim_{+\infty} B_{2,0}= b_{2,0},\quad
B_{2,0}- b_{2,0} \varphi \in \mathcal{Y},$$
where 
\begin{equation}\label{eq:20}
\forall \lambda\in (0,1),\quad  (1-\lambda) b_{2,0}\neq - \frac 1{6} b_{1,0}^3. 
\end{equation}
\end{lemma}
\begin{remark}
Existence of solutions to $\Omega_{2,0}$ follows from Proposition
\ref{prop:ex} below.
Note that we do not have to compute $A_{2,0}$ and $B_{2,0}$ to obtain
$b_{2,0}$.
The choice of that value of $\gamma_{2,0}$ above is justified in Section \ref{sec:32}.
\end{remark}
\begin{proof}
Recall that from Proposition \ref{prop:decomp}, the system $(\Omega_{2,0})$
writes:
\begin{align*}
(LB_{2,0})' &= (3-\lambda) A_{2,0}'' + 2 Q A_{2,0} + a_{2,0} (2\lambda-3) Q''\\
&\quad   +\frac{a_{1,0}}{2}\left\{(6\lambda-9)A_{1,0}'
+(\lambda-3)B_{1,0}''-2QB_{1,0}\right\}' \\ &
\quad +A_{1,0}^2+(A_{1,0}B_{1,0})'+A_{1,0}+\frac{\lambda-2}{1-\lambda}B_{1,0}'
+\frac{3}{2}(1-\lambda)a_{1,0}^2Q'',\\
(LA_{2,0})'&= a_{2,0}((\lambda-3)Q'' - Q^2)' + f'_{2,0} -\frac{2}{1-\lambda}QB_{1,0} \\
\text{where} \quad 
f_{2,0}&=a_{1,0}\{(\lambda-3)A_{1,0}''-2QA_{1,0}-Q\}
+(3-2\lambda)a_{1,0}^2Q''
+A_{1,0}^2\\
&\quad -\frac{1}{1-\lambda}A_{1,0}
+\frac{\lambda-3}{1-\lambda}B_{1,0}'.
\end{align*}
Let us integrate the equation of $B_{2,0}$ on $\R$. Since $b_{2,0}
=\lim_{+\infty} B_{2,0}
=\pm\lim_{\pm\infty} LB_{2,0}$, we have:
\begin{equation}\label{eq:b20}
2b_{2,0}=\int (2QA_{2,0} + A_{1,0}^2 +A_{1,0}) + \frac {\lambda - 2}{1-\lambda} 2b_{1,0}. 
\end{equation}
To determine $\int QA_{2,0}$,
we multiply the equation of $A_{2,0}$ by $\int_0^x P$ and use
$LP=-2Q$. We find
\begin{align*}
\int (L A_{2,0})' \int_0^x P& =  2\int A_{2,0} Q +\gamma_{2,0} \int P
= -a_{2,0} \int ((\lambda-3) Q''-Q^2) P +\int F_{2,0} \int_0^x P.
\end{align*}
Using 
$
-\int ((\lambda-3) Q''-Q^2) P= \frac {\lambda+1}2 \int Q^2,
$
the expression of $a_{2,0}$ in Claim \ref{lem:deta}, 
and then the expression of $a_{1,0}$ in Lemma \ref{lem:omega10}, we find
\begin{equation}\label{eq:firstA2Q}\begin{split}
2 \int A_{2,0} Q & = \frac {\lambda+1} 2 a_{2,0} \int Q^2  + \int F_{2,0} \int_0^x P
- \gamma_{2,0} \int P\\
&=  -\gamma_{2,0} \int (P-a_{1,0}  P_\lambda)
- a_{1,0} \int G_{2,0} Q + \int F_{2,0} \int_0^x (P- a_{1,0} P_\lambda)\\
&=  -\gamma_{2,0} \int (P-a_{1,0}  P_\lambda) 
- a_{1,0} \int G_{2,0} Q - \int f_{2,0}  (P- a_{1,0} P_\lambda)\\ &\quad 
- \frac 2{1-\lambda} \int QB_{1,0} \int_0^x (P- a_{1,0} P_\lambda).
\end{split}\end{equation}
Using Claim \ref{cl:phi}, we find
$$
\int_0^x (P- a_{1,0} P_\lambda) = \left(1+\frac {\lambda -3} 2 a_{1,0}\right) (xQ)+
b_{1,0} \frac {Q'} Q.
$$
Thus, using $\int Q'B_{1,0}=0$ 
we find
\begin{equation}\label{eq:secondA2Q}\begin{split}
2 \int A_{2,0} Q 
&=   2 \gamma_{2,0}  b_{1,0}
- a_{1,0} \int G_{2,0} Q - \int f_{2,0}  (P- a_{1,0} P_\lambda)\\ &\quad 
-\left( \frac 2{1-\lambda}\right) \left(1+\frac {\lambda -3} 2 a_{1,0}\right)\int xQ^2B_{1,0}.
\end{split}\end{equation}
Using \eqref{eq:secondA2Q} and \eqref{eq:b20}, we are able to find an expression of
$b_{2,0}=b_{2,0}(\lambda)$.
Using Mathematica (or by a long explicit computation using Claim \ref{cl:qint}),
we obtain the following 
\begin{equation}
b_{2,0}(\lambda)=\frac{4\,q(\lambda)}{5(1-\lambda) \left(15+10 \lambda-\lambda^2\right)^3},
\end{equation}
where
\begin{align*}
q(\lambda)&= 5625 -13500\lambda^2 -3375\lambda^3
+75(125-3 \pi ^2) \lambda^4+300(16-\pi ^2)\lambda^5
\\ & -10(18+7 \pi ^2)\lambda^6-5(45-4\pi ^2)\lambda^7-\pi^2 \lambda^8.
\end{align*}

\begin{claim}
  \label{le:nz}
For all $\lambda\in(0,1)$,
$$b_{2,0}\neq -\frac{b_{1,0}^3}{6(1-\lambda)}.$$
\end{claim}
\begin{proof}
  By a simple computation, we see that
  \begin{equation*}
b_{2,0}+\frac{b_{1,0}^3}{6(1-\lambda)}
=\frac{-4\lambda^2g(\lambda)}{5(1-\lambda)(15+10\lambda-\lambda^2)^3},
  \end{equation*}
where
\begin{align*}
  g(\lambda)&= 3375+3375\lambda-75(44-3\pi^2)\lambda^2
-300(16-\pi^2)\lambda^3
\\ &-5(207-14 \pi^2)\lambda^4+5(45-4 \pi ^2)\lambda^5+\pi^2\lambda^6.
\end{align*}
Since $$g^{(5)}(\lambda)=720\pi^2\lambda+600(45-4\pi^2)\ge g^{(5)}(0)
=600(45-4\pi^2)>0
\quad\text{for  $\lambda\in [0,1]$,}
$$
we see that $g^{(3)}(\lambda)$ is convex  on $[0,1]$.
We compute
$$g^{(3)}(\lambda)120\pi^2\lambda^3 +300(45-4\pi^2)\lambda^2-120(207-14\pi^2)\lambda
-1800(16-\pi^2),
$$
and $g^{(3)}(0)=1800(\pi ^2-16)$ and $g^{(3)}(1)=60(40\pi^2-669)$.
Since $g^{(3)}(\lambda)$ is convex on $[0,1]$, we have
$g^{(3)}(\lambda)\le \max(g^{(3)}(0),g^{(3)}(1))<0$.
Hence $g''(\lambda)$ is monotone decreasing and
$$g''(\lambda)\le g''(0)=-150(44-3 \pi^2)<0
\quad\text{for $\lambda\in[0,1]$.}$$
Thus we show that $g(\lambda)$ is concave on $[0,1]$.
By the concavity of $g$ and  the fact that $g(0)=3375>0$ and
$g(1)=144(4 \pi^2-15)>0$, we conclude that $g(\lambda)$ has no zero on $[0,1]$.
\end{proof}
The proof of Lemma \ref{lem:b20} is complete.
\end{proof}

\subsection{Resolution of the   systems $(\Omega_{1,1})$ and $(\Omega_{k,l})$ for $k+l=3$}\label{sec:2.3}

A main difference with the previous section is that we do not need to solve explicitely the systems $(\Omega_{1,1})$ and 
$(\Omega_{k,l})$ for $k+l=3$. Indeed, the existence of solutions satisfying some properties is sufficient for our purposes. In fact, to prove Theorem \ref{th:1}, this section is not necessary.
However, it is useful to give the sharp estimates of Proposition \ref{prop:approx}.

We claim the following result.

\begin{proposition}[Existence for a model problem]\label{prop:ex}
Let $F\in \mathcal{Y}$, odd and $G\in \mathcal{Y}$, even.
Let $\gamma \in \R$. Then, there exists $a,b\in \R$,   
$\widetilde A\in \mathcal{Y}$, even, and $\widetilde B\in \mathcal{Y}$, odd
such that
$$
A=\widetilde A+ \gamma,\quad B=\widetilde B+ b \varphi
$$
satisfy
\begin{equation*}
  (\Omega)\quad 
  \left\{
    \begin{aligned}
& (LA)'+a((3-\lambda)Q''+Q^2)'=F,\\
& (LB)'+a(3-2\lambda)Q'' -(3-\lambda)A''-2QA=G
    \end{aligned}
\right.
\end{equation*}
\end{proposition}

\begin{proof}
The proof is exactly the same as the one of Proposition 2.2 in \cite{MMcol1}, and even simpler since we deal with $F,G\in \mathcal{Y}$.
Set $
A=\widetilde A+ \gamma,$ $ B=\widetilde B+ b \varphi
$, where $\gamma$ is given, while $b$ is a parameter.
Since $(L 1)'=(1-2Q)'=-2Q'$, we obtain for $\widetilde A$, $\widetilde B$:
\begin{equation*}
    \left\{
    \begin{aligned}
& (L\widetilde A)'+a((3-\lambda)Q''+Q^2)'=F+2\gamma Q',\\
& (L\widetilde B)'+a(3-2\lambda)Q'' -(3-\lambda)\widetilde A''-2Q\widetilde A=G
+2\gamma Q -b (L\varphi)'.
    \end{aligned}
\right.
\end{equation*}
The function $F\in \mathcal{Y}$ begin odd, $\mathcal{H}(x)=\int_{-\infty}^x F(z) dz
+ 2 \gamma Q(x)$ belong to $\mathcal{Y}$ and is even.
By integration of the first line, we are reduced to solve
\begin{equation*}
    \left\{
    \begin{aligned}
&  L\widetilde A +a((3-\lambda)Q''+Q^2) =\mathcal{H},\\
& (L\widetilde B)'+a(3-2\lambda)Q'' -(3-\lambda)\widetilde A''-2Q\widetilde A=G
+2\gamma Q -b (L\varphi)'.
    \end{aligned}
\right.
\end{equation*}
Since $\int \mathcal{H} Q'=0$ (by parity) and $\mathcal{H}\in \mathcal{Y}$, 
by Lemma \ref{surL}, there exists $\overline H \in \mathcal{H}$, even, such that
$L\overline H= \mathcal{H}$. The function $V_\lambda$ being defined in \eqref{eq:sola10}, it follows that $\widetilde A=-a V_\lambda + \overline H$
is even, belongs to $\mathcal{Y}$ and solves the first line of the previous system.
Note that at this stage, the parameters $a$ and $b$ are still free.

Now, we only need to find $\widetilde  B\in \mathcal{Y}$, odd, such that
$(L\widetilde  B)'=-a Z_0 + D - b (L\varphi)',$
where 
$$
D = (3-\lambda) \overline{H}'' + 2Q \overline H + G + 2\gamma Q \in \mathcal{Y}, \text{ even, }
Z_0=(3-2 \lambda) Q'' + (3-\lambda) V_\lambda'' + 2 Q V_\lambda \in \mathcal{Y}, \text{ even}.
$$
Let 
$$
E=\int_0^x (D-aZ_0)(z) dz - b L \varphi.
$$
\begin{claim}\label{cl:ndg}
There exist $a$ and $b$ such that $E\in \mathcal{Y}$ and  $\int  E Q'=0$.
\end{claim}
Assuming Claim \ref{cl:ndg}, we fix  $a,b$ so that $E\in \mathcal{Y}$ and  $\int  E Q'=0$.
It follows from Lemma \ref{surL} that there exists $\widetilde B\in \mathcal{Y}$, odd,
such that $L\widetilde  B=E$. The solution is then given by $A=\widetilde A + \gamma$
and $B=\widetilde B + b \varphi$.
\end{proof}

\begin{proof}[Proof of Claim \ref{cl:ndg}]
First, we  check that $\int Z_0Q\neq 0$. Indeed,
by \eqref{eq:sola10}  and  \eqref{eq:on}
\begin{align*}
\int Z_0 Q &=(2\lambda -3) \int (Q')^2 - (LP_\lambda,V_\lambda)= (2\lambda -3) \int (Q')^2 + \int P_\lambda ((\lambda-3) Q'' - Q^2)\\&
= - \frac {15+10 \lambda -\lambda^2}{20} \int Q^2\neq 0.
\end{align*}
Thus, it suffices to choose $a=\int DQ/(\int Z_0 Q),$
and $b=\int_0^{+\infty} (D-aZ_0)(z) dz$ (note that
$\lim_{\pm \infty} L \varphi=\lim_{\pm \infty} \varphi=\pm 1$).
\end{proof}

\begin{lemma}[Resolution of $(\Omega_{1,1})$]\label{lem:b11}
There exists a solution $a_{1,1}$, $A_{1,1}$, $B_{1,1}$ of 
$(\Omega_{1,1})$ such that $A_{1,1}$ is even, $B_{1,1}$ is odd and
$$ \lim_{+\infty} A_{1,1}  =\gamma_{1,1}=\frac 12 b_{1,0}^2=-(1-\lambda) \gamma_{2,0} ,\quad
A_{1,1} -\gamma_{1,1} \in \mathcal{Y},$$
$$ \lim_{+\infty} B_{1,1}= b_{1,1},\quad
B_{1,1}- b_{1,1} \varphi \in \mathcal{Y}.$$
\end{lemma}

\begin{proof} From Proposition \ref{prop:decomp}, it is clear that
$F_{1,1}$ and $G_{1,1}$ satisfy the assumptions of Proposition \ref{prop:ex}.
The choice of $\gamma_{1,1}$ is justified in Section \ref{sec:32}.
In the rest of this paper, we will not need the expression of $b_{1,1}$
(note that it would be possible to compute it as in the proof of Lemma \ref{lem:b20}.
\end{proof}

From now on, 
we consider $(a_{k,l},A_{k,l},B_{k,l})$ defined for all $(k,l)\in \Sigma_0$, $1\leq k+l\leq 2$ in Lemmas \ref{lem:omega10}, \ref{lem:b20}, \ref{lem:b11}.
We now solve the systems $(\Omega_{k,l})$ for $k+l=3$ .

\begin{lemma}[Resolution of $(\Omega_{k,l})$ for $k+l=3$]\label{lem:2.5}
Let
\begin{align*}
&\gamma_{3,0}=\frac 5 {36(1-\lambda)^2} b_{1,0}^4+\frac {10}{3(1-\lambda)} d(\lambda)b_{1,0},\\
&\gamma_{2,1}= \frac 1{24(1-\lambda)} b_{1,0}^4 +\frac {\lambda}{2(1-\lambda)}b_{1,0}^2-\frac 1{1-\lambda} b_{1,0}b_{1,1}-4 d(\lambda) b_{1,0}
,\\& 
\gamma_{1,2}=-\frac 3{24} b_{1,0}^4 + b_{1,0}b_{1,1},
\end{align*}
where $d(\lambda)=b_{2,0}(\lambda)+\frac 1{6(1-\lambda)} b_{1,0}^3(\lambda)$.
For all $(k,l)\in \Sigma_0$   such that $k+l=3$,
there exists a solution $a_{k,l}$, $A_{k,l}$, $B_{k,l}$ of 
$(\Omega_{k,l})$ such that $A_{k,l}$ is even, $B_{k,l}$ is odd and
$$ \lim_{+\infty} A_{k,l}  =\gamma_{k,l},\quad
A_{k,l} -\gamma_{k,l} \in \mathcal{Y},$$
$$ \lim_{+\infty} B_{k,l}= b_{k,l},\quad
B_{k,l}- b_{k,l} \varphi \in \mathcal{Y}.$$
\end{lemma}

\begin{proof}
We claim the following
\begin{equation}\label{eq:cl}
\text{For all $k,l\in \Sigma_0$   such that $k+l=3$, we have
$F_{k,l}\in \mathcal{Y}$ is odd, $G_{k,l} \in \mathcal{Y}$ is even. }
\end{equation}
Assuming \eqref{eq:cl}, 
Lemma \ref{lem:2.5} is a direct consequence of Proposition \ref{prop:ex}.

Proof of \eqref{eq:cl}. 
Note that \eqref{eq:cl} follows from an algebraic "miracle", which we do not explain here. Indeed, several (bounded but) nonlocalized terms appear in the expression of 
$F_{k,l}$ and  $G_{k,l}$ for $k+l=3$ (see below), but all these terms eventually cancel.

To prove \eqref{eq:cl}, we look into more details   the proof of  Lemmas \ref{lem:SQ}--\ref{lem:Sint}.
First, the parity properties are clear by using the parity properties of
$Q$, $A_{k',l'}$, $B_{k',l'}$ for $k'+l'\leq 2$, the proof of  Lemmas \ref{lem:SQ}--\ref{lem:Sint} and Claims \ref{cl:SKdV1}, \ref{cl:SBBM0}.

Now, we recollect all the nonlocalized terms (due to $B_{1,0}$, $A_{1,1}$ and $A_{2,0}$) in $S(z)$ of order $\sigma^2 \wqs$ or $\sigma^2 \wqs'$. Note that terms containing derivatives of $B_{1,0}$, $A_{1,1}$ and
$A_{2,0}$ are in $\mathcal{Y}$ as well as terms of the kind $QB_{1,0}$.
Thus, we focus on the terms containing only $B_{1,0}$, $A_{1,1}$ and
$A_{2,0}$ without derivatives or multiplication by $Q$.
We skip the variables $\ys$ and $y$.

First, $S(Q)$ contains only localized terms.
Second, $\delta S_{KdV}(w)$ contains (see Claim \ref{cl:SKdV1}) 
$$
(\mu_\sigma-1) \wqs'' B_{1,0}+ (\mu_\sigma-1)(\wqs^2)'A_{2,0}
+(\mu_\sigma-1)\sigma \wqs'A_{1,1}
+ \wqs''''B_{1,0}+ (\wqs^2)''' A_{2,0}
+\sigma \wqs''' A_{1,1}.
$$
Third, from $S_{BBM}(w)$ we get (see Claim \ref{cl:SBBM0})
$$
\lambda \mu_\sigma \wqs''''(-B_{1,0}) + \lambda \mu_\sigma (\wqs^2)''' (-A_{2,0})
+\lambda \mu_\sigma \sigma\wqs''' (-A_{1,1}).
$$
Finally, from $S_{int}$, we obtain
$$
2 (\wqs \wqs')' B_{1,0} + 2 (\wqs^3)' A_{2,0}+ 2 \sigma (\wqs^2)' A_{1,1}
+ ((\wqs')^2)' B_{1,0}^2 .
$$
In the above formulas, we replace $\mu_\sigma$ by $1$ and
$\mu_\sigma-1$ by $(\lambda-1)\sigma$.
The contribution of $B_{1,0}$ at this order is
$$
(1) = B_{1,0} ((\lambda-1) \sigma \wqs + (1-\lambda) \wqs'' + \wqs^2)''.
$$
By \eqref{eq:qsigma}, $\wqs'' = \sigma \wqs -\frac 1{1-\lambda} \wqs^2 +O(\sigma^3)$,
$(1)$ gives only a lower order contribution.
The other terms are
\begin{align*}
  (2) &= A_{2,0} ((\lambda-1)\sigma \wqs^2 + (1-\lambda) (\wqs^2)'' + 2\wqs^3)'
\\& +A_{1,1} ((\lambda-1) \sigma^2 \wqs + (1-\lambda) \sigma \wqs'' + 2 \sigma \wqs^2)'
 +B_{1,0}^2 ((\wqs')^2)'.
\end{align*}
Using $\wqs'' = \sigma \wqs -\frac 1{1-\lambda} \wqs^2 +O(\sigma^3)$ and
$(\wqs')^2 = \sigma \wqs^2 -\frac 2{3(1-\lambda)} \wqs^3 + O(\sigma^4)$,
we find
$$
(2)=\sigma (\wqs^2)' (3 (1-\lambda)A_{2,0}+A_{1,1}+B_{1,0}^2)
+(\wqs^3)'(-\tfrac 43 A_{2,0} -\tfrac 2{3(1-\lambda)} B_{1,0}^2). 
$$
Using the following relations between the limits of $A_{2,0}$, $A_{1,1}$
and $B_{1,0}^2$ at $+\infty$ (see Lemmas \ref{lem:b20} and \ref{lem:b11}):
$$
\lim_{+\infty} A_{2,0}=-\frac 1{2(1-\lambda)} \lim_{+\infty} B_{1,0}^2,\quad
\lim_{+\infty} A_{1,1}=-(1-\lambda) \lim_{+\infty} A_{2,0},
$$
we observe that the functions of $y$ in $(2)$ are in fact all localized.
\end{proof}

\subsection{Recomposition of the approximate solution after the collision}\label{sec:32}

Let $1<c_2<c_1$, where $0<c_2-1<\epsilon_0$ is small and set
$$
\lambda=\frac {c_1-1}{c_1}, \quad \sigma=\frac {c_2-1}{c_2 \lambda}.
$$

We consider the function $z(t,x)$ defined by \eqref{eq:df1}--\eqref{eq:df3}
where for all $(k,l)\in \Sigma_0$, $a_{k,l}$, $A_{k,l}$, $B_{k,l}$  are chosen as in Lemmas \ref{lem:omega10}--\ref{lem:2.5}.

 We set
\begin{equation}\label{eq:dtau}
        \tau_{\sigma}=\sigma^{-\frac 12 -\frac 1{100}}
        =\left(\frac {c-1}{c\lambda}\right)^{-\frac 12 -\frac 1{100}}, \quad
        d(\lambda)=b_{2,0}(\lambda)+\frac 1{6(1-\lambda)} b_{1,0}^3(\lambda). 
\end{equation}
We claim the following result on $z$.
\begin{lemma}\label{lem:z}
\begin{equation*}
\forall t,x,\quad z(t,x)=z(-t,-x),
\end{equation*}
\begin{equation}\label{eq:z2}
 \forall t,\quad
\left\|  (1-\lambda \partial_x^2) \partial_t z  
+\partial_x(\partial_x^2 z -z  + z^2) \right\|_{H^1}\leq C \sigma^{\frac {15}4}.
\end{equation}
\begin{equation}\label{eq:z1}
\|z(\tau_\sigma)-\{Q(.-\tfrac 12 \delta)+\wqs(. +\mu_\sigma \tau_\sigma -\tfrac 12 \delta_\sigma) - d(\lambda) (\wqs^2)'(.+\mu_\sigma \tau_\sigma -\tfrac 12 \delta_\sigma)\}\|_{H^1}
\leq C \sigma^{\frac {13} 4},
\end{equation}
\begin{equation}\label{eq:z4}
\|z(-\tau_\sigma)-\{Q(.+\tfrac 12 \delta)+\wqs(. -\mu_\sigma \tau_\sigma +\tfrac 12 \delta_\sigma) + d(\lambda) (\wqs^2)'(.-\mu_\sigma \tau_\sigma +\tfrac 12 \delta_\sigma)\}\|_{H^1}
\leq C \sigma^{\frac {13} 4},
\end{equation}
where
\begin{equation}\label{eq:z3}
\delta=\sum_{(k,l)\in \Sigma_0} a_{k,l}\, \sigma^l \int \wqs^k,\quad
\tilde b_{1,1}=b_{1,1} -\frac 16 b_{1,0}^3,
\quad \delta_\sigma= 2 (b_{1,0}+ \sigma \tilde b_{1,1}).
\end{equation}
\end{lemma}
See Appendix C for the proof of Lemma \ref{lem:z}.

\subsection{Existence of the approximate $2$-soliton solution}\label{se:26}

The fact that $d(\lambda)\neq 0$ (see Claim \ref{le:nz}) in Lemma \ref{lem:z} means formally that the collision is not elastic and that the residue due to the collision is of order 
$(\wqs^2)'$.
However, the approximate solution $z(t,x)$ given in Lemma \ref{lem:z} being symmetric, it contains the residue at both $-\tau_\sigma$ and $+\tau_\sigma$ (see \eqref{eq:z1}, \eqref{eq:z4}). To match the solution $u(t)$ considered in Theorem \ref{th:1}, which is pure at $-\infty$, we need to introduced a different approximate solution, which, at the main orders, will contain a residue only at $+\tau_\sigma$. 

\begin{proposition}\label{prop:approx}
There exists a function $z_{\#}$ of the form \eqref{eq:df1}--\eqref{eq:df2} such that
\begin{equation}\label{eq:zd}
 \forall t\in [-\tau_\sigma,\tau_\sigma],\quad
\left\|  (1-\lambda \partial_x^2) \partial_t z_{\#} 
+\partial_x(\partial_x^2 z_{\#} -z_{\#} + z_{\#}^2) \right\|_{H^1}\leq C \sigma^{3},
\end{equation}
\begin{equation}\label{eq:zT}\begin{split}
        &\left \|z_{\#}(\tau_\sigma)-\{ Q(.-\tfrac 12 \delta) 
        + \wqs(.+\mu_\sigma \tau_\sigma - \tfrac 12 \delta_\sigma)
        -2 d(\lambda) (\wqs^2)'(.+\mu_\sigma \tau_\sigma-\tfrac 12 \delta_\sigma) \}\right\|_{H^1}
                  \\
        & +\left\|z_{\#}(-\tau_\sigma)- \{Q(.+\tfrac 12 \delta) 
        + \wqs(.-\mu_\sigma \tau_\sigma+\tfrac 12 \delta_\sigma)\}
         \right\|_{H^1}\leq  C \sigma^{\frac {11}4},
\end{split} \end{equation}
where 
$$
\forall \lambda\in (0,1),\quad 
d(\lambda)\neq 0\quad \text{(see Lemma \ref{lem:b20})},
$$
\begin{equation}\label{eq:dt}
        \left|\delta - \sigma^{\frac 12} \frac {10(1-\lambda^2)}{15+10 \lambda - \lambda^2}\int Q\right|
        \leq C \sigma^{\frac 32},       \quad
        \left|\delta_\sigma - \frac {-60 + 36 \lambda^2}{15 + 10 \lambda - \lambda^2}\right| \leq C \sigma.
\end{equation}
\end{proposition}
Using the change of variable \eqref{eq:chv}, we define
\begin{equation}\label{ztov}
v(t,x)=\frac {\lambda} {1-\lambda} z_{\#}(t',x'),       \quad 
        D=\frac {1-\lambda}{\lambda^{\frac 32}} d(\lambda),
\end{equation}
\begin{equation}\label{eq:DT}
        T =\frac {1-\lambda}{\lambda^{\frac 32}} \tau_\sigma
        = \left(\frac {c_2-1}{c_2}\right)^{-\frac 12- \frac 1{100}}
        \left(\frac {1-\lambda}{\lambda}\right) \lambda^{\frac 1{100}}.
\end{equation}
We obtain the following consequence of Proposition \ref{prop:approx}.

\begin{proposition}\label{cor:1}
The function  $v$ defined by \eqref{ztov} where $z_{\#}$ is as in
 Proposition \ref{prop:approx} satisfies, for some constant $C=C(c_1)>0$,
 \begin{equation*}\label{eq:approx1}
 \forall t\in [-T,T],\quad
\left\| (1-  \partial_x^2) \partial_t v +\partial_x(v+v^2) \right\|_{H^1}\leq C (c_2-1)^{3}.
\end{equation*}
\begin{equation*}\begin{split} 
&
\left\| v(T) -\{ \varphi_{c_1}(.-c_1 T -\tfrac 12 \Delta_1)
+\varphi_{c_2}(.-c_2T-\tfrac 12 \Delta_2)
-2 D (\varphi_{c_2}^2)'(.-c_2T-\tfrac 12 \Delta_2)\}
\right\|_{H^1} 
\\&
+ \left\| v(-T) - \{\varphi_{c_1}(.+c_1 T+ \tfrac 12 \Delta_2)
+\varphi_{c_2}(.+c_2T+\tfrac 12\Delta_2)\}
\right\|_{H^1}\leq  K (c_2-1)^{\frac {11}4}, 
\end{split}\end{equation*}
where 
$$\forall c_1>1, \quad D=D(c_1)\neq 0,$$
and
\begin{equation}\label{eq:Dt}
        \left|\Delta_1 - \sqrt{  {c_2-1} }  \frac {10(1-\lambda^2)}{\lambda(15+10 \lambda - \lambda^2)}\int Q\right|
        \leq C (c_2-1)^{\frac 32},      \quad
        \left|\Delta_2 - \frac {-30 + 18 \lambda^2}{\sqrt{\lambda}(15 + 10 \lambda - \lambda^2)}\right| \leq C (c_2-1).
\end{equation}
\end{proposition}
\begin{remark}\label{rk:no}
Note that  comparing \eqref{eq:zd} and \eqref{eq:z2}, there is a loss of $\sigma^{\frac 34}$
when changing the conditions at $\pm \tau_\sigma$ for nonsymmetric conditions.
In fact, while \eqref{eq:z2} can be improved by refining further the function $z$
(i.e. taking a larger set of indices $\Sigma_0$), it seems that by the method of this paper, one cannot improve estimate \eqref{eq:zd} on an approximate solution satisfying the conditions 
at $\pm\tau_\sigma$ as in Proposition \ref{prop:approx}.
See the remark after Proposition 5.2 in \cite{MMcol1} for a similar problem for the case of the quartic gKdV equation. As a consequence, we do not obtain optimal estimates in \eqref{eq:th:1:3}.
\end{remark}
\noindent\emph{Proof of Proposition \ref{cor:1} from Proposition \ref{prop:approx}.}
The first estimate is a consequence of
\begin{equation*}
 (1-\partial_x^2)\partial_t v + \partial_x(v+v^2)=
  \frac {\lambda^{\frac 52}} {(1-\lambda)^2} \left\{
(1-\lambda \partial_{x'}^2) \partial_{t'} z
+\partial_{x'} (\partial_{x'}^2 z - z + z^2)\right\}.
\end{equation*}

By \eqref{eq:v1} for $c=c_1$ and $c=c_2$ and by \eqref{eq:v2} for $c=c_2$, we have
\begin{align*}
          &Q(x'-\tfrac 12 \delta) 
        + \wqs(x'+\mu_\sigma \tau_\sigma - \tfrac 12 \delta_\sigma)
        -2 d(\lambda) (\wqs^2)'(x'+\mu_\sigma \tau_\sigma-\tfrac 12 \delta_\sigma) \\&  =\frac {1-\lambda}{\lambda} \left\{
        \varphi_{c_1}(x-c_1 T-\tfrac \delta{2 \sqrt{\lambda}}) + 
        \varphi_{c_2}(x-c_2 T -\tfrac {\delta_\sigma}{2 \sqrt{\lambda}}) -
        2 \frac {1-\lambda}{\lambda^{\frac 32}} d(\lambda)
        (\varphi_{c_2}^2)'(x-c_2 T -\tfrac {\delta_\sigma}{2 \sqrt{\lambda}}) \right\},
         \\
        & Q(x'+\tfrac 12 \delta) 
        + \wqs(x'-\mu_\sigma \tau_\sigma+\tfrac 12 \delta_\sigma)\\&
        =\frac {1-\lambda}{\lambda} \left\{
        \varphi_{c_1}(x+c_1T+\tfrac \delta{2 \sqrt{\lambda}}) + 
        \varphi_{c_2}(x+c_2 T +\tfrac {\delta_\sigma}{2 \sqrt{\lambda}}) \right\}.
 \end{align*}
 Using these identities and the estimates on $z_{\#}$, we
finish the proof of the proposition.

\medskip

\noindent\emph{Proof of Proposition \ref{prop:approx}.}
As in Proposition 5.2 of \cite{MMcol1}, we modify the function $z(t,x)$ constructed in Lemma \ref{lem:z} in the following way: let
$$z_{\#}(t,x)=z(t,x)+w_{\#}(t,x),\quad
w_{\#}(t,x)=-d(\lambda)(\wqs^2)'(\ys)(1-P(y)),$$ where $P$ is defined in \eqref{eq:defP}. 

\medskip

Proof of \eqref{eq:zT}.  
Replacing $z=z_{\#}-w_{\#}$ in \eqref{eq:z1}, we have
\begin{equation*}\begin{split}
&\left\|z_{\#}(\tau_\sigma)-\{ Q(.-\tfrac 12 \delta) 
        + \wqs(.+\mu_\sigma \tau_\sigma - \tfrac 12 \delta_\sigma)
        - d(\lambda) (\wqs^2)'(.+\mu_\sigma \tau_\sigma-\tfrac 12 \delta_\sigma) \}
        -w_{\#}(\tau_\sigma)\right\|_{H^1}\\&
\leq C \sigma^{\frac {13}4}.
\end{split}\end{equation*}
Thus, using \eqref{eq:de} ($P\in \mathcal{Y}$)
\begin{equation*} \begin{split}
        &\left \|z_{\#}(\tau_\sigma)-\{ Q(.-\tfrac 12 \delta) 
        + \wqs(.+\mu_\sigma \tau_\sigma - \tfrac 12 \delta_\sigma)
        -2 d(\lambda) (\wqs^2)'(.+\mu_\sigma \tau_\sigma-\tfrac 12 \delta_\sigma) \}\right\|_{H^1}\\& \leq C\sigma^{\frac {13}4} + 
        \left\| d(\lambda) (\wqs^2)'(.+\mu_\sigma \tau_\sigma-\tfrac 12 \delta_\sigma) 
        +w_{\#}(\tau_\sigma)\right\|_{H^1}
        \\& \leq C\sigma^{\frac {13}4} + 
        C \left\|  (\wqs^2)'(.-\tfrac 12 \delta_\sigma) 
        -(\wqs^2)' \right\|_{H^1}\leq C \sigma^{\frac {11}4}.
        \end{split}\end{equation*}
        Similarly,
                 \begin{equation*} \begin{split} 
        & \left\|z_{\#}(-\tau_\sigma)- \{Q(.+\tfrac 12 \delta) 
        + \wqs(.-\mu_\sigma \tau_\sigma+\tfrac 12 \delta_\sigma)
        + d(\lambda) (\wqs^2)'(.-\mu_\sigma \tau_\sigma+\tfrac 12 \delta_\sigma)\}
        -w_{\#}(-\tau_\sigma)\right\|_{H^1}\\& \leq  C \sigma^{\frac {13}4},
\end{split} \end{equation*}
so that
 \begin{equation*}  \begin{split} 
         & \left\|z_{\#}(-\tau_\sigma)- \{Q(.+\tfrac 12 \delta) 
        + \wqs(.-\mu_\sigma \tau_\sigma+\tfrac 12 \delta_\sigma)\}
        \right\|_{H^1}\\&  \leq  C \sigma^{\frac {13}4}+C \left\|  (\wqs^2)'(.+\tfrac 12 \delta_\sigma) 
        -(\wqs^2)' \right\|_{H^1}\leq C \sigma^{\frac {11}4}.
\end{split}  \end{equation*}

Note that \eqref{eq:dt} is a consequence of \eqref{eq:z3}.

\medskip

Proof of \eqref{eq:zd}. 
Let
$$
S_{\#}(t,x)=  (1-\lambda \partial_x^2) \partial_t z_{\#} 
+\partial_x(\partial_x^2 z_{\#} -z_{\#} + z_{\#}^2) 
=S(t,x)+ \delta S(w_{\#}) + \partial_x((z+w_{\#})^2 - z^2 - 2 Qw_{\#}).
$$
(See \eqref{eq:sion2} for the notation $\delta S$.)
We claim
\begin{equation}\label{eq:dS}
\|\delta S(w_{\#}) \|_{H^1}\leq C \sigma^{3}.
\end{equation}
Indeed, from Claims \ref{cl:SKdV1} and \ref{cl:SBBM0}, the lower order term in $\delta S(w_{\#})
=\delta S_{KdV}(w_{\#}) + S_{BBM}(w_{\#})$ is
$d(\lambda) (\wqs^2)'(\ys)(L(1-P))'$, but this term is zero since
$(L(1-P))'=(1-2Q+2Q)'=0$. All the other terms are controlled in $H^1$ by $\sigma^3$. For example, the next term (in increasing order of powers of $\sigma$) is
\begin{equation}\label{eq:fu}
(\wqs^2)''(\ys)((3-\lambda)(1-P)'' + 2(1-P)Q)(y)= (\wqs^2)''(\ys)(-(3-\lambda)P'' + 2Q-2PQ)(y).
\end{equation}
By \eqref{eq:qsigma}, this term is exactly of size $\sigma^3$ in $H^1$.
Indeed, for any $f\in \mathcal{Y}$, we have
$$
\|(\wqs^2)''(\ys)f(y)\|_{L^2}\leq \|(\wqs^2)'' \|_{L^\infty}\|f(y)\|_{L^2}
\leq C \sigma^3,
$$
and similarly for the $H^1$ norm.
Note that the function $(-(3-\lambda)P'' + 2Q-2PQ)$ being not orthogonal to $Q$, we cannot remove it by adding a further term $(\wqs^2)''$ in $z_{\#}$ (see proof of Proposition \ref{prop:ex}).

Finally, we claim the following, which completes the proof of Proposition \ref{prop:approx}.
\begin{equation}\label{eq:in}
\|\partial_x((z+w_{\#})^2 - z^2 - 2 Qw_{\#})\|_{H^1}\leq 
C\sigma^{\frac 72}.
\end{equation}
Indeed, note that
$\partial_x((z+w_{\#})^2 - z^2 - 2 Qw_{\#})=\partial_x(2(z-Q)w_{\#} + w_{\#}^2)$, and \eqref{eq:in} follows easily from the expression of $z$ and $w_{\#}$.

\section{Preliminary results for stability of the 2-soliton structure}\label{sec:3}

In this section, we gather several stability results. Section \ref{sec:3.1} concerns the
stability of $v(t)$ by the BBM equation during the interaction, i.e. on the time interval $[-T,T]$.
Indeed, we control the difference between the approximate solution $v(t)$ constructed in Proposition \ref{cor:1} and a solution of \eqref{eq:BBM}. We have to use long time stability arguments since 
$T\to +\infty$ as $c_2\to 1^+$. We use a functional which is at the first order the functional introduced by Weinstein \cite{We3}Êto prove the stability of one solitary wave for the BBM equation.
Some nonlinear corrective terms are added and another corrective term is needed to take into account the trajectory of the soliton $\varphi_{c_1}$, which is not a straight line (see formula of $y_1$ in \eqref{eq:y1} below).

Section \ref{sec:3.2}
concerns the large time behavior after interaction, i.e. for $t>T$. 
Here the techniques are related to the stability of the dynamics of two solitons and
involve monotonicity properties and Liouville theorems.
These results are refinements of the following works: \cite{Mi2}, \cite{Di}, \cite{DiMa}, \cite{Ma} concerning the BBM equation.  See also \cite{MMcol1} and references therein for the gKdV case.

\subsection{Dynamic stability in the interaction region}\label{sec:3.1}
For any $c>1$ sufficiently close to $1$, we consider the function $z_{\#}(t)$ of the form
\begin{equation*}
    z_{\#} (t',x')=Q(y)+\wqs(y_\sigma)+\sum_{(k,l)\in \Sigma_0} 
        \sigma^l \left\{\wqs^k(y_\sigma) A_{k,l}(y)+(\wqs^k)'(y_\sigma) B_{k,l}(y)\right\}
\end{equation*}
defined in Proposition \ref{prop:approx}
(recall that   $y$, $y_\sigma$ are defined in \eqref{eq:df1}). As in Proposition \ref{cor:1}, we set
\begin{equation}\label{eq:ddz}\begin{split}
    v (t,x)& =\frac {\lambda}{1-\lambda} z_{\#} (t',x')\\& 
    =\varphi_{c_1}(y_1)+\varphi_{c_2}(y_2)+\sum_{(k,l)\in \Sigma_0} 
        \sigma^l \left\{\varphi_{c_2}^k(y_2) \widetilde A_{k,l}(y_1)+(\varphi_{c_2} ^k)'(y_2) \widetilde B_{k,l}(y_1)\right\},
        \end{split}
\end{equation}
where (see proof of Claim \ref{cl:2.1})
\begin{equation}\label{eq:y1}
y_2= x-c_2 t,
\quad 
y_1=\frac y{\sqrt{\lambda}}=x - c_1 t
-  \widetilde  \alpha(y_2),\quad \widetilde \alpha(y_2)=\frac 1{\sqrt{\lambda}} \alpha(\sqrt{\lambda}y_2), 
\end{equation}
$$
\widetilde A_{k,l}(y_1)=\left(\frac {1-\lambda}{\lambda}\right)^{k-1}
A_{k,l}(\sqrt{\lambda} y_1),\quad
\widetilde B_{k,l}(y_1)=\left(\frac {1-\lambda}{\lambda}\right)^{k-1}
\frac 1 {\sqrt{\lambda}} B_{k,l}(\sqrt{\lambda} y_1).
$$
Now, we set
$$S(t)=(1-  \partial_x^2) \partial_t v +\partial_x(v+v^2).$$

\begin{proposition}[Exact solution close to the approximate solution $v$]\label{prop:I}
      Let $\theta>1$. There exists $\epsilon_0>0$ such that the following holds for any $0<c_2-1<\epsilon_0$.
    Suppose that
    \begin{equation}\label{INTkl}
        \forall t\in [-T,T],\quad 
        \left\| S(t)\right\|_{H^1(\mathbb{R})}\leq K \frac {(c_2-1)^\theta}{T},    \end{equation}
        ($T$ being defined in \eqref{eq:DT}),
    and   for some $T_0\in [-T,T]$,
    \begin{equation}\label{hypINT}
        \| u(T_0) - v(T_0) \|_{H^1(\mathbb{R})}\leq K (c_2-1)^{\theta},
    \end{equation}
where  $u(t)$ is an $H^1$ solution of \eqref{eq:BBM}. 
    Then, there exist $K_0=K_0(\theta,K,\lambda)$ and a function $\rho:[-T,T]\rightarrow \mathbb{R}$ such that, for all $t\in [-T,T]$,
    \begin{equation}\label{INT41}
        \|u(t)-v(t,.-\rho(t)) \|_{H^1} \leq K_0 (c_2-1)^{\theta},\quad |\rho'(t)|\leq K_0 (c_2-1)^{\theta}.
    \end{equation}
\end{proposition}
\begin{remark}\label{rk:th}
Note that the proof of Proposition \ref{prop:I} is based on long time stability methods. It is similar to the proof of Proposition 4.1 in \cite{MMcol1}.
\end{remark}

\begin{proof}
We  prove the result on $[T_0,T]$. By using the transformation $x\to -x$, $t\to -t$, the proof is the same on $[-T,T_0]$.
Let $K^*>K$ be a constant to be fixed later. Since $\|u(T_0)-v(T_0)\|_{H^1}\leq K(c_2-1)^\theta$, by continuity in time in $H^1(\mathbb{R})$, there exists $T^*>T_0$ such that
\begin{equation*}
    T^*=\sup\left\{T_1\in [T_0,T]\,|\,\text{$\exists r\in C^1([T_0,T_1])$
s.t. }\sup_{t\in[T_0,T_1]}
\|u(t){-}v(t,.{-}r(t))\|_{H^1}\leq K^* (c_2{-}1)^{\theta} \right\}.
\end{equation*}
Note that the translation direction is degenerate and without the freedom in the translation parameter, the result would not be correct.
The objective is to prove that $T^*=T$ for $K^*$ large. For this, we argue by contradiction, assuming that $T_0<T^*<T$ and reaching a    contradiction with the definition of $T^*$ by proving independent estimates on $\|u(t)-v(t,.-r)\|_{H^1}$ on $[T_0,T^*]$.

First, we claim some estimates related to $v$.

\begin{claim}[Preliminary estimates]\label{TRANScl}
       \begin{equation}\label{Tcl1}
        \|(1-\partial_x^2)(\partial_t v+c_1\partial_xv)(t)\|_{L^\infty}
        +  \|\partial_t^2 \partial_x^2 v(t)+c_1\partial_t \partial_x^3 v(t)\|_{L^\infty} \leq K (c_2-1),
    \end{equation}
    \begin{equation}\label{Tcl2}\begin{split}
&        \|\partial_t v(t)+c_1\partial_xv(t)+ (c_1-c_2)\widetilde \alpha'(y_2) \varphi_{c_1}'(y_1)\|_{L^2}\leq  K (c_2-1)^{\frac 54},\\
&        \|\partial_t v(t)+c_1\partial_xv(t) + (c_1-c_2)\widetilde \alpha'(y_2) \varphi_{c_1}'(y_1)\|_{L^\infty}\leq  K (c_2-1)^{\frac 32},
\end{split}    \end{equation}
        \begin{equation}\label{Tcl4}
        \|\partial_x v -  \varphi_{c_1}'(y_1) \|_{L^2}\leq K (c_2-1) ,
    \end{equation}
    \begin{equation}\label{Tcl5}
        \|\widetilde \alpha''(y_2)\|_{L^\infty}+\frac 1{c_2-1} \|\widetilde \alpha^{(4)}(y_2)\|_{L^\infty}\leq K (c_2-1)^{\frac 32}.
    \end{equation}\end{claim}

\begin{proof}[Proof of Claim \ref{TRANScl}]
These estimates are simple consequences of \eqref{eq:ddz} and elementary calculations.
\end{proof}

\emph{Step 1.} Choice of the translation parameter.

\begin{lemma}[Modulation]\label{TRANS}
    There exists a $C^1$ function $\rho:[T_0,T^*]\rightarrow \mathbb{R}$ such that, for all $t\in [T_0,T^*]$,
    the function $\varepsilon(t)$ defined by
    $
        \varepsilon(t,x)=u(t,x+\rho(t))-v(t,x)
    $
    satisfies, $ \forall t\in [T_0,T^*],$
    $$             \int \varepsilon(t,x) (1-\partial_x^2)(\varphi_{c_1}'(y_1)) dx =0,\label{TRANS1}
     $$ and for $K$ independent of $K^*$,
        \begin{equation}\label{TRANS3}\begin{split}
          &  \|\varepsilon(t)\|_{H^1}\leq  2 K^* (c_2-1)^{\theta},\quad
            |\rho(T_0)|+\|\varepsilon(T_0)\|_{H^1}\leq K (c_2-1)^{\theta},\\
& |\rho'(t)|\leq K \|\varepsilon(t)\|_{H^1}+K \|S(t)\|_{H^1}.
 \end{split}       \end{equation}
\end{lemma}
\begin{proof}[Proof of Lemma \ref{TRANS}]
The result follows from a standard argument. Let $t\in [T_0,T^*]$ and
\begin{equation*}
    \zeta(U,r)=\int (U(x+r)-v(t,x)) (1-\partial_x^2)(\varphi_{c_1}'(y_1)) dx.
\end{equation*}
Then $
    \frac{\partial \zeta}{\partial r}(U,r)= \int U'(x+r) (1-\partial_x^2)(\varphi_{c_1}'(y_1)) dx, 
$
so that
from Claim \ref{TRANScl} (see \eqref{Tcl4}),  for $(c_2-1)$ small enough,
\begin{equation*}\begin{split}
    &\frac{\partial \zeta}{\partial r}(v(t),0)= \int (\partial_x v)(t,x) (1-\partial_x^2)(\varphi_{c_1}'(y_1)) dx \\
&    >\int [(\varphi_{c_1}'')^2+(\varphi_{c_1}')^2]dx -K(c_2-1) 
   >\frac 12 \int  [(\varphi_{c_1}'')^2+(\varphi_{c_1}')^2]dx .
\end{split}\end{equation*} Since    $\zeta(v,0)=0$, for $U$ close to $v(t)$ in $L^2$ norm, the existence    of a unique 
$\rho(U)$ satisfying $\zeta(U(x-\rho(U)),\rho(U))=0$ is  a consequence of the Implicit Function Theorem.

From the definition of $T^*$, it follows that there exists $\rho(t)=\rho(u(t))$, such that 
$\zeta(u(x-\rho(t)),\rho(t))=0$.
We set
\begin{equation}\label{defz}
    \varepsilon(t,x)=u(t,x+\rho(t))-v(t,x),
\end{equation}
so that $\int \varepsilon(t)(1-\partial_x^2)(\varphi_{c_1}'(y_1))=0$ follows from the definition of $\rho(t)$. Estimate $\|\varepsilon(t)\|_{H^1}\leq  2 K^* (c_2-1)^{\theta}$ follows from the Implicit
Function Theorem and the definition of $K^*$.
Moreover, since $\|u(T_0)-v(T_0)\|\leq (c_2-1)^{\theta}$, we have $
    |\rho(T_0)|+\|\varepsilon(T_0)\|_{H^1}\leq K (c_2-1)^{\theta}$,
where $K$ is independent of $K^*$.

Now, let us prove that
\begin{equation}\label{estxxX}
    |\rho'(t)|\leq K \|\varepsilon(t)\|_{H^1} + K \|S(t)\|_{H^1}.
\end{equation}
From the definition of $\varepsilon(t)$ and $S(t)$ and $u(t)$ being a solution of the (BBM) equation, we have
\begin{align}
    (1-\partial_x^2)\partial_t \varepsilon + \partial_x (\varepsilon + (\varepsilon+v)^2 -v^2) 
    & =-[(1-\partial_x^2)\partial_{t} v + \partial_{x} (v + v^2)]+  \rho'(t) (1-\partial_x^2)\partial_x (v+\varepsilon) \nonumber\\
    & = -S(t) +  \rho'(t)  (1-\partial_x^2) \partial_x(v+\varepsilon).
    \label{eqz}
\end{align}
Since
$
 \int \varepsilon(t)(1-\partial_x^2)(\varphi_{c_1}'(y_1))dx = 0,
$
we have
\begin{equation*}\begin{split}
0 & = \frac d{dt} \int  \varepsilon(t)(1-\partial_x^2)(\varphi_{c_1}'(y_1)) 
= \int [(1-\partial_x^2) \partial_t \varepsilon ]\varphi_{c_1}'(y_1) 
+\int \varepsilon (1-\partial_x^2) [\partial_t (\varphi_{c_1}'(y_1) )]\\
& = -\int \partial_x( \varepsilon+(\varepsilon+v)^2 - v^2) \varphi_{c_1}'(y_1)
-\int S(t) \varphi_{c_1}'(y_1) + \rho'(t) \int [(1-\partial_x^2)\partial_x(v+\varepsilon)] \varphi_{c_1}'(y_1) \\&+
\int \varepsilon (1-\partial_x^2) [-c_1 \varphi_{c_1}''(y_1)+c_2 \widetilde \alpha'(y_2) \varphi_{c_1}''(y_1)].
\end{split}\end{equation*}
Thus, on the one hand
\begin{equation}\label{exprxxX}
    \begin{split}
&         \rho'(t)\int (v+\varepsilon) [(1-\partial_x^2)\partial_x (\varphi_{c_1}'(y_1))]
=-\int S(t) \varphi_{c_1}'(y_1)
\\ &
+\int \varepsilon [(1+ 2 v+\varepsilon) \partial_x (\varphi_{c_1}'(y_1) +(1-\partial_x^2) (-c_1 \varphi_{c_1}''(y_1) + c_2 \widetilde \alpha'(y_2) \varphi_{c_1}''(y_1)) ],
\end{split}
\end{equation}
and so
\begin{equation}\label{eq:34}
\left| \rho'(t)\int (v+\varepsilon) [(1-\partial_x^2)\partial_x (\varphi_{c_1}'(y_1))] \right|
\leq C \|\varepsilon(t)\|_{L^2} + \|S(t)\|_{L^2}.
\end{equation}
On the other hand,
\begin{align*}
 \int (v+\varepsilon) [(1-\partial_x^2)\partial_x (\varphi_{c_1}'(y_1))]
& = \int  [(1-\partial_x^2)(\varphi_{c_1}(y_1))] [\partial_x (\varphi_{c_1}'(y_1))]\\&
+\int (v-\varphi_{c_1}(y_1) +\varepsilon)  [(1-\partial_x^2)\partial_x (\varphi_{c_1}'(y_1))],
\end{align*}
and for $c_2-1<\epsilon_0$ small enough,
$$-\int   [(1-\partial_x^2)(\varphi_{c_1}(y_1))] [\partial_x (\varphi_{c_1}'(y_1))] \geq -\frac 34 
\int (\varphi_{c_1}-\varphi_{c_1}'')\varphi_{c_1}''=\frac 34 \int  (\varphi_{c_1}')^2+(\varphi_{c_1}'')^2>0,$$
so  that $|\int (v+\varepsilon) (1-\partial_x^2)\partial_x (\varphi_{c_1}'(y_1))|\geq \frac 12 
\int (\varphi_{c_1}')^2+(\varphi_{c_1}'')^2  $ for $c_2-1$ small and \eqref{estxxX} follows from \eqref{eq:34}.
\end{proof}

\emph{Step 2.} Control of the direction $\int \varepsilon  (1-\partial_x^2)(\varphi_{c_1}(y_1)) $. In this proof, the use of the
invariant $N(u(t))$ (see \eqref{eq:cv}) replaces a modulation argument in the scaling parameter.

\begin{lemma}[Control of the negative direction]\label{Qdir}
    For all $t\in [T_0,T^*]$,
    \begin{equation}\label{Qdir1}
        \left|\int  \varepsilon(t) (1-\partial_x^2)\varphi_{c_1}(y_1)dx\right|\leq  K (c_2-1)^\theta + K (c_2-1)^{\frac 34} \|\varepsilon(t)\|_{L^2}+ K \|\varepsilon(t)\|_{H^1}^2.
    \end{equation}
\end{lemma}

\noindent\emph{Proof of Lemma \ref{Qdir}.}
Remark that since $v(t)$ is an approximate solution of \eqref{eq:BBM}, $N(v(t))$ has a small variation.
Indeed, by multiplying the equation $S(t)=(1-\partial_x^2) \partial_t v + \partial_x(\partial_x^2 v -v + v^2)$ by $v$ and integrating,
we obtain
$
    \left| \frac d {dt} N(v(t))\right|=\left| \int S(t,x) v(t,x) dx\right|
    \leq K \|S(t)\|_{L^2}.
$
Thus,    
\begin{equation}\label{vL2}
\forall t\in [T_0,T^*],\quad
    \left|N(v(t))-N(v(T_0))\right| \leq K T \sup_{t\in [-T,T]} \|S(t)\|_{H^1}
    \leq K (c_2-1)^\theta. 
\end{equation}
Since $u(t)$ is a solution of the (gKdV) equation, we have
\begin{equation}\label{uL2}
   N(u(t))=N(v(t)+\varepsilon(t)) = N(u(T_0))=N(v(T_0)+\varepsilon(T_0)). 
\end{equation}
By expanding \eqref{uL2} and using \eqref{vL2} and \eqref{TRANS3}, we obtain:
\begin{equation*}\begin{split}
        2 \left|\int ((1-\partial_x^2) v(t)) \varepsilon(t)\right|
        &\leq K (c_2-1)^\theta + 2 \left|\int ((1-\partial_x^2) v(T_0))   \varepsilon(T_0)\right| + \|\varepsilon(T_0)\|_{H^1}^2 + \|\varepsilon(t)\|_{H^1}^2 
        \\& \leq K (c_2-1)^\theta  + \|\varepsilon(t)\|_{H^1}^2.
\end{split}\end{equation*}
 Using this  and $\|(1-\partial_x^2)(v(t)-\varphi_{c_1}(y_1))\|_{L^2}\leq K (c_2-1)^{\frac 34}$, we obtain:
\begin{equation*}\begin{split}
     \left|\int  \varepsilon(t) (1-\partial_x^2)\varphi_{c_1}(y_1)dx\right|
     &\leq \left| \int \varepsilon(t) [(1-\partial_x^2)(v(t)-\varphi_{c_1}(y_1)]\right|+ \left|\int \varepsilon(t) ((1-\partial_x^2) v(t))\right|
     \\ &\leq  K (c_2-1)^\theta + K (c_2-1)^{\frac 34} \|\varepsilon(t)\|_{L^2}+
     K  \|\varepsilon(t)\|_{H^1}^2.
\end{split}\end{equation*}

\medskip

\emph{Step 3.} Energy functional for  $\varepsilon(t)$.
We set
\begin{equation*}
    \mathcal{F}(t)=\frac 12 \int ((c_1-1)\varepsilon^2 + c_1(\partial_x \varepsilon)^2     -\tfrac 23   \left((v+\varepsilon)^3- v^{3}-3 v^2 \varepsilon\right)
    + \frac 12(c_1-c_2) \int \widetilde \alpha'(y_2) ((\partial_x\varepsilon)^2+\varepsilon^2).
\end{equation*}

\begin{claim}[Coercivity of $\mathcal{F}$]\label{POSf}
    There exists $\kappa_0>0$ such that, for $c_2-1$ small enough,
    \begin{equation}\label{posf1}
       \|\varepsilon(t)\|_{H^1}^2\leq  \kappa_0 \mathcal{F}(t)
       + {\kappa_0} \left|\int  \varepsilon(t) (1-\partial_x^2)\varphi_{c_1}(y_1)dx\right|^2.
    \end{equation}
\end{claim}
The proof of Claim \ref{POSf} follows from classical arguments and we omit it.
See e.g. \cite{We3}, \cite{DiMa} and  \cite{MMcol1}, Appendix D.1.

\medskip

Next, we claim the following control of the variation of $\mathcal{F}(t)$ through time.

\begin{lemma}[Control of the variation of the energy fonctional]\label{varF}
    \begin{equation}\label{varF1}
        |\mathcal{F}'(t)|\leq  K (c_2-1)^{\frac 32}\|\varepsilon(t)\|_{H^1}^2
        + K \|\varepsilon(t)\|_{H^1}\|S(t)\|_{H^1}.
    \end{equation}
    where $K$ is independent of $c_2$ and $K^*$.
\end{lemma}
\begin{proof}
First, we compute $\mathcal{F}'(t)$:
\begin{align*}
\mathcal{F}'(t) & = \int (\partial_t \varepsilon) ((c_1-1)  \varepsilon - c_1   \varepsilon_{xx} - ((v+ \varepsilon)^2 - v^2))\\
& - \int (\partial_t v) \varepsilon^2\\
& + \frac 12 (c_1-c_2)\left\{-c_2 \int \widetilde \alpha''(y_2) ( \varepsilon_x^2+  \varepsilon^2) + \int \widetilde \alpha'(y_2) \partial_t(  \varepsilon_x^2+  \varepsilon^2)\right\} = \mathbf{F}_1+\mathbf{F}_2+\mathbf{F}_3.
\end{align*}
We claim
\begin{align}
& \left|\mathbf{F}_1+\mathbf{F}_2 -  \left\{
\rho'(t) \int  \varepsilon  [(1-\partial_x^2)(\partial_tv+c_1 \partial_xv)]
-\int  \varepsilon^2 (\partial_tv+c_1 \partial_xv)  \right\} \right|\nonumber\\
&\leq K \|\varepsilon(t)\|_{L^2} \|S(t)\|_{H^1}, \label{eq:f1}\\
&  \left|\mathbf{F}_3 - (c_1-c_2)\left\{
\rho'(t) \int \varepsilon  [ (1-\partial_x^2) (\widetilde\alpha'(y_2)  \varphi_{c_1}') ]- \int\varepsilon^2 \widetilde \alpha'(y_2)  \varphi_{c_1}' )\right\}
\right| \nonumber\\
&  \leq K (c_2-1)^{\frac 32} \|\varepsilon\|_{H^1}^2 +K \|\varepsilon(t)\|_{H^1} \|S(t)\|_{H^1}.
\label{eq:f3}
\end{align}
Note that Lemma \ref{varF} follows from \eqref{eq:f1}, \eqref{eq:f3} and \eqref{Tcl2}.
Thus, we only have to prove \eqref{eq:f1}, \eqref{eq:f3} to complete the proof of the lemma.

\medskip

\emph{Proof of \eqref{eq:f1}.}
Using the equation of $\varepsilon(t)$ (i.e. \eqref{eqz}), we find
\begin{align*}
\mathbf{F}_1 &= c_1 \int \varepsilon ((1-\partial_x^2) \partial_t \varepsilon)
- \int (\partial_t \varepsilon) ( \varepsilon + ((v+\varepsilon)^2 - v^2))\\
& = c_1 \left\{ \int (-\partial_x(\varepsilon + (v+\varepsilon)^2 - v^2)) \varepsilon
-\int S(t) \varepsilon + \rho'(t) \int [(1-\partial_x^2)\partial_x(v+\varepsilon)] \varepsilon\right\}\\
&+\int [(1-\partial_x^2)^{-1} \partial_x (\varepsilon + (v+\varepsilon)^2 - v^2)]
(\varepsilon+ ((v+\varepsilon)^2 - v^2))\\&
+ \int [(1-\partial_x^2)^{-1}S(t)](\varepsilon+(v+\varepsilon)^2 - v^2))
-\rho'(t) \int (\partial_x( v+\varepsilon))(\varepsilon+(v+\varepsilon)^2 - v^2)),
\end{align*}
\begin{align*}
\mathbf{F}_1 
& = -c_1 \int \varepsilon^2 (\partial_x v) - c_1 \int S(t)\varepsilon 
+c_1 \rho'(t) \int [(1-\partial_x^2)(\partial_x v)] \varepsilon 
\\ & +\int [(1-\partial_x^2)^{-1}S(t)](\varepsilon+(v+\varepsilon)^2 - v^2))
 -\rho'(t) \int  (\partial_x v)\varepsilon
+ \rho'(t) \int (\partial_x \varepsilon) v^2.
\end{align*}
Thus,
$$\left|\mathbf{F}_1 - \left(-c_1 \int \varepsilon^2 v_x 
+\rho'(t) \int \varepsilon \partial_x(c_1 (1-\partial_x^2)v - v- v^2)\right)\right|\leq K \|\varepsilon(t)\|_{L^2} \|S(t)\|_{H^1}.$$
Using $S=(1-\partial_x^2) \partial_t v + \partial_x (v+v^2)$, we find
$$
\left|\mathbf{F}_1 - \left(-c_1 \int \varepsilon^2 v_x 
+\rho'(t) \int \varepsilon [(1-\partial_x^2)(\partial_t v+c_1 \partial_x v)]\right)\right|
\leq K \|\varepsilon(t)\|_{L^2} \|S(t)\|_{H^1}.
$$
and \eqref{eq:f1} follows from the definition of $\mathbf{F}_2$.

\medskip

\emph{Proof of \eqref{eq:f3}.}
First, from \eqref{Tcl5}, we have
$$\left|\int \widetilde \alpha''(y_2) ( \varepsilon_x^2+  \varepsilon^2)\right|
\leq \|\alpha''\|_{L^\infty} \|\varepsilon\|_{H^1}^2 \leq K (c_2-1)^{\frac 32} \|\varepsilon\|_{H^1}^2.$$
Second,
\begin{equation*}
\frac 12 \int \widetilde \alpha'(y_2) \partial_t(  \varepsilon_x^2+  \varepsilon^2)
= \int \widetilde \alpha'(y_2) (\partial_t(  \varepsilon -\partial_x^2  \varepsilon))
  \varepsilon - \int \widetilde \alpha''(y_2) (\partial_x \varepsilon)
(\partial_t \varepsilon).
\end{equation*}
As before, and using the equation of $\varepsilon$ and Lemma \ref{TRANS},
 we have
$$\left|\int \widetilde \alpha''(y_2) (\partial_x \varepsilon)
(\partial_t \varepsilon)\right| \leq C (c_2-1)^{\frac 32}
 \|\varepsilon\|_{H^1}( \|\varepsilon\|_{H^1}+\|S\|_{H^1}).$$
For the other term, we use the equation of $\varepsilon$,
\begin{align*}
& \int \widetilde \alpha'(y_2) (\partial_t(  \varepsilon -\partial_x^2  \varepsilon))
 \varepsilon \\& 
= \int \widetilde \alpha'(y_2) (-\varepsilon_x \varepsilon
 -(\varepsilon^2)_x \varepsilon - 2(v\varepsilon)_x \varepsilon - S(t) \varepsilon
 + \rho'(t) \varepsilon((1-\partial_x^2)\partial_x (v+\varepsilon) ))\\
 & =\int \widetilde \alpha''(y_2) (\tfrac 12 \varepsilon^2 + \tfrac 23 \varepsilon^3
 +v\varepsilon^2 -\rho'(t)  \tfrac 12 (\varepsilon^2+3\varepsilon_x^2))
 +\rho'(t)\int \widetilde \alpha^{(4)}(y_2) \tfrac 12 \varepsilon^2 \\
 & +\int \widetilde \alpha'(y_2) (-\varepsilon^2 v_x +  \rho'(t) \varepsilon [(1-\partial_x^2)v_x]  - S(t)\varepsilon).
\end{align*}
In the expression above, the term in $\widetilde \alpha''(y_2)$ 
and $\widetilde \alpha^{(4)}(y_2)$ 
is controlled as before, and we get \eqref{eq:f3}.
\end{proof}

\emph{Step 4.} Conclusion of the proof.
By Lemma \ref{Qdir}, and then \eqref{TRANS3}, we have
\begin{align*}
\left|\int \varepsilon(T^*) (1-\partial_x^2)\varphi_{c_1}(y_1)\right|&\leq 
K (c_2-1)^\theta+ K (c_2-1)^{\frac 34} \|\varepsilon(T^*)\|_{H^1}
+ K \|\varepsilon(T^*)\|_{H^1}^2\\&
\leq (K+1) (c_2-1)^\theta,
\end{align*}
for $0<c_2-1<\epsilon_0$ small enough, depending on $K^*$.
Thus, by Claim \ref{Qdir}, we obtain
$$
\|\varepsilon(T^*)\|_{H^1}^2 \leq \kappa_0 \mathcal{F}(T^*)
+ K (c_2-1)^{2\theta}.
$$
Next, integrating \eqref{varF1} on $[T_0,T^*]$,  
by \eqref{TRANS3} and then \eqref{INTkl}, there exists $K_1>0$ 
independent of $K^*$ such that 
\begin{align*}
|\mathcal{F}(T^*)|& \leq
|\mathcal{F}(T_0)| + K (c_2-1)^{\frac 54} T \sup_{t\in [T_0,T^*]} \|\varepsilon(t)\|_{H^1}^2 + K T \sup_{t\in [T_0,T^*]}(\|\varepsilon(t)\|_{H^1}\|S(t)\|_{H^1}) 
\\ & \leq K_1 (c_2-1)^{2\theta} +K(K^*)^2(c_2-1)^{2\theta+\frac 12}  + K_1 K^* (c_2-1)^{2\theta},
\end{align*} 
Thus, for $0<c_2-1<\epsilon_0$ small enough, depending on $K^*$, we obtain
\begin{equation*}
    \|\varepsilon(T^*)\|_{H^1}^2 \leq C (c_2-1)^{2 \theta}\left(2+ K^*\right).
\end{equation*}
Next, fix $K^*$ such that $C (2+K^*)<\frac 12 (K^*)^2$. Then 
$
    \|\varepsilon(T^*)\|_{H^1}^2 \leq \frac 12 (K^*)^2 (c_2-1)^{2 \theta}.
$
This contradict the definition of $T^*$, thus proving that $T^*=T$.
Thus estimate \eqref{INT41} is proved on $[T_0,T]$.
\end{proof}

\subsection{Stability and asymptotic stability for large time}\label{sec:3.2}
In this section, we consider the stability of the $2$-soliton structure after the collision.  For $v\in H^1(\mathbb{R})$, denote 
$$
\|v\|_{H^1_{c_2}}=\left(\int_{\mathbb{R}} \left( (v'(x))^2 + (c_2-1) v^2(x) \right)dx \right)^{\frac 12},$$
which corresponds to the natural norm to study the stability of $\varphi_{c_2}$.
Let $T$ be defined in \eqref{eq:DT}.

\begin{proposition}[Stability of two decoupled solitons]\label{ASYMPTOTIC}
Let $c_1>1$.
There exist  $K>0$,  $\epsilon_0>0$  such that for any $1<c_2<1+\epsilon_0$, the following holds.
Let $u(t)$ be an $H^1$ solution of \eqref{eq:BBM} such that
for some   $\omega>0$, $X_0\geq \frac 12 (c_1-c_2)T$,
\begin{equation}
\label{D25}
\|u(0)-\varphi_{c_1}-\varphi_{c_2}(.+X_0)\|_{H^1}\leq (c_2-1)^{\frac 54 +\omega}.
\end{equation}
Then there exist   $C^1$  functions $\rho_1(t)$, $\rho_2(t)$ defined on $[0,+\infty)$  such that
\begin{enumerate}
\item Stability.
\begin{equation}\label{huit}
\sup_{t\geq 0} \| u(t)-(\varphi_{c_1}(.-\rho_1(t))+\varphi_{c_2}(.-\rho_2(t)))  \|_{H^1_{c_2}}
\le K (c_2-1)^{\frac 54+\omega},
\end{equation}
\begin{equation}\label{suppl}\begin{split}
& 
\forall t\geq 0,\ \tfrac {c_1}2\leq \rho_1'(t)-\rho_2'(t) \leq \tfrac {3c_1}2,
\\ &|\rho_1(t_1)|\leq K (c_2-1)^{\frac 54+\omega},\quad
|\rho_2(t_1)-X_0|\leq K (c_2-1)^{\omega}.
\end{split}\end{equation}
\item Asymptotic stability.
There exist $c_1^+, c_2^+ > 1$  such that
\begin{equation}\label{neuf}
\lim_{t\rightarrow +\infty}\|u(t)-(\varphi_{c_1^+}(x-\rho_1(t))+\varphi_{c_2^+}(x-\rho_2(t)))\|
_{H^1(x> \frac 12 (1+c_2) t)}=0.
\end{equation}
\begin{equation}\label{sept}
|c_1^+-c_1|\leq   K  (c_2-1)^{\frac 54+\omega} ,\quad
 \left| {c_2^+}-{c_2}\right|\le K  (c_2-1)^{1+ \omega+\min(\frac 12, \omega)} .
\end{equation}
\end{enumerate}
\end{proposition}
 
The proof of Proposition \ref{ASYMPTOTIC} is essentially the same as the one of Theorem 1.1 in \cite{DiMa}, combined with Theorem 2 in \cite{Ma}
(see also the previous works \cite{Mi2}, \cite{Di}).
See also Proposition 4 in \cite{MMas2}.

\begin{proof}[Sketch of the proof of Proposition \ref{ASYMPTOTIC}]
Let $u(t)$ satisfying the assumption \eqref{D25} of the proposition.

\medskip

\emph{1. Stability.} For $D_0>2$ to be chosen later, we define
\begin{align*}
t^*(\varepsilon)=\sup& \Big\{t\geq 0 \ | \ \forall t'\in [0,t), \ \exists y_1,y_2\in \R \ | \ y_1-y_2 >\tfrac 14 (c_1-c_2) T\\Ê& \text{ and } 
  \|u(t')-\varphi_{c_1}(.-y_1)-\varphi_{c_2}(.-y_2) \|_{H^1_{c_2}} \leq 
D_0 (c_2-1)^{\frac 54 + \omega}\Big\}.
\end{align*}
Note that $t_0>0$ is well-defined by continuity of $t\mapsto u(t)$ in $H^1$. We assume for the sake of contradiction that $t^*<+\infty$.

First, we decompose the solution using modulation theory.

\begin{claim}\label{cl:md}
        For $0<c_2-1<\epsilon_0$ small enough, there exist $\rho_1(t)$, $\rho_2(t)$, $\bar c_1(t)$, $\bar c_2(t)$, defined on $[0,t^*]$ such that
        \begin{equation}\label{eq:ta}
                \eta(t,x)=u(t,x)-R_1(t,x)-R_2(t,x), \quad \text{where}\quad
                R_j(t,x)=\varphi_{\bar c_j(t)} (x-\rho_j(t)),
        \end{equation}
         satisfies 
        \begin{align}   
                \forall t\geq 0, \quad &        \int ((1-\partial_x^2) R_j(t))\eta(t)
                        =\int ((1-\partial_x^2) \partial_x R_j(t))\eta(t)=0, \ (j=1,2), 
                        \label{eq:o1}\\
                &       \|\eta(t)\|_{H^1}+|\bar c_1(t)-c_1|+(c_2-1)^{-\frac 14} |\bar c_2(t)-c_2|\leq CD_0(c_2-1)^{\frac 34 + \omega},
                        \label{eq:o2}\\
                &       |\rho'_j(t)-c_j|\leq \frac 12 (c_2-1),\quad \rho_1(t)-\rho_2(t)\geq \frac 14 (c_1-c_2) T + \frac 12 c_1 t,
                        \label{eq:o3}\\
                &       \|\eta(0)\|_{H^1} + |\bar c_1(0)-c_1|
 + (c_2-1)^{-\frac 14} |\bar c_2(0)-c_2|\leq C (c_2-1)^{\frac 54 + \omega},
                        \label{eq:o4}\\
                &       |\rho_1(0)|+ (c_2-1)^{\frac 54} |\rho_2(0)+X_0|\leq C(c_2-1)^{\frac 54 + \omega}.               \label{eq:o5}
        \end{align}
\end{claim}

\begin{proof}
For the proof of this claim, we refer to  proof of Claim 2.1 in \cite{MMas2} and proof of Proposition~2.1 in \cite{DiMa}.
We only observe that since for $c>1$, $\varphi_c(x)=(c-1)Q\left(\sqrt{\frac {c-1} c} x\right)$, we have
$$
\frac {d \varphi_c} {dc}(x)= Q\left(\sqrt{\frac {c-1} c} x\right)+ \frac {(c-1)^{\frac 12}}{c^{\frac 32}} \frac x 2 Q'\left(\sqrt{\frac {c-1} c} x\right) 
= \frac 1{c-1} \left(\varphi_c + \frac x{2c} \partial_x \varphi_c\right).
$$
and thus, setting
$P_j(t,x)={\frac {d \varphi_c} {dc}} _{|c=\bar c_j(t)} (x-\rho_j(t)),$
we check that    $\eta(t,x)$ satisfies the following
\begin{equation}\label{eq:ea}\begin{split}
&(1-\partial_x^2) \partial_t \eta + \partial_x \eta + \sum_{j=1,2} \bar c_j'(t)(1-\partial_x^2) P_j - \sum_{j=1,2} (\rho_j'(t)-\bar c_j(t)) (1-\partial_x^2) \partial_x R_j\\&+\partial_x(2R_1R_2 + 2 \eta (R_1+R_2) + \eta^2)=0.
\end{split}\end{equation}
\end{proof}
Note that for $t=0$, the estimates \eqref{eq:o4} and \eqref{eq:o5} are independent of $D_0$. In the rest of the proof the objective is to prove estimates on $\eta(t)$ at $t=t^*$, independent of $D_0$ by using conservation laws and monotonicity properties on localized versions of these conservation laws, thus contradicting the definition of $t^*$ for $D_0$ large enough.

Indeed, we claim the following.

\begin{lemma}\label{le:in}
        There exist $D_0,\varepsilon_0>0$ such that for $0<c_2-1<\epsilon_0$,
        \begin{equation}\label{eq:ts}
                \sup_{t\in [0,t^*]}\|u(t)-\varphi_{c_1}(.-\rho_1(t))-\varphi_{c_2}(.-\rho_2(t))\|_{H^1_{c_2}} \leq \tfrac 12  D_0 (c_2-1)^{\frac 54 + \omega}.
        \end{equation}
\end{lemma}

Assuming $t^*<+\infty$, by Lemma \ref{le:in} and the continuity of $u(t)$ in $H^1$ we obtain a contradiction. Therefore, we only have to prove Lemma \ref{le:in}.

\medskip

With respect to the classical proof of stability of one soliton by Weinstein \cite{We3}, the main additionnal argument of the proof of estimate \eqref{eq:ts} is the following monotonicity property.
Let
\begin{equation}\label{eq:ph}\begin{split}
        & \psi(x)=\frac 2 \pi \arctan(\exp(x/\kappa)), \quad \text{so that }
        \lim_{-\infty} \psi=0,\ \lim_{\infty} \psi=1,\\
        & \forall x\in \R,\quad \psi(-x)=1-\psi(x),\quad
 \psi'(x)=\frac 1{\pi\kappa \cosh(x/\kappa)},
 \quad  |\psi'''(x)|\leq \frac 1{\kappa^2} |\psi'(x)|,\\
        & {\mathcal{N}_1}(t)=\frac12\int (u^2(t,x)+u_x^2(t,x))\psi(x-m(t)) dx,
\quad   m(t)=\frac 12 (\rho_1(t)+\rho_2(t)).
\end{split}\end{equation}

\begin{claim}\label{cl:mo}
For $0<c_2-1<\epsilon_0$ small enough and $\kappa$ large enough,
\begin{equation*}
        \forall t\in [0,t^*],\quad
                {\mathcal{N}_1}(t)-{\mathcal{N}_1}(0)\leq C (c_2-1)^{10+\omega}.
\end{equation*}
\end{claim}

The proof of Claim \ref{cl:mo} is based on the following identity
($g$ any $C^1$ function):
\begin{equation*}
        \frac d{dt} \int (u^2(t)+u_x^2(t)) g(x) = -\int u^2(t) g'(x) - \frac 23 \int u^3 g' + 2 \int u  [(1-\partial_x^2)^{-1} (u+u^2)] g'(x),
\end{equation*}
and the arguments of the proof of Lemma 2.1 in \cite{DiMa}. Note that such monotonicity results in the context of the BBM equation were first introduced in \cite{Mi} and \cite{Di}. We also refer to Appendix \ref{se:xC} in the present paper for similar monotonicity arguments.

One can actually obtain an estimate of the type 
$\exp(- (c_2-1)^{-\gamma})$, for some $\gamma>0$, but the estimate in Claim \ref{cl:mo} will be sufficient for our purposes.

\medskip

The rest of the proof of Lemma \ref{le:in} is   similar to the proof of Lemma 2.2 in \cite{MMas2} and Theorem 1.1 in \cite{DiMa}.
Let 
\begin{equation}\label{eq:dg}\begin{split}
&       g(t)=  \int \left[\eta_x^2(t,x) + (c(t,x)-1) \eta^2(t,x)\right]dx ,
 \\Ê&c(t,x)=\bar c_2(t) + (\bar c_1(t)-\bar c_2(t)) \psi(x-m(t)).
\end{split}\end{equation}

First, using the two invariant quantities $N(u(t))$ and $E(u(t))$, one proves
\begin{equation}\label{eq:qd}
\forall t\in [0,t^*],\quad
|\bar c_1(t)-\bar c_1(0)|+ (c_2-1)^{\frac 32} |\bar c_2(t)-\bar c_2(0)|\leq
C (g(t)+g(0)+(c_2-1)^{10+\omega}),
\end{equation}

Second, using the monotonicity property (see Claim \ref{cl:mo}) and the related quantity
\begin{equation}\label{eq:dF} 
\mathcal{F}(t)= \int \left[ c(t,x) u_x^2(t,x) + (c(t,x)-1) u^2(t,x) - \tfrac 23 u^3(t,x)\right] dx,  \end{equation}
we claim the following 
\begin{equation}\label{eq:hc}
\forall t\in [0,t^*],\quad
g(t)\leq C g(0) + C(c_2-1)^{10+\omega}.
\end{equation}
Note   that $\mathcal{F}(t)$ is a functional of $u(t)$ which is locally around each soliton $R_1$, $R_2$ adapted to the proof the stability of one soliton,

\medskip

Finally, combining \eqref{eq:qd}, \eqref{eq:hc} and \eqref{eq:o4}, \eqref{eq:o5}, we obtain $\forall t\in [0,t^*],$
\begin{equation}\label{eq:fi}\begin{split}
&
\|\eta(t)\|_{H^1_{c_2}}^2+
g(t)+|\bar c_1(t)-\bar c_1(0)|+ (c_2-1)^{\frac 32} |\bar c_2(t)-\bar c_2(0)|\leq
C (c_2-1)^{\frac 52+2\omega}, \\
& |\bar c_1(t)-c_1|\leq C (c_2-1)^{\frac 54 + \omega},
\quad |\bar c_2(t)-c_2|\leq C (c_2-1)^{\omega+\min(\frac 12, \omega)},\\
& \|u(t)-\varphi_{c_1}(.-\rho_1(t))-\varphi_{c_2}(.-\rho_2(t))\|_{H^1_{c_2}}\leq
C_0 (c_2-1)^{\frac 54+\omega},
\end{split}\end{equation}
where $C_0>0$ is independent of $D_0$. Choosing now
$D_0=4 C_0$, we obtain Lemma \ref{le:in}.

\medskip

\emph{2. Asymptotic stability.} For this part, we refer to section 4 of \cite{DiMa}. Recall that the main argument of the proof is the following rigidity result, combined with monotonicity arguments, such as Claim \ref{cl:mo}.

\begin{proposition}\label{pr:li} Let $c_0>1$. There exists $\alpha_0=\alpha_0(c_0)>0$ such that if $u(t)$ is an $H^1$ solution of \eqref{eq:BBM} satisfying
\begin{equation}
        \|u(0)-\varphi_{c_0}\|_{H^1} \leq \alpha_0,
\end{equation}
and 
\begin{equation}
        \forall \delta>0,\ \exists B_\delta>0, \text{ s.t. } \forall t\in \R,\quad
        \int_{|x|> B_\delta} (u^2+u_x^2)(t,x+y(t)) dx <\delta,
\end{equation}
for some function $y(t)$, then there exists $x_1\in \R$, $c_1>0$, such that
\begin{equation}
        \forall t,x\in \R,\quad
        u(t,x)=\varphi_{c_1}(x-x_1-c_1t).
\end{equation}
\end{proposition}

Since the linear Liouville theorem which underlies Proposition \ref{pr:li} has been extended to any $c_0>1$ in \cite{Ma} (see Theorem 2 in \cite{Ma}), 
Theorem 4.1 in \cite{DiMa} applies for any $c_0>1$ (note that \cite{MW} and 
\cite{Di,Mi2} could not exclude countably many exceptions). 
The convergence of $\bar c_j(t)$ to some limit value $c_j^+$ is obtained as in \cite{DiMa} using monotonicity results such as Claim \ref{cl:mo}.
Finally, estimate \eqref{sept} follows from passing to the limit as $t\to +\infty$ in estimate \eqref{eq:fi}.
\end{proof}

\section{Proof of Theorem \ref{th:1}}

First, we recall the following existence and uniqueness result of asymptotic $2$-soliton solutions for the BBM equation. Recall that $T$ is defined in \eqref{eq:DT}.

\begin{proposition}
\label{prop:5.1}
    Let $c_1>1$ and $1<c_2<1+\epsilon_0$, for $\epsilon_0>0$ small enough.
    \begin{enumerate}
    \item  Let $x_1,x_2\in \mathbb{R}$.
      There exists a unique $H^1$ solution  
    $u(t)=u_{c_1,c_2,x_1,x_2}(t)$
     of \eqref{eq:BBM} such that 
    \begin{equation}\label{eq:sol}
        \lim_{t\to -\infty} \|u(t)-\varphi_{c_1}(.-c_1t-x_1)-\varphi_{c_2}(.-c_2t-x_2)\|_{H^1}=0.
    \end{equation}
    Moreover,  
    for all  $t\leq - \frac {T} {32}$,
        \begin{equation}\label{2SOL1}
            \|u(t)-\varphi_{c_1}(.-c_1t-x_1)-\varphi_{c_2}(.-c_2t-x_2)\|_{H^1}
            \leq K e^{\frac 14{\sqrt{c_2-1}(c_1-1)t}}.
        \end{equation}
        \item  
        If $w(t)$ is an $H^1$ solution of \eqref{eq:BBM} satisfying
        \begin{equation}\label{2SOL2}
            \lim_{t\to -\infty}\|w(t)- \varphi_{c_1}(.-\rho_1(t))-\varphi_{c_2}(.-\rho_2(t))\|_{H^1}=0,
        \end{equation}
        for some $\rho_1(t)$, $\rho_2(t)$,  then there exist $x_1,x_2\in \mathbb{R}$ such that 
        $w(t) \equiv u_{c_1,c_2,x_1,x_2}(t)$.
    \end{enumerate}
\end{proposition}

Proposition \ref{prop:5.1} is essentially the same as Theorem 1.3 in \cite{DiMa}. Recall that such result was first proved for the gKdV equations in \cite{Ma2}, refining techniques introduced in \cite{MM1}, \cite{MMT}.
The second statement of the proposition slightly improves the original result in \cite{DiMa} and is easily proved by the same techniques (see \cite{MMcol1}, Appendix D.2 for the case of the gKdV equation).

\medskip

Next, we claim the following lemma, concerning the variation of $c_1$ and $c_2$ the speeds of the two solitons after the collision.

\begin{lemma}\label{le:ii}
Let $c_1>1$.  There exists $\epsilon_0=\epsilon_0(c_1)>0$ such that the following holds. Let $1<c_2<1+\epsilon_0$.
Suppose that  $u(t)$ is a solution of \eqref{eq:BBM}  satisfying, for some $\rho_j(t)$ ($j=1,2$)
\begin{equation}\label{eq:mi}
        \lim_{t\to -\infty} \|u(t)-\varphi_{c_1}(.-\rho_1(t)) - \varphi_{c_2}(.-\rho_2(t))\|_{H^1}=0,
\end{equation}
\begin{equation}\label{eq:pi}
        \lim_{t\to +\infty} \|u(t)-\varphi_{c_1^+}(.-\rho_1(t)) - \varphi_{c_2^+}(.-\rho_2(t))- w^+(t)\|_{H^1}=0,
\end{equation}
where $|c_j^+-c_j|\leq \epsilon_0 |c_j-1|$ and 
\begin{equation}\label{eq:wi}
\lim_{t\to +\infty} \|w^+(t)\|_{H^1(x>\frac 12 (c_2+1) t)}=0,\quad 
\limsup_{t\to +\infty}\|w^+(t)\|_{H^1}\leq  \epsilon_0 |c_2-1|.
\end{equation}
Then, there exists $C=C(c_1)>0$ such that
\begin{equation}\label{eq:vc}\begin{split}
&\frac 1 C \limsup_{t\to +\infty} \|w^+(t)\|_{H^1_{c_2}}^2
\leq c_1^+-c_1 \leq C \liminf_{t\to +\infty} \|w^+(t)\|_{H^1_{c_2}}^2,\\
&\frac 1 C  (c_2-1)^{-\frac 12}  \limsup_{t\to +\infty} \|w^+(t)\|_{H^1}^2
\leq c_2-c_2^+ \leq C  (c_2-1)^{-\frac 12}  \liminf_{t\to +\infty} \|w^+(t)\|_{H^1}^2.
\end{split}\end{equation}
\end{lemma}
\begin{remark}
This kind of property was first observed for the quartic gKdV equation in \cite{MMcol1}.
See also the general discussion in \cite{MMcol2}.
\end{remark}
\begin{proof}
By \eqref{eq:cv}, \eqref{eq:mi}, \eqref{eq:pi} and \eqref{eq:wi}, we have for $t$ large
\begin{align}
& N(u(0))=N(\varphi_{c_1})+ N(\varphi_{c_2}) = N(\varphi_{c_1^+})+ N(\varphi_{c_2^+}) + N(w^+(t)) + o(1),\label{eq:N2}\\
& E(u(0))=E(\varphi_{c_1})+ E(\varphi_{c_2}) = E(\varphi_{c_1^+})+ E(\varphi_{c_2^+}) + E(w^+(t)) + o(1).\label{eq:E2}
\end{align}
Let ($j=1,2$)
\begin{equation}\label{eq:da}
\bar a_j=\frac { E(\varphi_{c_j^+}) -E(\varphi_{c_j})}{ N(\varphi_{c_j^+}) -N(\varphi_{c_j})},\quad
\text{so that}\quad 
|\bar a_j-c_j|\leq C |c_j^+-c_j|.
\end{equation}
Indeed, by \eqref{eq:A2},
$$\frac { E(\varphi_{c_j^+}) -E(\varphi_{c_j})}{ N(\varphi_{c_j^+}) -N(\varphi_{c_j})}
={\frac {\frac d{dc} E(\varphi_c)}{\frac d{dc} N(\varphi_c)}} { _{|c=c_j}}+O(|c_j^+-c_j|)=c_j+O(|c_j^+-c_j|).$$
In particular,
$\bar a_2-1\geq (c_2-1)-|c_2-\bar a_2|\geq (c_2-1)-C |c_2^+-c_2|
\geq (1-C\epsilon_0)(c_2-1)\geq \frac 12 (c_2-1)$, for $\epsilon_0$ small.

Considering $\bar a_2 \times$\eqref{eq:N2}$-$\eqref{eq:E2} and then $\bar a_1 \times$\eqref{eq:N2}$-$\eqref{eq:E2}, we find, for $t$ large,
\begin{equation}\label{eq:l1}
[ E(\varphi_{c_1^+})-\bar a_2 N(\varphi_{c_1^+})]
- [E(\varphi_{c_1})- \bar a_2 N(\varphi_{c_1})]
=  \bar a_2 N(w^+(t)) - E(w^+(t))  + o(1),
\end{equation}
\begin{equation}\label{eq:l2}
[\bar a_1 N(\varphi_{c_2})- E(\varphi_{c_2})]
- [\bar a_1 N(\varphi_{c_2^+})- E(\varphi_{c_2^+})]
=\bar a_1 N(w^+(t)) - E(w^+(t)) + o(1).
\end{equation}
Note that $\int |w^+|^3 \leq C \|w^+\|_{H^1} \int (w^+)^2\leq C \epsilon_0 |c_2-1|\int (w^+)^2$ so that
$\bar a_2N(w^+)-E(w^+)> \frac 14 ((c_2-1) \int (w^+)^2 + \int (w_x^+)^2)$.

Now, let $\beta_1=\frac d{dc} N(\varphi_c)_{|c=c_1}>0$. 
By \eqref{eq:A2}, we have
$(\frac d{dc} E(\varphi_c)- \bar a_2 \frac d{dc} N(\varphi_c))_{|c=c_1}=(c_1-\bar a_2)\frac d{dc} N(\varphi_c)_{|c=c_1}$, and so 
$\frac 12 (c_1-1) \beta_1<(\frac d{dc} E(\varphi_c)- \bar a_2 \frac d{dc} N(\varphi_c))_{|c=c_1} <(c_1-1) \beta_1$.
Thus, from \eqref{eq:l1}, we obtain, for $t$ large,
\begin{equation*}\begin{split}
& c_1^+-c_1\geq C \left[(c_2-1)\int (w^+(t))^2 
+c_2 \int (w_x^+(t))^2\right]+o(1)\geq C  \|w^+(t)\|_{H^1_{c_2}}^2+o(1),\\
& c_1^+-c_1\leq C'  \|w^+(t)\|_{H^1_{c_2}}^2+o(1).
\end{split}\end{equation*}
Similarly, using 
$\beta_2=(c_2-1)^{-\frac 12} \frac d{dc} N(\varphi_c)_{|c=c_2}>K>0$, independent of $c_2$, it follows  from \eqref{eq:l2} that for $t$ large,
\begin{equation*}\begin{split}
&
c_2-c_2^+ 
\geq \frac C{(c_2-1)^{\frac 12}} \left[(c_1-1) \int (w^+(t))^2  + c_1 \int (w_x^+(t))^2\right]+o(1)\geq \frac {C \|w^+(t)\|_{H^1}^2} {(c_2-1)^{\frac 12}}+o(1),\\
& c_2-c_2^+ 
\leq  \frac {C' \|w^+(t)\|_{H^1}^2} {(c_2-1)^{\frac 12}}+o(1).
\end{split}\end{equation*}
Estimates \eqref{eq:vc} follow.
\end{proof}

\begin{proof}[Proof of Theorem \ref{th:1}]
Let $c_1>1$ and let $\epsilon_0=\epsilon_0(c_1)$ small enough, so that the results of Sections 2 and 3 apply.
Let $1<c_2<1+\epsilon_0$.
Let $T$ be defined by \eqref{eq:DT}.
Let $\widetilde u(t)$ be the unique solution of \eqref{eq:BBM} such that
$$
\lim_{t\to -\infty} \|\widetilde u(t)-\varphi_{c_1}(.-c_1t)-\varphi_{c_2}(.-c_2t)\|_{H^1}=0.
$$

\emph{1. Behavior at $-T$.} Proposition \ref{prop:5.1} implies that
        \begin{equation}\label{eq:so}
        \forall t\leq - \frac {T} {32},\quad 
            \|\widetilde u(t)-\varphi_{c_1}(.-c_1t)-\varphi_{c_2}(.-c_2t)\|_{H^1}            \leq K e^{\frac 14{\sqrt{c_2-1}(c_1-1)t}}.
        \end{equation}
Let $\Delta_1$, $\Delta_2$ be defined in Proposition \ref{cor:1} and let
$$
T^-=T+ \frac 12 \frac {\Delta_1-\Delta_2}{c_1-c_2}.
$$
Since $|\Delta_1|+|\Delta_2|\leq C=C(c_1)$, 
and $c_1-c_2>c_1-1-\varepsilon_0 \geq \frac 12 (c_1-1)$, 
we have $-T^-<-\frac 1{32} T$, for $c_2$ small and so
\begin{equation}\label{eq:sq}
            \|\widetilde u(-T^-)-\varphi_{c_1}(.+c_1 T^-)-\varphi_{c_2}(.+c_2T^-)\|_{H^1}            \leq K e^{-\frac 14{\sqrt{c_2-1}(c_1-1)T^-}}\leq (c_2-1)^{10},
        \end{equation}
for $\epsilon_0$ small enough.

Let 
\begin{equation}\label{eq:ut}
u(t,x)=\widetilde u(t+T-T^-,x+\tfrac 12 \Delta_1+c_1(T-T^-)).
\end{equation}
Then, $u(t)$ is solution of \eqref{eq:BBM} and satisfies
\begin{equation}\label{eq:sp}
            \|  u(-T)-\varphi_{c_1}(.+c_1 T+\tfrac 12 \Delta_1)-\varphi_{c_2}(.+c_2T+\tfrac 12 \Delta_2)\|_{H^1}            \leq   (c_2-1)^{10}.
        \end{equation}
In what follows, we work with $u(t)$. It is easily checked that the results obtained for $u(t)$ imply the desired results on $\widetilde u(t)$.

\medskip

\emph{2. Behavior at $T$.} 
By Proposition \ref{cor:1} and \eqref{eq:sp}, we have
$$
\|u(-T)-v(-T)\|_{H^1}\leq K (c_2-1)^{\frac {13}4}.
$$
By Proposition \ref{prop:approx} and the above estimate,
we can apply Proposition \ref{prop:I} with $\theta=3-\frac 12 -\frac 1{100}=\frac 52 - \frac 1{100}.$ It follows that there exists $\rho(t)$ such that
$$
\forall t\in [-T,T],\quad 
\|u(t)-v(t,.-\rho(t))\|_{H^1}+|\rho'(t)| \leq C (c_2-1)^{\frac 52- \frac 1{100}}.
$$
In particular, for $r=\rho(T)$, $|r|\leq C (c_2-1)^{2-\frac 1{50}}$, we have
$\|u(T)-v(t,.-r)\|_{H^1}\leq C (c_2-1)^{\frac 52- \frac 1{100}}$, and using
Proposition \ref{cor:1}, we obtain
\begin{equation}\label{eq:sr}
\|u(T)-\{\varphi_{c_1}(.-r_1)+\varphi_{c_2}(.-r_2)-2 D(\varphi_{c_2}^2)'(.-r_2)\}\|_{H^1}
\leq C (c_2-1)^{\frac 52- \frac 1{100}},
\end{equation}
where 
$r_1=c_1 T +\tfrac 12 \Delta_1+r$ and 
$r_2=c_2 T +\tfrac 12 \Delta_2+r$, so that $$\frac 12 (c_1-c_2)T  \leq r_1-r_2
\leq \frac 32 (c_1-c_2)T.$$
Moreover, since $\|(\varphi_{c_2}^2)'\|_{H^1}\leq  C (c_2-1)^{\frac 94},$ we also obtain
\begin{equation}\label{eq:su}
\|u(T)-\{\varphi_{c_1}(.-r_1)+\varphi_{c_2}(.-r_2)\}\|_{H^1}
\leq C (c_2-1)^{\frac 94}.
\end{equation}
In what follows, \eqref{eq:su} will serve us to prove that $u(t)$ is close to the sum of two solitons for $t>T$, whereas \eqref{eq:sr} will allow us to prove that $u(t)$ is not a pure $2$-soliton solution at $+\infty$.

\medskip

\emph{3. Behavior as $t\to +\infty$}. We use Proposition \ref{ASYMPTOTIC} with $\omega=1$.  It follows from \eqref{sept}, \eqref{huit} and \eqref{neuf} that there exists $\rho_1(t)$, $\rho_2(t)$, $c_1^+$, $c_2^+$ such that
\begin{equation}\label{eq:xx}\begin{split}
&c_1^+=\lim_{t\to +\infty} \bar c_1(t),\quad
c_2^+=\lim_{t\to +\infty} \bar c_2(t),\\
& |c_1^+-c_1|\leq C (c_2-1)^{\frac 94},\quad
|c_2^+-c_2|\leq C (c_2-1)^{\frac 32},\quad \text{and}\\
& w^+(t,x)=u(t,x)-\{\varphi_{c_1^+}(x-\rho_1(t))+\varphi_{c_2^+}(x-\rho_2(t))\} \quad \text{satisfies}  \\
&\sup_{t\geq T} \|w^+(t)\|_{H^1_{c_2}} \leq C (c_2-1)^{\frac 94},\quad
\lim_{t\to +\infty} \|w^+(t)\|_{H^1(x> \frac 12 (1+c_2) t)}=0.  
\end{split}\end{equation}
 From Lemma \ref{le:ii}, we have the following more precise estimates on $c_j^+-c_j$
 $$
 0\leq c_1^+-c_1 \leq C (c_2-1)^{\frac 92},\quad
 0\leq c_2-c_2^+ \leq C (c_2-1)^{4}.
 $$

\emph{4. Lower bound on $w^+(t)$ for $t>T$ large.}
Consider the decomposition of $u(T+.,.+r_1)$ defined in Claim~\ref{cl:md}, i.e.
the center of mass $\rho_1(t)$, $\rho_2(t)$ defined before and
$\bar c_1(t)$, $\bar c_2(t)$, $\eta(t)$ for $t>T$.
In particular,
\begin{equation}\label{eq:tt}\begin{split}
        &\sup_{t\geq T}\|\eta(t)\|_{H^1_{c_2}}\leq C (c_2-1)^{\frac 94},\quad 
        |\rho_1(T)-r_1|+|\bar c_1(T)-c_1| \leq C (c_2-1)^{\frac 94},\\
        & |\rho_2(T)-r_2|\leq C (c_2-1),\quad 
        |\bar c_2(T)-c_2|\leq C(c_2-1)^2.
\end{split}\end{equation}

First, as a consequence of \eqref{eq:sr}, we claim the following lower bound at $t=T$: for $C_0>0$, independent of $c_2$,
\begin{equation}\label{eq:aT}
        \int_{x<\rho_2(T)+\frac 14 T} \eta^2(T,x) dx \geq C_0 (c_2-1)^{\frac 92}.
\end{equation}
Proof of \eqref{eq:aT}. Replacing $u(T,x)=\varphi_{\bar c_1(T)}(x-\rho_1(t)) + \varphi_{\bar c_2(T)}(x-\rho_2(t))+\eta(t,x)$ in \eqref{eq:sr}, we find
\begin{equation*}\begin{split}
& \|[\varphi_{\bar c_1(T)}(.-\rho_1(T))-\varphi_{c_1}(.-r_1)]
+[\varphi_{\bar c_2(T)}(.-\rho_2(T))-\varphi_{c_2}(.-r_2)]
\\ & + \eta(T)+2D  (\varphi_{c_2}^2)'(.-r_2)\|_{H^1} 
 \leq C (c_2-1)^{\frac 52-\frac 1{100}}.
\end{split}\end{equation*}
By the decay properties of $\varphi_{c_1}$ and $r_1-r_2\geq \frac 12 (c_1-c_2) T$,
we obtain
\begin{equation*} 
  \|[\varphi_{\bar c_2(T)}(.-\rho_2(T))-\varphi_{c_2}(.-r_2)]
 + \eta(T)+2D  (\varphi_{c_2}^2)'(.-r_2)\|_{L^2(x<\rho_2(T)+\frac 14 T)} 
 \leq C (c_2-1)^{\frac 52-\frac 1{100}}.
\end{equation*}
Assuming, to the contrary,  that for any $\alpha>0$ there exist $c_2$ arbitrary close to $1$ such that
$$\| \eta(T) \|_{L^2(x<\rho_2(T)+\frac 14 T)} 
 \leq \alpha (c_2-1)^{\frac 9 4}.
 $$
Then
\begin{equation*} 
  \|[\varphi_{\bar c_2(T)}(.-\rho_2(T))-\varphi_{c_2}(.-r_2)]
  +2D  (\varphi_{c_2}^2)'(.-r_2)\|_{L^2(x<\rho_2(T)+\frac 14 T)} 
 \leq 2\alpha (c_2-1)^{\frac 9 4}.
\end{equation*}
By scaling and translation, and decay of $Q$, we obtain
\begin{equation*} 
  \|[\bar Q-Q + 2 D(c_2-1)^{\frac 32}  (Q^2)'\|_{L^2} \leq 2\alpha (c_2-1)^{\frac 32},
\end{equation*}
where $\bar Q(x)= \lambda Q\left( \mu x - \xi\right),$ and
$$
\lambda=\frac {\bar c_2(T)-1}{c_2-1}, \
\mu=\sqrt{\frac {\bar c_2(T)-1}{\bar c_2(T)}} \sqrt{\frac {c_2}{c_2-1}},\
\xi=\sqrt{\frac {\bar c_2(T)-1}{\bar c_2(T)}}(\rho_2(T)-r_2).
$$
Note that by \eqref{eq:tt}, we have
$$
|\lambda-1|\leq C (c_2-1)^{\frac 12},\quad
|\xi|\leq C (c_2-1)^{\frac 32}.
$$
Expanding $\bar Q$ in $\lambda-1$, $\mu$ and $\xi$, and using parity properties, we find
$$
\| \xi Q'  + 2 D(c_2-1)^{\frac 32}  (Q^2)'\|_{L^2} \leq 3\alpha (c_2-1)^{\frac 32},
$$
so that for some constant $\bar \xi\in \R$,
$$
\| \bar \xi Q'  + 2 D  (Q^2)'\|_{L^2} \leq 4\alpha.
$$
But since $D\neq 0$, whatever is the value of $\bar \xi$, this is not true for $\alpha>0$ small enough. This contradiction proves \eqref{eq:aT}.

\medskip

Now, we finish the proof of the lower bound by proving the following.
There exists $K_0>0$ such that 
\begin{equation}\label{eq:lw}
        \liminf_{t\to +\infty} \|w^+(t)\|_{H^1_{c_2}} \geq K_0 (c_2-1)^{\frac {11}4}.
\end{equation}
Indeed, note that \eqref{eq:lw} combined with Lemma \ref{le:ii} prove the lower bounds in \eqref{eq:th:1:2}. Thus, we are now reduced to prove \eqref{eq:lw}.

\medskip

\noindent Proof of \eqref{eq:lw}.
We argue by contradiction. Assume that for any $\alpha>0$, there exist
arbitrarily large $T_0$ and $c_2$ arbitrarily close to $1$ such that
\begin{equation}\label{eq:ct}
\|w^+(T_0)\|_{H^1_{c_2}} \leq \alpha (c_2-1)^{\frac {11}4}.
\end{equation}
By \eqref{eq:xx}, we can choose $T_0>T$ large enough so that
\begin{equation}\label{eq:yy}
        \|\eta(T_0)\|_{H_{c_2}^1(x< m(T_0)+\frac{T_0}4)} \leq 2 \alpha (c_2-1)^{\frac {11}4}.
\end{equation}

We consider the same functions $\psi(x)$ and $m(t)$ as in \eqref{eq:ph}, with $\kappa= \sqrt{\frac {c_1+7}{c_1-1}}$.
Let 
$$
a_2= \frac {E(\varphi_{\bar c_2(T_0)}) - E(\varphi_{\bar c_2(T)})}
{N(\varphi_{\bar c_2(T_0)}) - N(\varphi_{\bar c_2(T)})}.
$$
We set 
\begin{equation}\label{eq:vq}
\begin{split}
\mathcal{G}(t)  &       =a_2\int (u_x^2+u^2)(t,x) (1-\psi(x-m(t))) dx 
- \int (u^2+\tfrac 23 u^3)(t,x) (1-\psi(x-m(t))) dx \\
                                &       =a_2 N(u(t)) - E(u(t)) - (a_2 \mathcal{N}_1(t) - \mathcal{E}_1(t)),
\end{split}\end{equation}
where 
$$
\mathcal{N}_1(t) = \frac 12 \int (u_x^2+u^2)(t,x)  \psi(x-m(t)) dx,\quad 
\mathcal{E}_1(t) = \frac 12 \int (u^2+\tfrac 23 u^3)(t,x) \psi(x-m(t)) dx.
$$
We claim the following results on $\mathcal{G}(t)$.

\begin{lemma}\label{le:mo}
For $0<c_2-1<\epsilon_0$ small enough,
$$
\mathcal{G}(T)-\mathcal{G}(T_0) \leq C (c_2-1)^{10}.
$$
\end{lemma}

\begin{lemma}\label{le:qd}
For $0<c_2-1<\epsilon_0$ small enough,
$$
\mathcal{G}(T_0)        -\mathcal{G}(T) = \frac 12 (\mathcal{H}(T_0)-\mathcal{H}(T))
+O(\sup_{[T,T_0]}\|\eta\|_{H^1}^3)+O((c_2-1)^{10}), 
$$
where
\begin{equation}\begin{split}
\mathcal{H}(t)&
=\int \left((a_2-1) \eta^2 + a_2 \eta_x^2 - 2 R_2 \eta^2\right)(t,x) (1-\psi(x-m(t))) dx \\
& \geq  \sigma_0 \int \left[(c_2-1) \eta^2 + \eta_x^2 \right](t,x) (1-\psi(x-m(t))) dx.
\end{split}\end{equation}
for some $\sigma_0>0$ independent of $c_2$.
\end{lemma}
See proofs of Lemmas \ref{le:mo} and \ref{le:qd} in Appendix \ref{se:xC}.

\medskip

Combining Lemmas \ref{le:mo} and \ref{le:qd}, we find
\begin{equation*}\begin{split}
        &\int \left[(c_2-1) \eta^2 + \eta_x^2 \right](T,x) (1-\psi(x-m(T))) dx\\
        & \leq C \mathcal{H}(T_0) + O(\|\eta\|_{H^1}^3)+O((c_2-1)^{10})
        +\mathcal{G}(T)-\mathcal{G}(T_0)\\
        & \leq C \alpha (c_2-1)^{\frac {11} 2} + C (c_2-1)^{\frac {21}4},
\end{split}\end{equation*}
by using \eqref{eq:yy} and \eqref{eq:xx}.
But this estimate contradicts \eqref{eq:aT} for $\alpha>0$ small enough and $0<c_2-1<\epsilon_0>0$ small enough.
\end{proof}

\appendix

\section{Appendix -- Identities on $\varphi_c$}\label{sec:A}

\begin{claim}[Identities on $Q$]
  \label{cl:qint}
  \begin{align}
  & \int Q=\int Q^2,\quad 
\int Q^3=\frac{6}{5}\int Q^2, \quad \int Q'^2=\frac15\int Q^2,\label{eq:qint1}\\
&\int x^2 Q^3= \frac 65 \int x^2 Q^2 -\frac 35 \int Q^2, \quad 
\int x^2 (Q')^2 = \frac 15 \int x^2 Q^2 + \frac 25 \int Q^2,\label{eq:qint2}\\
& \int Q^2=6, \quad \int x^2 Q^2 = 2 \pi^2 - 12.
\label{eq:qint3}
\end{align}
\end{claim}
\begin{proof}
The values of $\int Q^2$ and $\int x^2 Q^2$ are easily computed
using the fact that
$\int_0^\infty \frac{x}{e^x+1}dx=\frac{\pi^2}{12}$.
The relations \eqref{eq:qint1}, \eqref{eq:qint2} are obtained by using
\eqref{eq:Q}.
\end{proof}

\begin{claim}\label{cl:A.2}
For all $c>1$,
\begin{equation}\label{eq:A2}\begin{split}
&
\int \varphi_c^2 = (c-1)^{\frac 32} c^{\frac 12} \int Q^2, \quad
\int \varphi_c^3=\frac 65 (c-1) \int \varphi_c^2, \quad
\int (\varphi_c')^2= \frac 15 \left(\frac {c-1}{c}\right) \int \varphi_c^2,\\
&
E(\varphi_c)=\frac 12 \left(1+\frac 45 (c-1)\right)  \int \varphi_c^2=
 \frac 12 (c-1)^{\frac 32}c^{\frac 12} \left(1+\frac 45 (c-1)\right) \int Q^2,\\
&
N(\varphi_c)=\frac 12  \left( \frac 15 \left(\frac {c-1}c\right) +1\right) \int \varphi_c^2 =
\frac 1{2} (c-1)^{\frac 32}c^{-\frac 12} \left( \frac 15 (c-1)+c\right) \int Q^2,\\
& E(\varphi_c)-c N(\varphi_c)
=-\frac 15 (c-1) \int \varphi_c^2
=-\frac 15 (c-1)^{\frac 52} c^{\frac 12} \int Q^2,\quad
\frac d{dc} E(\varphi_c)= c \frac d{dc} N(\varphi_c)>0.
\end{split}\end{equation}
\end{claim}
\begin{proof}
Recall that $\varphi_c(x)=(c-1) Q\left(\sqrt{\frac {c-1}{c}} x\right)$. Thus,
$$
\int \varphi_c^2 = (c-1)^2 \sqrt{\frac {c}{c-1}} \int Q^2= (c-1)^{\frac 32} c^{\frac 12} \int Q^2.
$$
Next, we have
$$c\varphi_c''-(c-1)\varphi_c + \varphi_c^2=0,\quad
c (\varphi_c')^2 - (c-1) \varphi_c^2 + \frac 23 \varphi_c^3=0.$$

Thus,
$$
-c \int (\varphi_c')^2 - (c-1) \int \varphi_c^2 + \int \varphi_c^3 =0,\quad
c \int (\varphi_c')^2 - (c-1) \int \varphi_c^2 +\frac 23 \int \varphi_c^3 =0.
$$
Combining the above identities, we find $(c-1) \int \varphi_c^2 = \frac 56 \int \varphi_c^3$
and $c \int (\varphi_c')^2 = \frac 15 (c-1) \int \varphi_c^2$.
The formulas concerning $E(\varphi_c)$ and $N(\varphi_c)$ then follow from direct computations.

For the last relation, multiply by $\frac d{dc} \varphi_c$ the equation of $\varphi_c$ written under the form
$-c (1-\partial_x^2) \varphi_c + (\varphi_c+\varphi_c^2)=0$, and integrate on $\mathbb{R}$.

\end{proof} 

\section{Appendix -- Proof of Proposition \ref{prop:decomp}}\label{sec:B}

\begin{claim}
\label{cl:hderiv}
Let $g$ be a $C^3$-function and $h(t,x)=g(y)=g(x-\alpha(\ys))$.
Then
  \begin{align*}
& \pd_th=-\mu_\sigma\beta(\ys)g'(y),\\
& \pd_xh=(1-\beta(\ys))g'(y),\\
& \pd_x^2h=(1-\beta(\ys))^2g''(y)-\beta'(\ys)g'(y),\\
& \pd_x\pd_th=-\mu_\sigma(1-\beta(\ys))\beta(\ys)g''(y)
-\mu_\sigma\beta'(\ys)g'(y),\\
& \pd_x^3h=(1-\beta(\ys))^3g'''(y)
-3(1-\beta(\ys))\beta'(\ys)g''(y)-\beta''(\ys)g'(y),\\
& \pd_x^2\pd_th=\mu_\sigma\{
-(1-\beta(\ys))^2\beta(\ys)g'''(y)
+3\beta(\ys)\beta'(\ys)g''(y)-2\beta'(\ys)g''(y)-
\beta''(\ys)g'(y)\}.
\end{align*}
\end{claim}
\begin{proof}
Differentiating $h(t,x)=g(x-\alpha(\ys))$ with respect to $t$
and $x$ respectively, we have
\begin{align*}
& \pd_th(t,x)=\frac{\pd y}{\pd t}g'(y)=-\alpha'(\ys)
\frac{\pd \ys}{\pd t}g'(y)=-\mu_\sigma\beta(\ys)g'(y),
\\
& \pd_xh(t,x)=\frac{\pd y}{\pd x}g'(y)=(1-\beta(\ys))g'(y).
  \end{align*}
Here we use $\pd_t\ys=\mu_\sigma$ and $\pd_ty=-\mu_\sigma\beta(\ys)$.
We compute $\pd_x^2h(t,x)$, $\pd_x\pd_th(t,x)$ and $\pd_x^3h(t,x)$
in the same way.
By using the first and the third formulas, we compute
\begin{align*}
& \pd_t\pd_x^2h  = (1-\beta(\ys))^2\pd_tg''(y)
+g''(y)\pd_t(1-\beta(\ys))^2 
-\beta'(\ys)\pd_t g'(y)-g'(y)\pd_t\beta'(\ys)
\\&= 
-\mu_\sigma(1-\beta(\ys))^2\beta(\ys)g'''(y)
+3\mu_\sigma\beta(\ys)\beta'(\ys)g''(y)
-2\mu_\sigma\beta'(\ys)g''(y)-\mu_\sigma\beta''(\ys)g'(y).
\end{align*}
\end{proof}

We follow the notation introduced in \eqref{eq:df1}--\eqref{eq:df4} and we also set
\begin{equation*}\begin{split}
& S(z):=(1-\lambda\pd_x^2) \pd_t z+\pd_x(\pd_x^2z-z+z^2)=S_{KdV}(z)+S_{BBM}(z),\\
&  S_{KdV}(z):=\pd_t z+\pd_x(\pd_x^2z-z+z^2),\quad
S_{BBM}(z):=-\lambda\pd_t\pd_x^2z.
\end{split}\end{equation*}
Then
\begin{equation}\label{eq:sion}
S(z(t,x))=S(Q(y))+S(\wqs(\ys))
+\delta S(w(t,x))+S_{int}(t,x),
\end{equation}
where 
\begin{equation}\label{eq:sion2}\begin{split}
& \mathcal{L}=-\partial_x^2 + 1- 2 Q(y),\\
&
\delta S(w):=\delta S_{KdV}(w)+S_{BBM}(w),\quad 
\delta S_{kdV}(w):=\pd_t w-\pd_x\mathcal{L}w,
 \\
&S_{int}(t,x)=\pd_x\left\{w^2(t,x)
+2\wqs(\ys)(Q(y)+w(t,x))\right\}.
\end{split}\end{equation}
Since $\wqs(\ys)$ is a solution to \eqref{eq:BBM2},
we have $S(\wqs)=0.$

\begin{claim}
  \label{cl:SKdV1}
Let $A$ and $q$ be $C^3$-functions. Then 
\begin{align*}
& \quad  \delta S_{KdV}(A(y)q(\ys))  \\    &=
q(\ys)\bigl\{-(LA)'(y)+\beta(\ys)(-3A''-2AQ+(1-\mu_\sigma)A)'(y)-\beta'(\ys)(3A'')(y)
\\ & +\beta^2(\ys)(3A''')(y)+(\beta^2)'(\ys)(3A''/2)(y)
-\beta''(\ys)A'(y)-\beta^3(\ys)A'''(y)\bigr\}
\\ &+
q'(\ys)\bigl\{3A''(y)+2A(y)Q(y)+(\mu_\sigma-1)A(y)-\beta(\ys)(6A'')(y)
-\beta'(\ys)(3A')(y)\\&+\beta^2(\ys)(3A'')(y)\bigr\}
  +q''(\ys)\{3(1-\beta(\ys))A'(y)\}+q'''(\ys)A(y).
\end{align*}
\end{claim}
\begin{proof}
In the proof, we omit the variable $y$ of $A(y)$.
Using Claim \ref{cl:hderiv}, we compute
$$\pd_t (A(y)q(\ys))=-\mu_\sigma\beta(\ys)A'q(\ys)+\mu_\sigma A
q'(\ys),$$
and
\begin{align*}
 -\pd_x\mathcal{L}(A(y)q(\ys))&= 
\pd_x\{(\pd_x^2A-A+2AQ)q(\ys)+2(\pd_xA)q'(\ys)+Aq''(\ys)\}
\\&\{\pd_x(\pd_x^2A-A+2AQ)\}q(\ys)+(\pd_x^2A-A+2AQ)q'(\ys)
\\ &+2(\pd_x^2A)q'(\ys)+3(\pd_xA)q''(\ys)+Aq'''(\ys)
\\&=
q(\ys)\bigl\{(1-\beta(\ys))^3A'''
-3(1-\beta(\ys))\beta'(\ys)A''-\beta''(\ys)A'
\\ & \qquad\qquad
-(1-\beta(\ys))A'+2(1-\beta(\ys))(AQ)'\bigr\}
\\ &
+q'(\ys)\{3(1-\beta(\ys))^2A''-3\beta'(\ys)A'-A+2AQ\}
\\ & +q''(\ys)\{3(1-\beta(\ys))A'\}+q'''(\ys)A.
\end{align*}
Combining the above, we conclude Claim \ref{cl:SKdV1}.
\end{proof}
\begin{claim}
\label{cl:SBBM0}
  Let $q$ and $A$ be $C^3$-functions.
Then
  \begin{align*}
 &\quad S_{BBM}(A(y)q(\ys))= \lambda\mu_\sigma q(\ys)
\{\beta(\ys)A'''(y)+\beta'(\ys)(2A''(y))\}
\\ & +
\lambda\mu_\sigma q(\ys)\{\beta^2(\ys)(-2A''')(y)
+(\beta^2)'(\ys)(-3A''/2)(y)+\beta''(\ys)A'(y)+\beta^3(\ys)A'''(y)\}
\\ & +
\lambda\mu_\sigma q'(\ys)
\{-A''(y)+\beta(\ys)(4A'')(y)+\beta'(\ys)(3A')(y)+\beta^2(\ys)(-3A'')(y)\}
\\ & +
\lambda\mu_\sigma q''(\ys)\{-2A'(y)+\beta(\ys)(3A')(y)\}
+\lambda\mu_\sigma q'''(\ys)(-A)(y).
  \end{align*}
\end{claim}
\begin{proof}
  \begin{align*}
\pd_x^2\pd_t(A(y)q(\ys)) &=
(\pd_x^2\pd_tA(y))q(\ys)+2(\pd_x\pd_tA(y))\pd_xq(\ys)
+(\pd_tA(y))\pd_x^2q(\ys),\\
& + (\pd_x^2A(y))\pd_tq(\ys)+2(\pd_xA(y))\pd_x\pd_tq(\ys)
+A(y)\pd_x^2\pd_tq(\ys).
\end{align*}
By Claim \ref{cl:hderiv},
\begin{align*}
&(\pd_x^2\pd_tA(y))q(\ys)+2(\pd_x\pd_tA(y))\pd_xq(\ys)
+(\pd_tA(y))\pd_x^2q(\ys)\\
   &= \mu_\sigma \{-(1-\beta(\ys))^2\beta(\ys)A'''
+3\beta(\ys)\beta'(\ys)A''-2\beta'(\ys)A''
-\beta''(\ys)A'\}q(\ys)
\\ &+
2\mu_\sigma\{-(1-\beta(\ys))\beta(\ys)A''-\beta'(\ys)A'\}
q'(\ys)-\mu_\sigma\beta(\ys)A'q''(\ys),
\end{align*}
and
\begin{align*}
 & (\pd_x^2A(y))\pd_tq(\ys)+2(\pd_xA(y))\pd_x\pd_tq(\ys)
+A(y)\pd_x^2\pd_tq(\ys)\\&=
\mu_\sigma\{(1-\beta(\ys))^2A''-\beta'(\ys)A'\}q'(\ys)
  +2\mu_\sigma(1-\beta(\ys))A'q''(\ys)
+\mu_\sigma Aq'''(\ys).
\end{align*}
Combining the above, we obtain
  \begin{align*}
  \pd_t\pd_x^2(A(y)q(\ys))
 &= \mu_\sigma\{-(1-\beta(\ys))^2\beta(\ys)A'''
+\frac{3}{2}(\beta^2)'(\ys)A''
-2\beta'(\ys)A''-\beta''(\ys)A'\}q(\ys)
\\ &+
\mu_\sigma\{(3\beta(\ys)^2-4\beta(\ys)+1)A''
-3\beta'(\ys)A'\}q'(\ys)
\\ &+
\mu_\sigma(2-3\beta(\ys))A'q''(\ys)
+\mu_\sigma Aq'''(\ys).
  \end{align*}
Thus  Claim \ref{cl:SBBM0} is proved.
\end{proof}
\begin{claim}\label{cl:musigma}
\begin{equation}
  \label{eq:musigma}\begin{split}
& \mu_\sigma=\frac{1-\sigma}{1-\lambda\sigma}=1+(\lambda-1)\sigma
\sum_{j=0}^\infty(\lambda\sigma)^j,\\
& 
\frac 1 {\theta_\sigma}=\frac 1{1-\lambda} - \frac \lambda {1-\lambda} \sigma,\quad
\theta_\sigma=\frac {1-\lambda}{1-\lambda\sigma}=(1-\lambda)\sum_{j=0}^\infty(\lambda\sigma)^j.
\end{split} 
\end{equation}
\end{claim}

\begin{claim}
\label{cl:beta}
Let 
$$ \beta=a_{1,0}\wqs+a_{1,1}\sigma\wqs+a_{2,0}\wqs^2+a_{3,0}\wqs^3+
a_{2,1}\sigma \wqs^2+a_{1,2}\sigma^2 \wqs.$$
Then,
  \begin{align*}
& \beta'=a_{1,0}\wqs'+a_{1,1}\sigma\wqs'+a_{2,0}(\wqs^2)'+a_{3,0}(\wqs^3)'+
a_{2,1}\sigma (\wqs^2)'+a_{1,2}\sigma^2 (\wqs)',\\
& \beta''=\sigma\wqs a_{1,0}+\wqs^2\left(-\frac{a_{1,0}}{1-\lambda}\right)
+\sigma^2\wqs a_{1,1}
+\sigma\wqs^2\left(\frac{\lambda a_{1,0}}{1-\lambda}
-\frac{a_{1,1}}{1-\lambda}+4a_{2,0}\right)\\
& \qquad +\wqs^3\left(-\frac{10a_{2,0}}{3(1-\lambda)}\right)
+\sigma^3 O(\wqs),\\
& \beta^2=a_{1,0}^2\wqs^2+a_{1,1}^2\sigma^2\wqs^2
+2(a_{1,0}a_{1,1}\sigma\wqs^2
+a_{2,0}a_{1,0}\wqs^3)+\sigma^3 O(\wqs),\\
& (\beta^2)'=a_{1,0}^2(\wqs^2)'+a_{1,1}^2\sigma^2(\wqs^2)'
+2\left\{a_{1,0}a_{1,1}\sigma(\wqs^2)'
+a_{2,0}a_{1,0}(\wqs^3)'\right\}+\sigma^3 O(\wqs).
  \end{align*}
\end{claim}
\begin{proof}
The proof follows by elementary calculations from \eqref{eq:qsigma} and \eqref{eq:musigma}.
\end{proof}

In the next lemmas, we expand the various terms in \eqref{eq:sion}.

\begin{lemma}
\label{lem:SQ}
\begin{equation}
  \label{eq:SQ}\begin{split}
S(Q(y))= &\sum_{(k,l)\in \Sigma_0} \sigma^l\left( \wqs^k(\ys)a_{k,l}\{(\lambda-3)Q''-Q^2\}'(y)
+  (\wqs^k)'(\ys)a_{k,l}(2\lambda-3)Q''(y)\right)
\\ &+\sum_{(k,l)\in \Sigma_0}\sigma^l\left(\wqs^k(\ys)F_{k,l}^I(y)
+ (\wqs^k)'(\ys)G_{k,l}^I(y)\right)+\sigma^3 O(\wqs(\ys)),
\end{split}
\end{equation}
where
\begin{align*}
& F_{1,0}^I=0,\quad
G_{1,0}^I=0,\quad  F_{1,1}^I= \lambda(\lambda-1)a_{1,0}Q''',
\quad  
G_{1,1}^I= a_{1,0}\left(2\lambda(\lambda-1)Q''\right),
\\ &
F_{2,0}^I= (3-2\lambda)a_{1,0}^2Q'''+a_{1,0}Q',
\quad G_{2,0}^I= \frac{3}{2}(1-\lambda)a_{1,0}^2Q'',
\end{align*}
and for all $(k,l)\in \Sigma_0$ such that $k+l=3$, $F_{k,l}^I\in \mathcal{Y}$ is odd, $G_{k,l}^I\in \mathcal{Y}$ is even and depend only on $a_{k',l'}$ for $1\leq k'+l'\leq 2.$
\end{lemma}
\begin{lemma}
  \label{lem:dSKdVw}
  \begin{align*}
\delta S_{KdV}(w)&= \sum_{(k,l)\in\Sigma_0}\sigma^l
\left( \wqs^k(\ys)(-LA_{k,l})'(y)+ (\wqs^k)'(\ys) ((-LB_{k,l})'+3A_{k,l}''+2QA_{k,l})(y)\right)
\\&
+\sum_{(k,l)\in\Sigma_0}\sigma^l
\left(\wqs^k(\ys)F_{k,l}^{II}(y)+(\wqs^k)'(\ys)G_{k,l}^{II}(y)\right)+\sigma^3 O(\wqs(\ys)),
  \end{align*}
where
\begin{align*}
&  F_{1,0}^{II}=0,\quad
G_{1,0}^{II}=0,\quad 
 F_{1,1}^{II}= 3A_{1,0}'+3B_{1,0}''+2QB_{1,0},\quad
 G_{1,1}^{II}= \lambda A_{1,0}+3B_{1,0}',\\
& F_{2,0}^{II}= -a_{1,0}(3A_{1,0}''+2QA_{1,0})'
-\frac{1}{1-\lambda}(3A_{1,0}'+3B_{1,0}''+2QB_{1,0})\\
& G_{2,0}^{II}=-a_{1,0}\left(\frac{9}{2}A_{1,0}'+\frac{3}{2}B_{1,0}''+QB_{1,0}\right)'
-\frac{1}{1-\lambda}(A_{1,0}+3B_{1,0}'),
\end{align*}
and for $(k,l)\in \Sigma_0$ such that $k+l=3$, $F_{k,l}^{II}$, $G_{k,l}^{II}$ depend on 
$A_{k',l'}$, $B_{k',l'}$ for $1\leq k'+l'\leq 2$.
Moreover, if $A_{k',l'}$ are even and $B_{k',l'}$ are odd then 
$F_{k,l}^{II}$ are odd and $G_{k,l}^{II}$ are even.
\end{lemma}
\begin{lemma}
  \label{lem:SBBMw}
 \begin{align*}
S_{BBM}(w)&= \sum_{(k,l)\in \Sigma_0 }\sigma^l
  (\wqs^k)'(\ys) (-\lambda A_{k,l}'')(y)
\\ & + \sum_{(k,l)\in \Sigma_0 }\sigma^l
\left(\wqs^k(\ys)F_{k,l}^{III}(y)+(\wqs^k)'(\ys)G_{k,l}^{III}(y)\right)
+\sigma^3 O(\wqs(\ys)),
  \end{align*}
where
\begin{align*}
&  F_{1,0}^{III}=0,\quad G_{1,0}^{III}=0,\\
&  F_{1,1}^{III}=-2\lambda A_{1,0}'-\lambda B_{1,0}'' \quad
G_{1,1}^{III}= \lambda(1-\lambda)A_{1,0}''
-\lambda A_{1,0}-2\lambda B_{1,0}',\\
&   F_{2,0}^{III}=\lambda a_{1,0}A_{1,0}'''+\frac{2\lambda}{1-\lambda}A_{1,0}'+
\frac{\lambda}{1-\lambda}B_{1,0}'',\\
&
G_{2,0}^{III}= 3\lambda a_{1,0}A_{1,0}''
+\frac{\lambda}{1-\lambda}A_{1,0}+\frac{\lambda a_{1,0}}{2}B_{1,0}'''
+\frac{2\lambda}{1-\lambda}B_{1,0}',
\end{align*}
and for $(k,l)\in \Sigma_0$ such that $k+l=3$,  $F_{k,l}^{III}$, $G_{k,l}^{III}$ depend on 
$A_{k',l'}$, $B_{k',l'}$ for $1\leq k'+l'\leq 2$.
Moreover, if $A_{k',l'}$ are even and $B_{k',l'}$ are odd then 
$F_{k,l}^{III}$ are odd and $G_{k,l}^{III}$ are even. 
\end{lemma}

\begin{lemma}
  \label{lem:Sint}
  \begin{align*}
S_{int}(w)=\sum_{(k,l)\in \Sigma_0}\sigma^l
\left(\wqs^k(\ys)F_{k,l}^{int}(y)+(\wqs^k)'(\ys)G_{k,l}^{int}(y)\right)+\sigma^3 O(\wqs),
  \end{align*}
where
\begin{align*}
&  F_{1,0}^{int}=2Q',\quad G_{1,0}^{int}=2Q,\quad
F_{1,1}^{int}=G_{1,1}^{int}=0,\\
& F_{2,0}^{int}=(A_{1,0}^2)'+2A_{1,0}'-2a_{1,0}Q',
\quad
G_{2,0}^{int}=(A_{1,0}B_{1,0})'+A_{1,0}^2+2A_{1,0}+B_{1,0}',
\end{align*}
and for $(k,l)\in \Sigma_0$ such that $k+l=3$,   $F_{k,l}^{int}$, $G_{k,l}^{int}$ depend on 
$A_{k',l'}$, $B_{k',l'}$ for $1\leq k'+l'\leq 2$.
Moreover, if $A_{k',l'}$ are even and $B_{k',l'}$ are odd then 
$F_{k,l}^{int}$ are odd and $G_{k,l}^{int}$ are even.
\end{lemma}

Putting together Lemmas \ref{lem:SQ}--\ref{lem:Sint}, we obtain Proposition \ref{prop:decomp}, in
particular, the explicit expressions of $F_{k,l}$ and $G_{k,l}$ for $1\leq k+l\leq 2$.

\begin{remark}\label{rk:4}
In the proof of Lemmas \ref{lem:SQ}--\ref{lem:Sint},
we mainly focus on the computations of $F_{k,l}$, $G_{k,l}$
for $1\leq k+l\leq 2$, since the explicit expressions of these terms is fundamental in proving the inelasticity of the collision. The exact expressions of $F_{k,l}$ and $G_{k,l}$ for $k+l=3$ is not needed in what follows. The expression of the rest terms i.e. $\sigma^3 O(\wqs)$ is also not useful.
We will give a simple bound of the rest term in Section \ref{sec:2.3} after choosing the functions $A_{k,l}$, $B_{k,l}$, $(k,l)\in \Sigma_0$. For more details on the structure of $F_{k,l}$ and $G_{k,l}$ for $k+l=3$, we refer to \cite{MMcol1}, where similar computations are performed completely in the case of the quartic gKdV equation, see proof of Lemmas A.1--A.4 in \cite{MMcol1}.
\end{remark} 

\begin{proof}[Proof of Lemma \ref{lem:SQ}.]
By  Claim \ref{cl:SKdV1} with $A(y)=Q(y)$ and $q=1$,
\begin{align*}
  S_{KdV}(Q)&=\delta S_{KdV}(Q)-\pd_x(Q^2)
\\ &= (Q''-Q+Q^2)'+\beta(\ys)\{-3Q''-Q^2+(1-\mu_\sigma)Q\}'
-\beta'(\ys)3Q''
\\ & +\beta^2(\ys)(3Q''')+(\beta^2)'(\ys)(3Q''/2)
-\beta''(\ys)Q'-\beta^3(\ys)Q'''.
\end{align*}
Using $Q''=Q-Q^2$ and \eqref{eq:musigma}, we obtain
  \begin{align*}
  S_{KdV}(Q)&= \beta(\ys)(-3Q''-Q^2)'+\beta'(\ys)(-3Q'')
\\ &+\beta^2(\ys)(3Q''')+(\beta^2)'(\ys)(3Q''/2)-\sigma\beta(\ys)(\lambda-1)Q'\\
& -\beta''(\ys) Q' - \beta^3(\ys) Q'''-(\lambda-1)\lambda\sigma^2 \beta(\ys) Q'+\sigma^3 O(\wqs).
  \end{align*}
Next, Claim \ref{cl:SBBM0} and \eqref{eq:musigma} implies
\begin{align*}
S_{BBM}(Q)&= \lambda\mu_\sigma\bigl\{\beta(\ys)Q'''+\beta'(\ys)(2Q'')
+\beta^2(\ys)(-2Q''')+(\beta^2)'(\ys)(-3Q''/2)
\\ & +\beta''(\ys)Q' +\beta^3(\ys)Q'''\bigr\}
\\&= 
 \beta(\ys)(\lambda Q''')+\beta'(\ys)(2\lambda Q'')+ 
\beta^2(\ys)(-2\lambda Q''')+(\beta^2)'(\ys)(-3\lambda Q''/2)
\\ & 
+\sigma\beta(\ys)\{\lambda(\lambda-1) Q'''\}
+\sigma\beta'(\ys)\{2\lambda(\lambda-1)Q''\}
+\beta''(\ys) (\lambda Q') + \beta^3(\ys) \lambda Q'''\\&
+ \sigma \beta^2(\ys)\lambda (\lambda-1)(-2Q''')
+\sigma (\beta^2)'(\ys)\lambda(\lambda-1) (-3Q''/2) +\sigma \beta''(\ys)\lambda(\lambda-1)Q'\\&
+ \sigma^2 \beta(\ys) \lambda^2(\lambda-1) Q''' + \sigma^2 \beta'(\ys) \lambda^2(\lambda-1)(2Q'')+\sigma^3 O(\wqs).
\end{align*}
Combining the above, we obtain
\begin{align*}
S(Q)&= \beta(\ys)\{(\lambda-3)Q''-Q^2\}'+\beta'(\ys)(2\lambda-3)Q''
\\ & +
\beta^2(\ys)(3-2\lambda)Q'''+(\beta^2)'(\ys) (1-\lambda)(3Q''/2)
+\beta''(\ys)(\lambda-1)Q'
\\ & +\sigma\beta(\ys)(\lambda-1)\{\lambda Q''-Q\}'
+\sigma\beta'(\ys)\{2\lambda(\lambda-1)Q''\} 
+ \beta^3(\ys) (\lambda-1) Q'''\\&
+ \sigma \beta^2(\ys)\lambda (\lambda-1)(-2Q''')
+\sigma (\beta^2)'(\ys)\lambda(\lambda-1) (-3Q''/2) +\sigma \beta''(\ys)\lambda(\lambda-1)Q'\\&
+ \sigma^2 \beta(\ys) \lambda(\lambda-1) (\lambda Q'' -Q)' + \sigma^2 \beta'(\ys) \lambda^2(\lambda-1)(2Q'')+\sigma^3 O(\wqs).
\end{align*}

Hence using Claim \ref{cl:beta}, we can obtain
\begin{align*}
S(Q)&=\wqs(\ys)a_{1,0}\{(\lambda-3)Q''-Q^2\}'+
\wqs'(\ys)a_{1,0}(2\lambda-3)Q'' \\ & +
\wqs^2(\ys)\left(
a_{2,0}\{(\lambda-3)Q''-Q^2\}'
+(3-2\lambda)a_{1,0}^2Q'''+a_{1,0}Q'\right)
\\ &+
(\wqs^2)'(\ys)\left(
a_{2,0}(2\lambda-3)Q''+\tfrac{3}{2}(1-\lambda)a_{1,0}^2Q''\right)
\\ &+
\sigma\wqs(\ys)\left(a_{1,1}\{(\lambda-3)Q''-Q^2\}'
+\lambda(\lambda-1)a_{1,0}Q'''\right)
\\ & +
\sigma\wqs'(\ys)
\left(a_{1,1}(2\lambda-3)Q''
+2\lambda(\lambda-1)a_{1,0}Q''\right)\\
& + \sum_{k+l=3} \sigma^l\left( \wqs^k(\ys)a_{k,l}\{(\lambda-3)Q''-Q^2\}'(y)
+  (\wqs^k)'(\ys)a_{k,l}(2\lambda-3)Q''(y)\right)
\\
&+\sum_{k+l=3} \sigma^l\left(\wqs^k(\ys) F^I_{k,l}+(\wqs^k)'(\ys)G^I_{k,l} \right) +\sigma^3 O(\wqs),
\end{align*}
where for all $k+l=3$, $F_{k,l}^I\in \mathcal{Y}$ and $G_{k,l}^I\in \mathcal{Y}$ are as in the statement of the Lemma.
\end{proof}

In the proof of Lemmas \ref{lem:dSKdVw}--\ref{lem:Sint}, we compute explicitly only up to the order
of $c^l (\wqs^k)'$ for $1\leq k+l\leq 2$. See Remark \ref{rk:4}.

\begin{proof}[Proof of Lemma \ref{lem:dSKdVw}.]
  \begin{equation*}
\delta S_{KdV}(w)=\sum_{(k,l)\in \Sigma_0}\sigma^l
\left(\delta S_{KdV}(A_{k,l}(y)\wqs^k(\ys))
+\delta S_{KdV}(B_{k,l}(y)(\wqs^k)'(\ys))\right).
  \end{equation*}
First, we compute $\delta S_{KdV}(A_{1,0}(y)\wqs(\ys))$.
By Claim \ref{cl:SKdV1} and the definition of $\beta$, we have
\begin{align*}
&\quad \delta S_{KdV}(A_{1,0}(y)\wqs(\ys))=\\
&
\wqs(\ys)\bigl\{-(LA_{1,0})'+a_{1,0}\wqs(\ys)(-3A_{1,0}''-2A_{1,0}Q)'
-a_{1,0}\wqs'(\ys)(3A_{1,0}'')\bigr\}
\\ & + \wqs'(\ys)\left\{
3A_{1,0}''+2A_{1,0}Q-a_{1,0}\wqs(\ys)(6A_{1,0}'')\right\}
\\ & +\wqs''(\ys)(3A_{1,0}')+\wqs'''(\ys)A_{1,0}
+\sigma\wqs'(\ys)(\lambda-1)A_{1,0}+\sigma^2 O(\wqs(\ys)).
  \end{align*}
  Note in particular that we have used $(\wqs')^2 = \sigma^2  O(\wqs(\ys))$ from \eqref{eq:qsigma}.
Next, by \eqref{eq:qsigma} and \eqref{eq:musigma}, we have
\begin{align*}
& \wqs''(\ys)(3A_{1,0}')+\wqs'''(\ys)A_{1,0}
\\&=\left(\sigma\wqs(\ys)-\frac{1}{1-\lambda}\wqs^2(\ys)\right)(3A_{1,0}')
+\left(\sigma\wqs'(\ys)-\frac{1}{1-\lambda}(\wqs^2)'(\ys)\right)A_{1,0}
+\sigma^2 O(\wqs(\ys)).
\end{align*}
Thus,
\begin{equation}
  \label{eq:*1}
\begin{split}
&\quad \delta S_{KdV}(A_{1,0}(y)\wqs(\ys)) =  \wqs(\ys)(-LA_{1,0})'+\wqs'(\ys)(3A_{1,0}''+2A_{1,0}Q)
\\ & +
\wqs^2(\ys)\left(a_{1,0}(-3A_{1,0}''-2A_{1,0}Q)'-\frac{3A_{1,0}'}{1-\lambda}
\right)
+(\wqs^2)'(\ys)
\left(-\frac{9}{2}a_{1,0}A_{1,0}''-\frac{A_{1,0}}{1-\lambda}\right)
\\ & +
\sigma\wqs(\ys)(3A_{1,0}')+\sigma\wqs'(\ys)(\lambda A_{1,0})+\sigma^2 O(\wqs(\ys)).
\end{split}  
\end{equation}

Now, we compute $\delta S(B_{1,0}(y)\wqs'(\ys))$ in a similar way:
\begin{align*}
\delta S(B_{1,0}(y)\wqs'(\ys))&=  \wqs'(\ys)
\bigl\{-(LB_{1,0})'+a_{1,0}\wqs(\ys)(-3B_{1,0}''-2B_{1,0}Q)'\bigr\}
\\ & +\wqs''(\ys)(3B_{1,0}''+2B_{1,0}Q)+\wqs'''(\ys)(3B_{1,0}') + \sigma^2 O(\wqs(\ys)).
\\   &=\wqs'(\ys)(-LB_{1,0})'
+\wqs^2(\ys)\left(-\frac{1}{1-\lambda}(3B_{1,0}''+2B_{1,0}Q)\right)
\\ & +(\wqs^2)'(\ys)\left(-\frac{3}{1-\lambda}B_{1,0}'
-a_{1,0}\left(\frac{3}{2}B_{1,0}''+B_{1,0}Q\right)'\right)
\\ & +\sigma\wqs(3B_{1,0}''+2B_{1,0}Q)+\sigma\wqs'(3B_{1,0}')+ \sigma^2 O(\wqs(\ys)).
\end{align*}

Similarly, we have for all $(k,l)$  with $2\leq k+l\leq 3$,
\begin{align*}
& \delta S_{KdV}(\sigma^l\wqs^k(\ys)A_{k,l}(y))=\sigma^l\wqs^k(\ys)(-LA_{k,l})'+\sigma\wqs'(\ys)(3A_{k,l}''+2A_{k,l}Q)
+ \sigma^2 O(\wqs(\ys)),\\
&  \delta S_{KdV}(\sigma^l(\wqs^k)'(\ys)B_{k,l}(y))=\sigma^l(\wqs^k)'(\ys)(-LB_{k,l})'+ \sigma^2 O(\wqs(\ys)).
\end{align*}
Combining the above, we obtain Lemma \ref{lem:dSKdVw}.
\end{proof}

\begin{proof}[Proof of Lemma \ref{lem:SBBMw}.]
By definition,
    \begin{equation*}
S_{BBM}(w)=\sum_{(k,l)\in \Sigma_0}\sigma^l
\left(S_{BBM}(A_{k,l}(y)\wqs^k(\ys))
+S_{BBM}(B_{k,l}(y)(\wqs^k)'(\ys))\right).
  \end{equation*}
First, we compute $S_{BBM}(A_{1,0}(y)\wqs(\ys))$.
As in the proof of Lemma \ref{lem:dSKdVw}, it follows from
Claim \ref{cl:SBBM0} and then \eqref{eq:qsigma}, \eqref{eq:musigma} that
\begin{align*}
& \quad S_{BBM}(A_{1,0}(y)\wqs(\ys))\\ & = \lambda\mu_\sigma\wqs(\ys)\{\beta(\ys)A_{1,0}'''+\beta'(\ys)(2A_{1,0}'')\}
  +\lambda\mu_\sigma\wqs'(\ys)\{-A_{1,0}''+\beta(\ys)(4A_{1,0}'')\}
\\ & + \lambda\mu_\sigma\wqs''(\ys)(-2A_{1,0}')
+\lambda\mu_\sigma\wqs'''(\ys)(-A_{1,0})+ \sigma^2 O(\wqs(\ys))
\\ & =\lambda\wqs(\ys)\{a_{1,0}\wqs(\ys)A_{1,0}'''
+a_{1,0}\wqs'(\ys)(2A_{1,0}'')\}
\\ &  
+\lambda\{1+(\lambda-1)\sigma\}\wqs'(\ys)(-A_{1,0}'')
+(\wqs\wqs')(\ys)(4\lambda a_{1,0}A_{1,0}'')
\\ &  
+ \lambda\left(\sigma\wqs(\ys)-\frac{1}{1-\lambda}\wqs^2(\ys)\right)
(-2A_{1,0}') \\ &+\lambda\left(
\sigma\wqs'(\ys)-\frac{1}{1-\lambda}(\wqs^2)'(\ys)\right)(-A_{1,0})+ \sigma^2 O(\wqs(\ys))
\\& =  
\wqs'(\ys)(-\lambda A_{1,0}'')+\sigma\wqs(\ys)(-2\lambda A_{1,0}')
+\sigma\wqs'(\ys)\{\lambda(1-\lambda) A_{1,0}''-\lambda A_{1,0}\}
\\ &   + 
\wqs^2(\ys)\left(\lambda a_{1,0}A_{1,0}'''+\frac{2\lambda}{1-\lambda}
A_{1,0}'\right)
+(\wqs^2)'(\ys)\left(3\lambda a_{1,0}A_{1,0}''+\frac{\lambda}{1-\lambda}
A_{1,0}\right)+ \sigma^2 O(\wqs(\ys)).
\end{align*}
 
 Similarly, we obtain
\begin{align*}
 S_{BBM}(B_{1,0}\wqs'(\ys)) & =\lambda \mu_\sigma \wqs'(\ys) \beta(\ys) B_{1,0}''' +
 \lambda \mu_\sigma \wqs''(\ys)(-B_{1,0}'') + \lambda \mu_\sigma \wqs'''(\ys)(-2  B_{1,0}')
 \\ &=\sigma\wqs(-\lambda B_{1,0}'')+\sigma\wqs'(-2\lambda B_{1,0}')
+\wqs^2(\ys)\left(\frac{\lambda}{1-\lambda}B_{1,0}''\right)
\\ &  
+(\wqs^2)'(\ys)\left(\frac{\lambda a_{1,0}}{2}B_{1,0}'''
+\frac{2\lambda}{1-\lambda}B_{1,0}'\right) + \sigma^2 O(\wqs(\ys)).
\end{align*}
Finally, we check that for $(k,l)$ such that $2\leq k+l\leq 3$,
\begin{align*}
& S_{BBM}( \sigma^l\wqs^k(\ys) A_{k,l}(y))=\sigma^l(\wqs^k)'(\ys)(-\lambda A_{k,l}'')
+ \sigma^2 O(\wqs(\ys)),\\
&  S_{BBM}( \sigma^l (\wqs^k)'(\ys)B_{k,l}(y) )=\sigma^2 O(\wqs(\ys)).
\end{align*}
\end{proof}
\begin{proof}[Proof of Lemma \ref{lem:Sint}.]
We have  
 \begin{align*}
\pd_x(w^2)&= \pd_x(A_{1,0}(y)^2\wqs(\ys)^2
+2A_{1,0}(y)B_{1,0}(y)\wqs(\ys)\wqs'(\ys))+\sigma^2 O(\wqs(\ys))
\\&= (1-\beta(\ys))\{(A_{1,0}^2)'\wqs^2(\ys)+(A_{1,0}B_{1,0})'(\wqs^2)'(\ys)\}
\\ & +A_{1,0}^2(\wqs^2)'(\ys)+2A_{1,0}B_{1,0}(\wqs\wqs')'(\ys)+\sigma^2 O(\wqs(\ys)).
 \end{align*}
By the definition of $\beta$ and \eqref{eq:qsigma}, we obtain
$$\pd_x(w^2)= \wqs^2(\ys)(A_{1,0}^2)'
+(\wqs^2)'(\ys)((A_{1,0}B_{1,0})'+A_{1,0}^2)+\sigma^2 O(\wqs(\ys)).$$
Next, by similar arguments,
\begin{align*}
  2\pd_x(\widetilde{Q}_\sigma(\ys) Q)&=2(1-\beta(\ys))\wqs(\ys)Q'+2\wqs'(\ys)Q
\\&= \wqs(\ys)(2Q')+\wqs'(\ys)(2Q)+\wqs^2(\ys)(-2a_{1,0}Q')+\sigma^2 O(\wqs(\ys)).
\end{align*}
Finally,  
\begin{align*}
 2\pd_x(\wqs(\ys)w)=\wqs^2(\ys)(2A_{1,0}')+(\wqs^2)'(\ys)(2A_{1,0}+B_{1,0}')
+\sigma^2 O(\wqs(\ys)).
\end{align*}
\end{proof}

\section{Appendix -- Proof of Lemma \ref{lem:z}}

The symmetry property $z(t,x)=z(-t,-x)$ is clear from \eqref{eq:df1}--\eqref{eq:df3} since the transformation $x\to -x$, $t\to -t$ gives $\ys\to -\ys$ (by parity of $\alpha$) and $y\to -y$, and since the functions  $A_{k,l}$, $\wqs^k$ are
even and the functions  $B_{k,l}$, $(\wqs^k)'$ are odd.

\medskip

Proof of \eqref{eq:z2}.  Note that from Proposition \ref{prop:decomp}, and the choice of $a_{k,l}$, $A_{k,l}$, $B_{k,l}$ for $(k,l)\in \Sigma_0$, solving
$(\Omega_{k,l})$, we have $S(z)=\mathcal{E}(t,x)$. Moreover, from the proof of Proposition~\ref{prop:decomp} (see Appendix B), the rest term $\mathcal{E}(t,x)$
is a finite sum of terms of the type $\sigma^l\wqs^k(\ys)f(y)$ or $\sigma^l(\wqs^k)'(\ys)f(y)$, where $k+l\geq 4$ and  $f$ is a bounded function such that $f'\in \mathcal{Y}$.
It follows that
$$\|S(t)\|_{H^1} \leq C \sigma^3 \|\wqs\|_{H^1}\leq C \sigma^{\frac {15}4}.$$

Proof of \eqref{eq:z1}.
We begin with some preliminary estimates.

\begin{claim}\label{cl:30}
\begin{equation}\label{eq:al}
\|\alpha\|_{L^\infty}\leq C\sqrt{\sigma},\quad \|\alpha'\|_{L^\infty}\leq C \sigma.
\end{equation}
For $t=\tau_\sigma$, for   $f\in \mathcal{Y}$, we have, for $ \sigma>0$ small,
\begin{equation}\label{eq:de}
\|f(y)\wqs(\ys)\|_{H^1}\leq C \sigma^{10},
\end{equation}
\begin{equation}\label{eq:Qy}
\|Q(y)-Q(x-\tfrac 12 \delta)\|_{H^1}\leq C \sigma^{10}.
\end{equation}
\end{claim}
\begin{proof}[Proof of Claim \ref{cl:30}.]
By the definition of $\wqs$ (see Claim \ref{cl:2.1}), we have
\begin{equation}\label{eq:dw}
\forall x\in \R,\quad 0\leq \wqs(x)\leq C\sigma e^{-\sqrt{\sigma} |x|}.
\end{equation}
It follows that $\int \wqs \leq C \sqrt{\sigma}$ and the result for $\alpha$ follows.
Since $\|\wqs\|_{L^\infty}\leq C \sigma$, the result for $\alpha'$ is also clear.

Now, we prove \eqref{eq:de}. Let $f\in \mathcal{Y}$, so that $|f(y)|\leq C |y|^r e^{-|y|}$ on $\R$. Note that for $t=\tau_\sigma$, since $\mu_\sigma>\frac 12$, we have
$
\sqrt{\sigma} |y_{\sigma}| \geq \sqrt{\sigma} ( \mu_\sigma \tau_\sigma
-|y| - |\alpha(\ys)|) \geq \frac 12 \sigma^{-\frac {1} {100}} - \sqrt{\sigma} |y| -1.
$
Thus, by \eqref{eq:dw},
$$
|\wqs(\ys) f(y)|^2 \leq C \sigma e^{-\sigma^{-\frac 1{100}}} |y|^{2r} e^{-2 (1-\sqrt{\sigma}) |y|}
\leq C e^{-\sigma^{-\frac 1{100}}}  e^{- |y|}.
$$
Since $\int e^{-|y|} dx = \int e^{-|y|} \frac {dy}{1-\alpha'(\ys)} \leq C.$
We obtain
$$
\|\wqs(\ys)f(y)\|_{L^2} \leq C e^{-\frac 12 \sigma^{-\frac 1{100}}} \leq C \sigma^{10}.
$$
The result for the $H^1$ norm is obtained in the same way.

Finally, we prove \eqref{eq:Qy}. First, we remark that for $t=\tau$ and $x>-\tfrac 12 \tau_\sigma$,  we have 
$|\alpha(\ys)-\tfrac 12 \delta|\leq K \sigma^{10}$.
Indeed, $|\alpha(\ys)-\tfrac 12 \delta|\leq C \int_{\ys}^{+\infty} \wqs
\leq C \sqrt{\sigma } e^{-\sqrt{\sigma} \ys}.$ For $t=\tau_\sigma$ and $x>-\tfrac 12 \tau_\sigma$, we have $\ys\geq \frac 14 \tau_\sigma$ and so
$e^{-\sqrt{\sigma} \ys}\leq e^{-\frac 14 \sigma^{-\frac 1{100}}}\leq C \sigma^{10}$.
Using this remark, we obtain directly for $t=\tau_{\sigma}$,
$$
\|Q(y)-Q(.-\tfrac 12 \delta)\|_{H^1(x> - \frac 12 \tau_\sigma)}\leq C \sigma^{10}.
$$
To complete the proof of \eqref{eq:Qy}, it suffices to use the decay of $Q$.
Note that if $x<-\frac 12 \tau_\sigma$, since $|\alpha(\ys)|\leq 1$,
we have $y<-\frac 12 \tau_\sigma +1$ and thus 
$$
\|Q(y)-Q(.-\tfrac 12 \delta)\|_{H^1(x< - \frac 12 \tau_\sigma)}\leq 
\|Q(y)\|_{H^1(x<- \frac 12 \tau_\sigma+1)}
+\|Q(.-\tfrac 12 \delta)\|_{H^1(x< - \frac 12 \tau_\sigma)}\leq C \sigma^{10}.
$$
\end{proof}

Now, we continue the proof of \eqref{eq:z1}.
First, from the expression of $z(\tau_\sigma)$,
the structure of the functions $A_{k,l}$, $B_{k,l}$, \eqref{eq:de} and
$\lim_{-\infty} \varphi=-1$, we have
(for simplicity, we drop the variable   $\ys$):
\begin{equation}\begin{split}\label{eq:c1}
\| z(\tau_\sigma) &-\{Q(y)+ 
\wqs - b_{1,0} \wqs'  + \gamma_{2,0} \wqs^2
-b_{2,0} (\wqs^2)' + \gamma_{1,1} \sigma \wqs
-b_{1,1} \sigma \wqs'\\ & +\gamma_{3,0} \wqs^3 + \gamma_{2,1} \sigma \wqs^2 
+\gamma_{1,2} \sigma^2 \wqs\}\|_{H^1}
\leq K\sigma^{\frac {13}4}.
\end{split}\end{equation}
Note that  $\sigma^{\frac {13}4}$ corresponds to the size of 
$\sigma^l(\wqs^k)'(\ys)B_{k,l}(y)$, for $k+l=3$, where $B_{k,l}$ is bounded
(see Lemma \ref{lem:2.5}).

It is natural to combine the following two terms
$$\wqs(\ys) - b_{1,0} \wqs'(\ys) \sim
\wqs(\ys - b_{1,0}),$$
but in fact, most terms above can be viewed as translation terms. To see this, 
let us now expand $\wqs(\ys-b_{1,0}-\sigma \tilde b_{1,1})$
and $(\wqs^2)'(\ys-b_{1,0}-\sigma \tilde b_{1,1})$
up to the order $\sigma^{\frac {13}4}$ in $H^1$:
\begin{equation}\label{eq:c2}\begin{split}
\|\wqs(\ys-b_{1,0}-\sigma \tilde b_{1,1})
& -\{ \wqs  - b_{1,0} \wqs'-  \tilde b_{1,1}\sigma\wqs'
+ \tfrac 12 b_{1,0}^2 \wqs'' +b_{1,0} \tilde b_{1,1} \sigma \wqs''
 \\ &-\tfrac 16 b_{1,0}^3 \wqs'''
+ \tfrac 1{24} b_{1,0}^4 \wqs''''\}\|_{H^1}
\leq K \sigma^{\frac {13}4},
\end{split}\end{equation}
\begin{equation}\label{eq:c2bis}
\|(\wqs^2)'(\ys-b_{1,0}-\sigma \tilde b_{1,1})
- \{(\wqs^2)' - b_{1,0}(\wqs^2)''\}\|_{H^1}\leq 
 C \sigma^{\frac {13}4}.
\end{equation}
In \eqref{eq:c2}, \eqref{eq:c2bis}, we now replace (from \eqref{eq:qsigma}, \eqref{eq:musigma}):
$$
\wqs''=\sigma \wqs - \frac 1{1-\lambda} \wqs^2  +\frac {\lambda}{1-\lambda} \sigma\wqs^2,
\quad 
\wqs'''=\sigma \wqs' - \frac 1{1-\lambda} (\wqs^2)' +\sigma^{\frac 52} O(\wqs),
$$
\begin{align*}
\wqs''''
&=\sigma \wqs'' - \frac {1}{1-\lambda} 2(\wqs'' \wqs + (\wqs')^2)
+\sigma^{\frac 52} O(\wqs)\\
&=\sigma^2 \wqs - \frac 5{1-\lambda} \sigma \wqs^2 + \frac {10}{3(1-\lambda)^2}
\wqs^3 +\sigma^{\frac 52} O(\wqs),
\end{align*}
$$
(\wqs^2)'' = 2 (\wqs'' \wqs + (\wqs')^2)
=4 \sigma \wqs^2 -\frac {10}{3(1-\lambda)} \wqs^3+\sigma^{\frac 52} O(\wqs).
$$
We obtain
\begin{equation}\label{eq:c3}\begin{split}
&  \|\{\wqs(\ys-b_{1,0}-\sigma \tilde b_{1,1}) - d(\lambda)(\wqs^2)'(\ys - b_{1,0}-\sigma \tilde b_{1,1})\}\\&
-\{\wqs-b_{1,0}\wqs'+\tfrac 12 b_{1,0}^2 \sigma \wqs-(\tilde b_{1,1} + \tfrac 16 b_{1,0}^3) \sigma\wqs' 
  - \tfrac 1{2(1-\lambda)}  b_{1,0}^2 \wqs^2
-b_{2,0} (\wqs^2)' \\
& -(-\tfrac \lambda {2(1-\lambda)} b_{1,0}^2+b_{1,0}(4d(\lambda)+\tfrac 1{(1-\lambda)} \tilde b_{1,1}) + \tfrac 5{24(1-\lambda)} b_{1,0}^4) \sigma \wqs^2 \\
& + (b_{1,0}\tilde b_{1,1}+\tfrac 1{24} b_{1,0}^4) \sigma^2 \wqs + (\tfrac 5{36(1-\lambda)^2}b_{1,0}^4+
\tfrac {10}{3(1-\lambda)} b_{1,0} d(\lambda)) \wqs^3 \}\|_{H^1}\leq C \sigma^{\frac {13}4}.
\end{split}
\end{equation}
Combining \eqref{eq:c1} and \eqref{eq:c3}, we find
\begin{equation*} \begin{split}
&  \|z(\tau_\sigma) - \{Q(y)+\wqs(\ys-b_{1,0}-\sigma \tilde b_{1,1}) - d(\lambda)(\wqs^2)'(\ys - b_{1,0}-\sigma \tilde b_{1,1})\}\\&
+(\gamma_{1,1}-\tfrac 12 b_{1,0}^2) \sigma \wqs
+(-b_{1,1}+\tilde b_{1,1} + \tfrac 16 b_{1,0}^3) \sigma\wqs' 
 +(\gamma_{2,0}+  \tfrac 1{2(1-\lambda)}  b_{1,0}^2) \wqs^2
\\
& +(\gamma_{2,1}-\tfrac \lambda {2(1-\lambda)} b_{1,0}^2+b_{1,0}(4d(\lambda)+\tfrac 1{(1-\lambda)} \tilde b_{1,1}) + \tfrac 5{24(1-\lambda)} b_{1,0}^4) \sigma \wqs^2 \\
& +(\gamma_{1,2}- b_{1,0}\tilde b_{1,1}-\tfrac 1{24} b_{1,0}^4)\sigma^2 \wqs 
+ (\gamma_{3,0}-\tfrac 5{36(1-\lambda)^2}b_{1,0}^4-
\tfrac {10}{3(1-\lambda)} b_{1,0} d(\lambda)) \wqs^3 \}\|_{H^1}\leq C \sigma^{\frac {13}4}.
\end{split}
\end{equation*}
It follows that with the choice
\begin{align*}
&\gamma_{1,1}=\frac 12 b_{1,0}^2,\quad \tilde b_{1,1}=b_{1,1} -\frac 16 b_{1,0}^3,
\quad 
\gamma_{1,2}=b_{1,0}b_{1,1}- \frac 18 b_{1,0}^4,
\\
&\gamma_{2,0}=- \frac 1{2(1-\lambda)} b_{1,0}^2,\quad
\gamma_{3,0}= \frac 5{36(1-\lambda)^2} b_{1,0}^4+ \frac {10}{3(1-\lambda)} b_{1,0}d(\lambda) ,\\
&\gamma_{2,1}= \frac {\lambda}{2(1-\lambda)} b_{1,0}^2-4 b_{1,0}d(\lambda) 
- \frac 1{1-\lambda} b_{1,0}b_{1,1} - \frac 1{24(1-\lambda)} b_{1,0}^4,
\end{align*}
we obtain
\begin{equation}\label{eq:c9}
\|z(\tau_\sigma)-\{Q(y)+\wqs(\ys - b_{1,0}-\sigma \tilde b_{1,1}) - d(\lambda) (\wqs^2)'(\ys - b_{1,0}-\sigma \tilde b_{1,1})\}\|_{H^1}
\leq C \sigma^{\frac {13} 4}.
\end{equation}
Together with \eqref{eq:Qy}, this completes the proof of \eqref{eq:z1} 
This justifies in particular
the choices of $\gamma_{2,0}$, $\gamma_{1,1}$, $\gamma_{3,0}$,
$\gamma_{2,1}$ and $\gamma_{1,2}$ done in Lemmas \ref{lem:b20}, \ref{lem:b11}
 and \ref{lem:2.5}.

\section{Appendix -- Monotonicity properties}\label{se:xC}

In this Appendix, we prove Lemmas \ref{le:mo} and \ref{le:qd}.
\begin{proof}[Proof of Lemma \ref{le:mo}]
First, note that as in \eqref{eq:da}, we have $|a_2-c_2|\leq C (|\bar c_2(T)-c_2|+|\bar c_2(T_0)-c_2|)$, and then $a_2-1 \geq \frac 12 (c_2-1)>0$.

Now, following computations in \cite{DiMa}, p. 424 (see also \cite{Mi2} and \cite{Di}), we compute $\frac d{dt}\mathcal{G}(t)$.
Recall that $\psi$ and $m(t)$ are defined in \eqref{eq:ph} with $\kappa=\sqrt{\frac {c_1+3}{c_1-1}}$.

By \eqref{eq:vq} and conservation of $N(u(t))$ and $E(u(t))$,  we have
$$
\frac d{dt}\mathcal{G}(t)       = -\frac d{dt} (a_2 \mathcal{N}_1(t) - \mathcal{E}_1(t)).
$$
\begin{align*}
 & \frac d{dt} (a_2 \mathcal{N}_1(t) - \mathcal{E}_1(t))    = 
a_2 \int (u_{xt} u_x + u_t u) \psi(x-m(t))  
-\int (u_t u + u_t u^2) \psi(x-m(t))   \\
& - \tfrac 12 a_2 m'(t) \int (u_x^2+u^2) \psi'(x-m(t))  
+\tfrac 12 m'(t) \int (u^2 +\tfrac 23 u^3) \psi'(x-m(t))\\
& = a_2 \int ((1-\partial_x^2) u_t) u \psi(x-m(t)) 
- a_2 \int u_{xt} u \psi'(x-m(t)) - \int u_t (u+u^2) \psi(x-m(t))
  \\
& - \tfrac 12 a_2 m'(t) \int (u_x^2+u^2) \psi'(x-m(t))  
+\tfrac 12 m'(t) \int (u^2 +\tfrac 23 u^3) \psi'(x-m(t)).
\end{align*}
Now, we use $(1-\partial_x^2) u_t = - \partial_x (u+u^2)$ and   $u_t=-\partial_x h$, where
$h=(1-\partial_x^2)^{-1} (u+u^2)$.
We get
\begin{align*}
&  \frac d{dt} (a_2 \mathcal{N}_1(t) - \mathcal{E}_1(t))     
= -a_2 \int (u_x+2 u_x u ) u \psi(x-m(t)) 
+ a_2 \int h_{xx} u \psi'(x-m(t)) 
\\ & +\int h_x (h-h_{xx}) \psi(x-m(t))
 - \tfrac 12 a_2 m'(t) \int (u_x^2+u^2) \psi'(x-m(t))  
\\& +\tfrac 12 m'(t) \int (u^2 +\tfrac 23 u^3) \psi'(x-m(t)) \\
& = a_2 \int (\tfrac 12 u^2 + \tfrac 23 u^3) \psi'(x-m(t)) 
+a_2 \int (h-u-u^2) u \psi'(x-m(t))\\
&+ \int (-\tfrac 12 h^2 +\tfrac 12 h_x^2)  \psi'(x-m(t))
- \tfrac 12 a_2 m'(t) \int (u_x^2+u^2) \psi'(x-m(t))  
\\& +\tfrac 12 m'(t) \int (u^2 +\tfrac 23 u^3) \psi'(x-m(t)) \\
& = -\frac 12 (a_2+(a_2-1) m'(t)) \int u^2 \psi'(x-m(t))
- \frac 12 a_2 m'(t) \int u_x^2\psi'(x-m(t))\\
& + \frac 13 (m'(t)-a_2) \int u^3 \psi'(x-m(t))+a_2 \int h u \psi'(x-m(t)) + \int (-\tfrac 12 h^2 +\tfrac 12 h_x^2)  \psi'(x-m(t)).
\end{align*}
Using the estimate
$$
2 a_2 \left|\int hu \psi'(x-m(t)) \right|\leq a_2^2 \int u^2 \psi'(x-m(t))
+ \int h^2 \psi'(x-m(t)),
$$
we find
\begin{align*}
&  2 \frac d{dt} (a_2 \mathcal{N}_1(t) - \mathcal{E}_1(t)) 
     \leq -(a_2-1) (m'(t)-a_2) \int u^2 \psi'(x-m(t))
\\ &- a_2 m'(t) \int u_x^2 \psi'(x-m(t)) + \frac 23 (m'(t)-a_2) \int u^3 \psi'(x-m(t))+  \int  h_x^2  \psi'(x-m(t)).
\end{align*}
By the definition of $m(t)$ and \eqref{eq:fi}, we have, for $0<c_2-1<\epsilon_0$ small enough,
$$m'(t)\geq \frac 14 (c_1+c_2+2)\geq \frac 14(c_1+3),\qquad
\frac 14 (c_1-1)\leq  m'(t)-a_2\leq c_1-1
.$$
Moreover, $a_2-1\geq \frac 12 (c_2-1)$.
Thus,
\begin{align*}
&  2 \frac d{dt} (a_2 \mathcal{N}_1(t) - \mathcal{E}_1(t))   
  \leq -\frac 18 (c_2-1) (c_1-1)  \int u^2 \psi'(x-m(t))
- \frac 14 (c_1+3) \int u_x^2 \psi'(x-m(t))\\
& + \frac 23 (c_1-1)\int |u|^3 \psi'(x-m(t))+  \int  h_x^2  \psi'(x-m(t)).
\end{align*}

Now, we claim the following
\begin{equation}\label{eq:uh}
\int u_x^2 \psi'(x-m(t))
\geq \left(1-\frac 1{\kappa^2}\right) \int h_x^2 \psi'(x-m(t)) - 4 \int u_x^2(u+2u^2)\psi'(x-m(t)).
\end{equation}
Indeed, note that using $h-h_{xx}=u+u^2$, we have
\begin{align*}
& \int u_x^2 \psi'(x-m(t))  = \int (-u^2 + h - h_{xx})_x^2 \psi'(x-m(t))\\
 & =\int [4u_x^2u^2 + h^2_x + h_{xxx}^2 -2 h_x h_{xxx} -4 u_x u(h-h_{xx})_x  ] \psi'(x-m(t))\\
& =\int [h_x^2+h_{xxx}^2- 2h_x h_{xxx} + 4 u^2 u_x^2  - 4 u_x u (u+u^2)_x  ]\psi'(x-m(t)) 
\\
& =\int [h_x^2+h_{xxx}^2+ 2h_{xx}^2+ 4 u^2 u_x^2  -4 u_x u (u_x+2u_xu)  ]\psi'(x-m(t))
- \int h_{x}^2\psi'''(x-m(t))\\
& \geq \int [h_x^2 - 4 u_x^2 (u+2u^2)  ]\psi'(x-m(t))- \int h_{x}^2\psi'''(x-m(t)).
\end{align*}
Thus, by \eqref{eq:ph},   \eqref{eq:uh} is proved.

Therefore, we obtain
\begin{align*}
&\quad 2 \frac d{dt} (a_2 \mathcal{N}_1(t) - \mathcal{E}_1(t))   
\\ & \leq - \frac 18 (c_2-1) (c_1-1)  \int u^2 \psi'(x-m(t))
-\left[\frac {c_1+7} 8 \left(1-\frac 1{\kappa^2}\right) -1\right] \int h_x^2 \psi'(x-m(t))
\\ &
- \frac {c_1-1}8 \int u_x^2 \psi'(x{-}m(t))
+ \frac {c_1+7}2 \int u_x^2 (u+2u^2)\psi'(x{-}m(t))+ (c_1-1) \int |u|^3 \psi'(x{-}m(t)).
\end{align*}
Note that for our choice of $\kappa,$ we have $\frac {c_1+7} 8 \left(1-\frac 1{\kappa^2}\right) -1=0$.

Now, we treat the two nonlinear terms 
$\int u_x^2 (u+2u^2)\psi'(x-m(t))$ and $\int u^3 \psi'(x-m(t))$.
Let $I=[\rho_2(t)+\frac 18 T , \rho_1(t) - \frac 18 T]$.

First, for $x\in I$, by \eqref{eq:tt} and the decay of $\varphi_c$, we have,
for $t\geq T$,
$$|u(t,x)|\leq |R_1(t,x)|+|R_2(t,x)|+ |\eta(t,x)|
\leq C (c_2-1)^{2}.$$
Thus, for $0<c_2-1\leq \epsilon_0$ small enough, we obtain
\begin{equation*}\begin{split}
\left|4c_1 \int_I u_x^2 (u+2 u^2) \psi'(x-m(t)) \right|
&\leq 8c_1 (\|u\|_{L^\infty(I)}+\|u\|_{L^\infty(I)}^2) 
\int u_x^2 \psi'(x-m(t))\\ &
 \leq \frac 1{16} (c_1-1) \int u_x^2 \psi'(x-m(t)).
\end{split}\end{equation*}
\begin{equation*}\begin{split}
\left| (c_1-1) \int_I u^3 \psi'(x-m(t)) \right|
& \leq (c_1-1) \|u\|_{L^\infty(I)} \int u^2 \psi'(x-m(t))\\
& \leq \frac {(c_1-1)(c_2-1)} {16} \int u^2 \psi'(x-m(t)).
\end{split}\end{equation*}

Next, for $x\in \R\setminus I$, we have
$|x-m(t)|\geq \frac14c_1t+\frac 18(c_1-c_2-1)T$ and 
$\psi'(x-m(t)) \leq C e^{- \gamma t}$, where $\gamma=\gamma(c_1)>0$,
and so
$$
2 \frac d{dt} (a_2 \mathcal{N}_1(t) - \mathcal{E}_1(t))   
\leq C e^{- \gamma t}.
$$
By integration on $[T,T_0]$, we get
$$
\mathcal{G}(T_0) - \mathcal{G}_1(T)
\geq - C e^{-\gamma T}
\geq - C(c_2-1)^{10}.
$$
Thus, the lemma is proved.
\end{proof}

\begin{proof}[Proof of Lemma \ref{le:qd}]
We expand $u(t,x)=R_1(t,x)+R_2(t,x)+\eta(t,x)$ in the expression of $\mathcal{G}(t)$,
using the following estimates, for all $t\geq T$,
$$
\int R_1(t )(1-\psi(x-m(t)))  +
\int R_1(t )R_2(t )  + \int R_2(t ) (\psi'+\psi)(x-m(t))  \leq C (c_2-1)^{10}.
$$
We obtain
\begin{equation*}\begin{split}
        \mathcal{G}(t)&=a_2 N(R_2(t)) - E (R_2(t)) 
        +a_2 \int (\eta_x(t) R_{2x}(t) + \eta(t) R_2(t))\\
        & - \int (\eta(t) R_2(t) + \eta (t) R_2^2(t)) 
        + \mathcal{H}(t) + O(\|\eta(t)\|_{H^1}^3 + O ((c_2-1)^{10}).
\end{split}\end{equation*}
But, by the equation of $\varphi_{\bar c_2(t)}$ and \eqref{eq:o1},
$$
\int (\eta_x  R_{2x}  + \eta  R_2 )=\int \eta (1-\partial_x^2) R_2=0,$$ 
$$ \int (\eta(t)  R_2(t) +  \eta(t)   R_2^2(t) ) 
         =-\bar c_2(t) \int \eta (1-\partial_x^2) R_2=0,
$$
 and by the definition of $a_2$,
$$
a_2 N(R_2(T_0)) - E(R_2(T_0) = 
a_2 N(R_2(T))-E(R_2(T)).
$$
Thus, we get
$$
\mathcal{G}(T_0)        -\mathcal{G}(T) = \frac 12 (\mathcal{H}(T_0)-\mathcal{H}(T))
+O(\|\eta\|_{H^1}^3)+O((c_2-1)^{10}).
$$

Finally, the coercivity property
$$\mathcal{H}(t)
\geq  \sigma_0 \int \left[(c_2-1) \eta^2 + \eta_x^2 \right](t,x) (1-\psi(x-m(t))) dx,
$$
under the orthogonality conditions \eqref{eq:o1} is a standard fact, 
see for example Appendix A in \cite{DiMa}.
Note that the proof is based on 
the following positivity property (see \cite{We3}): there exists $\sigma>0$ such that
\begin{equation*}\begin{split}
& \int v (1-\partial_x^2) \varphi_c = \int v (1-\partial_x^2) \partial_x \varphi_c=0\\
&\Rightarrow \quad 
\int c v_x^2 + (c-1) v^2 - 2 \varphi_c^2 v^2 \geq \sigma \|v\|_{H^1_c}^2.
\end{split}\end{equation*}
\end{proof}

\end{document}